# Linear Representation of Lie Group

Aleks Kleyn

Aleks_Kleyn@MailAPS.org
http://AleksKleyn.dyndns-home.com:4080/
http://sites.google.com/site/AleksKleyn/
http://arxiv.org/a/kleyn_a_1
http://AleksKleyn.blogspot.com/

ABSTRACT. In this book, I explored differential equations for operation in Lie group and for representations of group Lie in a vector space.

# Contents

CHAPTER 1

# Introduction

## 1.1. About This Book

This book is based on the research project that I performed during my study in the Odessa university. I was very lucky that my teacher Gavrilchenko Michail Leonidovich provided lectures "Lie Group and Lie Algebra". These lectures were based on the book [5]. These lectures became the basis of my research for many years. Based on these lectures, I wrote down the equations of a Lie group for the particular case of the representation of the group $GL$ in the vector space of geometric objects.

Subsequently, I replaced the group $GL$ by an arbitrary Lie group. This is the form that I used to present theory in this book. Accordingly, chapters of the book have the following contents. Representation of a Lie group is closely connected with the theory of differential equations. In the chapter 2, I give a brief lecture of the theory of differential equations. I dedicated this lecture to completely integrable systems of differential equations. In the chapter 3, I write down the differential equations for a Lie group. I write down differential equations of the representation of Lie group in the chapter 4.

## 1.2. Conventions

CONVENTION 1.2.1. *Let $A$ be $\Omega_1$-algebra. Let $B$ be $\Omega_2$-algebra. Notation*

$$A \xrightarrow{\ *\ } B$$

*means that there is representation of $\Omega_1$-algebra $A$ in $\Omega_2$-algebra $B$.* □

Without a doubt, the reader may have questions, comments, objections. I will appreciate any response.

CHAPTER 2

# Differential Equations

## 2.1. Completely Integrable Systems

CONVENTION 2.1.1. *Let $F$ be complete field. In this book, we consider $F$-vector spaces. Correspondingly, we consider differential equations over field $F$.* □

Consider[2.1] a system of $nm$ partial differential equations

$$(2.1.1) \qquad \frac{\partial \theta^\alpha}{\partial x^i} = \psi_i^\alpha(\theta^1, ..., \theta^m, x^1, ..., x^n) = \psi_i^\alpha(\theta, x)$$

$$1 \leq \alpha \leq m \quad 1 \leq i \leq n \quad \theta^\alpha = \theta^\alpha(x^1, ..., x^n)$$

This system is solved relative all partial derivatives and is equivalent to system of total differential equations

$$d\theta^\alpha = \psi_i^\alpha(\theta, x) dx^i$$

To find conditions of integrability we differentiate (2.1.1) with respect to $x^j$.

$$\frac{\partial^2 \theta^\alpha}{\partial x^j \partial x^i} = \frac{\partial \psi_i^\alpha}{\partial x^j} + \frac{\partial \psi_i^\alpha}{\partial \theta^\sigma} \frac{\partial \theta^\sigma}{\partial x^j} =$$

$$= \frac{\partial \psi_i^\alpha}{\partial x^j} + \frac{\partial \psi_i^\alpha}{\partial \theta^\sigma} \psi_j^\sigma$$

$$\frac{\partial^2 \theta^\alpha}{\partial x^i \partial x^j} = \frac{\partial \psi_j^\alpha}{\partial x^i} + \frac{\partial \psi_j^\alpha}{\partial \theta^\sigma} \psi_i^\sigma$$

Because second derivative of continuous function does not depend on order of differentiation we expect that

$$(2.1.2) \qquad \frac{\partial \psi_i^\alpha}{\partial x^j} + \frac{\partial \psi_i^\alpha}{\partial \theta^\sigma} \psi_j^\sigma = \frac{\partial \psi_j^\alpha}{\partial x^i} + \frac{\partial \psi_j^\alpha}{\partial \theta^\sigma} \psi_i^\sigma$$

DEFINITION 2.1.2. *We call the system of differential equations* (2.1.1) **completely integrable system** *if condition* (2.1.2) *is satisfied identically.* □

Solution of completely integrable system (2.1.1) has Taylor expansion

$$(2.1.3) \quad \theta^\alpha = C^\alpha + \left(\frac{\partial \theta^\alpha}{\partial x^i}\right)_{x=x_0} (x^i - x_0^i) + \frac{1}{2} \left(\frac{\partial^2 \theta^\alpha}{\partial x^j \partial x^i}\right)_{x=x_0} (x^i - x_0^i)(x^j - x_0^j) + ...$$

Here $C^\alpha$ ($1 \leq \alpha \leq m$) are constants,

$$\left(\frac{\partial \theta^\alpha}{\partial x^i}\right)_{x=x_0} = \psi_i^\alpha(C^1, ..., C^m, x_0^1, ..., x_0^n)$$

We obtain the rest coefficients by differentiating of right part of (2.1.1) and substituting $C^\alpha$ and $x_0^i$. We tell about solution when expansion (2.1.3) converges.

---

[2.1]I wrote this chapter under influence of [5].

According to theorem by Caushi and Covalevskaya it has place when $\psi_i^\alpha$ are analitical by all arguments.

For the vicinity of point $(x_0^1, ..., x_0^n)$ for which series (2.1.3) are convergent we have a solution

$$\theta^\alpha = \varphi^\alpha(x^1, ..., x^n, C^1, ..., C^m) \tag{2.1.4}$$

determined by $m$ constants.

If the condition (2.1.2) is not satisfied identically, then system of differential equations (2.1.1) is not completely integrable. Therefore functions $\theta^\alpha$ satisfy (2.1.1) and (2.1.2). If constraints (2.1.2) are independent from system $S_1$, then equation (2.1.2) puts constraints on functions $\theta^\alpha$ and generates the system $S_1$ of equations. We differentiate equations of the system $S_1$ with respect to $x^k$ and express derivatives $\dfrac{\partial \theta^\alpha}{\partial x^i}$ by their values from equations (2.1.1). If new equations are not satisfied identically and get new system $S_2$ of equations. We can repeat this process soon.

Because we cannot have more than $m$ constraints, this chain breaks. If system of differential equations (2.1.1) has solution this solution satisfies all constraints, otherwise these constraints are conflicting.

THEOREM 2.1.3. *Necessary and sufficient condition in order that system of differential equations* (2.1.1) *has solution is existing of such $N$ that $S_1, ..., S_N$ are algebraically compatible, and $S_{N+1}$ is consequence of them. If we get $p$ constraints then common solution depends on $m - p$ arbitrary constants.*

PROOF. We proved that this condition is necessary. Now we show that this condition is sufficient. Assume we got $p$ constraints

$$\Phi^\gamma(\theta, x) = 0 \quad \gamma = 1, ..., p \tag{2.1.5}$$

Because

$$Rg \left\| \frac{\partial \Phi^\gamma}{\partial \theta^\alpha} \right\| = p$$

we can express $p$ variables $\theta^\alpha$ through others using (2.1.5). If we enumerate variables by right way we get

$$\theta^\nu = \phi^\nu(\theta^{p+1}, ..., \theta^m, x) \quad \nu = 1, ..., p \tag{2.1.6}$$

We substitute (2.1.6) into (2.1.1) and get

$$\frac{\partial \theta^\alpha}{\partial x^i} = \bar{\psi}_i^\alpha(\theta^{p+1}, ..., \theta^m, x^1, ..., x^n) \quad \alpha = p+1, ..., m \tag{2.1.7}$$

The condition (2.1.2) is satisfied. System of differential equations (2.1.7) is completely integrable and has solution that depends on $m - p$ arbitrary constants. First $p$ equations are fulfilled because of (2.1.6) when $\theta^{p+1}, ..., \theta^m, x$ are solutions of (2.1.7). $\square$

## 2.2. Linear Differential Operator

Let us have $r$ differential operators

$$X_a f = \xi_a^i \frac{\partial f}{\partial x^i} \tag{2.2.1}$$

Operator (2.2.1) is linear

$$X_a(f + g) = X_a f + X_a g$$

$$X_a(\lambda f) = \lambda X_a f \quad \lambda = const$$

The product of these operators has form

$$X_a X_b f = X_a(X_b f) = \xi_a^i \frac{\partial}{\partial x^i}\left(\xi_b^j \frac{\partial f}{\partial x^j}\right)$$

$$= \xi_a^i \frac{\partial \xi_b^j}{\partial x^i}\frac{\partial f}{\partial x^j} + \xi_a^i \xi_b^j \frac{\partial^2 f}{\partial x^i \partial x^j} =$$

$$= \frac{\partial f}{\partial x^j} X_a \xi_b^j + \xi_a^i \xi_b^j \frac{\partial^2 f}{\partial x^i \partial x^j}$$

We used (2.2.1) on the last step. This is not operator of type (2.2.1). However because second term is symmetric relative $a$ and $b$ we can introduce commutator

(2.2.2) $$(X_a, X_b)f = X_a X_b f - X_b X_a f = (X_a \xi_b^j - X_b \xi_a^j)\frac{\partial f}{\partial x^j}$$

This is linear operation and we call it Poisson brackets. Skew symmetry of this commutator follows from (2.2.2).

Simple calculation shows that

(2.2.3) $$((X_a, X_b), X_c)f = ((X_a, X_b)\xi_c^j - X_c(X_a \xi_b^j - X_b \xi_a^j))\frac{\partial f}{\partial x^j} =$$

$$= \left(\frac{\partial \xi_c^j}{\partial x^k} X_a \xi_b^k - \frac{\partial \xi_c^j}{\partial x^k} X_b \xi_a^k - \frac{\partial \xi_b^j}{\partial x^k} X_c \xi_a^k + \xi_c^p \xi_a^q \frac{\partial^2 \xi_b^j}{\partial x^p \partial x^q} + \frac{\partial \xi_a^j}{\partial x^k} X_c \xi_b^k + \xi_c^p \xi_b^q \frac{\partial^2 \xi_a^j}{\partial x^p \partial x^q}\right)\frac{\partial f}{\partial x^j}$$

If we change order of parameters we finally get equation

(2.2.4) $$((X_a, X_b), X_c)f + ((X_b, X_c), X_a)f + ((X_c, X_a), X_b)f = 0$$

(2.2.4) is Jacoby condition.

THEOREM 2.2.1. *Commutator of linear combination of operators*

$$X'_a = \lambda_a^b X_b$$

*where $\lambda_a^b$ are functions of $x$ expresses through these operators and their commutators*

(2.2.5) $$(X'_a, X'_c)f = \mu_{ac}^d X_d f - \mu_{ca}^d X_d f + \lambda_a^b \lambda_c^d (X_b, X_d)f$$

PROOF. Commutator has form

(2.2.6) $$\begin{aligned}(X'_a, X'_c)f &= (\lambda_a^b X_b, \lambda_c^d X_d)f \\ &= \lambda_a^b X_b(\lambda_c^d X_d f) - \lambda_c^d X_d(\lambda_a^b X_b f) \\ &= \lambda_a^b X_b \lambda_c^d X_d f + \lambda_a^b \lambda_c^d X_b X_d f - \lambda_c^d X_d \lambda_a^b X_b f - \lambda_c^d \lambda_a^b X_d X_b f\end{aligned}$$

(2.2.5) follows from (2.2.6) if we assume

$$\mu_{ac}^d = \lambda_a^b X_b \lambda_c^d$$

□

## 2.3. Complete System of Linear Partial Differential Equations

Now we want to learn system of differential equations

(2.3.1) $$X_a f = 0 \quad X_a = \xi_a^i \frac{\partial}{\partial x^i} \quad 1 \leq i \leq n \quad 1 \leq a \leq r$$

We assume that

(2.3.2) $$\text{rank} \, \|\xi_a^i\| = r \leq n$$

This means that equations (2.3.1) are line independent.

If $r = n$, the only solution is evidently $f = \text{const}$. Statement that every solution of equations (2.3.1) is also solution of equation

(2.3.3) $$(X_a, X_b)f = 0$$

follows from (2.2.2). If

(2.3.4) $$(X_a, X_b)f = \gamma_{ab}^c X_c f$$

where the $\gamma$ are functions of $x$, then system of differential equations (2.3.1) and (2.3.3) is equivalent to system of differential equations (2.3.1). We do not know in advance whether (2.3.3) set new constraints for (2.3.1). We adjoin to system of differential equations (2.3.1) those equations of system (2.3.3) which do not hold (2.3.4). We obtain $s$ equations, $s \geq r$. We repeat the same procedure for new system and obtain $t$ equations, $t \geq s$, and so on.

As result of the sequence of indicated actions we may get $n$ equations. In this case the system of differential equations (2.3.1) has the only solution $f = \text{const}$, because

$$\frac{\partial f}{\partial x^i} = 0$$

Otherwise we get $u$ equations, $u \leq n$, which hold condition (2.3.4). In this case obtained system of differential equations is called **complete** of order $u$. Solution of the system of differential equations (2.3.1) is solution of complete system of differential equations.

THEOREM 2.3.1. *If system of differential equations* (2.3.1) *is a complete system, so also is*

(2.3.5) $$X'_a f = \lambda_a^b X_b f = 0$$

*where* $\lambda$ *are functions of* $x$

(2.3.6) $$\det \|\lambda_a^b\| \neq 0 \quad 1 \leq a, b \leq r$$

PROOF. From (2.3.4) it follows that the right-hand member of (2.2.5) is linear in $X_a f$. According to (2.3.5), (2.3.6) $X_a f$ is linear in operators $X'_a f$. Therefore, we have from (2.2.5) expressions of the form (2.3.4), which proves the statement. □

Since (2.3.2), there is no loss in generality in assuming that

(2.3.7) $$\det \|\xi_a^i\| \neq 0 \quad 1 \leq i, a \leq r$$

Hence we may solve system of differential equations (2.3.1) for $\frac{\partial f}{\partial x^1}$, ..., $\frac{\partial f}{\partial x^r}$ and write the result in the form

(2.3.8) $$X'_a f = \frac{\partial f}{\partial x^a} + \psi_a^t \frac{\partial f}{\partial x^t} = 0 \quad a = 1, ..., r \quad t = r+1, ..., n$$

These equations are of the form (2.3.5) and, hence, form system of differential equations equivalent to system (2.3.1).

Complete system of differential equations expressed in the form (2.3.8) is called **Jacobian**.

THEOREM 2.3.2. *Complete system of differential equations* (2.3.1) *admits exactly* $n - r$ *independent solutions.*

PROOF. Similarly to equation (2.3.4), we have

$$(2.3.9) \qquad (X'_a, X'_b)f = \gamma'^c_{ab} X'_c f$$

From (2.2.2) and (2.3.8), it follows that

$$(2.3.10) \qquad (X'_a, X'_b)f = (X_a \psi^t_b - X_b \psi^t_a)\frac{\partial f}{\partial x^t} \quad t = r+1, ..., n$$

Since $\frac{\partial f}{\partial x^1}$, ..., $\frac{\partial f}{\partial x^r}$ do not appear in (2.3.10), $\gamma'^c_{ab} = 0$. Hence, we get for complete system in the Jacobian form

$$(X'_a, X'_b)f = 0$$
$$X_a \psi^t_b - X_b \psi^t_a = 0$$

$$(2.3.11) \qquad \frac{\partial \psi^p_b}{\partial x^a} + \psi^q_a \frac{\partial \psi^p_b}{\partial x^q} = \frac{\partial \psi^p_a}{\partial x^b} + \psi^q_b \frac{\partial \psi^p_a}{\partial x^q} \quad a, b = 1, ..., r \quad p, q = r+1, ..., n$$

Comparing (2.3.11) and (2.1.2), we see that system of differential equations

$$(2.3.12) \qquad \frac{\partial x^p}{\partial x^a} = \psi^p_a \quad a = 1, ..., r \quad p = r+1, ..., n$$

is completely integrable. According to theorem 2.1.3 system of differential equations (2.3.8), and consequently, system of differential equations (2.3.1) admit $n - r$ independent solutions, and system of differential equations (2.3.1) does not admit more than $n - r$ independent solutions. □

## 2.4. Essential Parameters in a Set of Functions

DEFINITION 2.4.1. Consider $n$ functions

$$f^i(x^1, ..., x^n, a^1, ..., a^r), i = 1, ..., n$$

of variables $x^1, ..., x^n$ and parameters $a^1, ..., a^r$. We assume that functions $f^i$ are continuous in $x^1, ..., x^n$ and $a^1, ..., a^r$. Also their derivatives up to any order are continuous. Parameters $a^1, ..., a^r$ are called **essential parameters** unless there are functions $A^1, ..., A^{r-1}$ of $a^1, ..., a^r$ such that we have identically

$$(2.4.1) \qquad f^i(x, a) = F^i(x, A)$$

□

THEOREM 2.4.2. *A necessary and sufficient condition that the* $r$ *parameters* $a^\alpha$ *be essential is that the functions* $f^i$ *do not satisfy an equation of the form*

$$(2.4.2) \qquad \phi^\alpha \frac{\partial f^i}{\partial a^\alpha} = 0$$

*where* $\phi^\alpha \neq 0$.

PROOF. Assume that parameters $a$ are not essential. Then there exist $A$ such that (2.4.1) is satisfied and

$$\operatorname{rank} \left\| \frac{\partial A^\sigma}{\partial a^\alpha} \right\| < r$$

Therefore there exist functions $\phi^\alpha(a) \neq 0$ such that

$$\phi^\alpha \frac{\partial A^\sigma}{\partial a^\alpha} = 0$$

Hence functions $A^1$, ..., $A^{r-1}$ and any their function $\Phi(A)$ also satisfy (2.4.2). Functions $F^i$ in (2.4.1) are example of function $\Phi$. Therefore if $f^i$ do not essentially depend on $a^\alpha$ they satisfy system (2.4.2).

Conversely, suppose that functions $f^i$ satisfy equation (2.4.1) for some functions $\phi^\alpha$. This equation admits $r-1$ independent solutions $A^1$, ..., $A^{r-1}$ which are functions of $a^\alpha$ alone and any solution of (2.4.2) is function of $A^1$, ..., $A^{r-1}$. Therefore, each of functions $f^i$ is a function of $x^1, ..., x^n$ and $A^1$, ..., $A^{r-1}$ and $a^1$, ..., $a^r$ are not essential. $\square$

REMARK 2.4.3. From theorem 2.4.2, it follows that, if the functions $f^i$ satisfy a complete system of $s$ differential equations

$$(2.4.3) \qquad \phi^\alpha_\sigma \frac{\partial f^i}{\partial a^\alpha} = 0 \quad \alpha = 1, ..., r \quad \sigma = 1, ..., s < r$$

then $f^i$ are functions of $x$ and of $r-s$ independent solutions of system (2.4.3), which are functions of the $a$, and consequently $f^i$ are expressible in terms of $r-s$ essential parameters. $\square$

We may interpret the system of equations (2.4.2) as system of linear equations for $\phi^\alpha$. From linear algebra, it follows, that solutions of system (2.4.3) form a vector space, and only in case of condition

$$(2.4.4) \qquad \mu_0 = \operatorname{rank}\left(\frac{\partial f^i}{\partial a^\alpha}\right) = r$$

the system of linear equations (2.4.3) has the only solution $\phi^\alpha = 0$. Thus, from remark 2.4.3, it follows, that parameters are essential, if condition (2.4.4) is satisfyed.

Suppose that $\mu_0 < r$. Differentiating the system of differential equations (2.4.2) with respect to the $x^j$, we have

$$(2.4.5) \qquad \phi^\alpha \frac{\partial f^i}{\partial a^\alpha} = 0 \quad \phi^\alpha \frac{\partial f^i_{,j}}{\partial a^\alpha} = 0 \quad f^i_{,j} = \frac{\partial f^i}{\partial x^j}$$

Let

$$(2.4.6) \qquad \mu_1 = \operatorname{rank}\left(\frac{\partial f^i}{\partial a^\alpha} \quad \frac{\partial f^i_{,j}}{\partial a^\alpha}\right)$$

Evidently $\mu_1 \geq \mu_0$. If $\mu_1 = r$, the system of linear equations (2.4.5) admits the only solution $\phi^\alpha = 0$. Consequently the parameters are essential.

If we put

$$f^i_{,j_1...j_s} = \frac{\partial^s f^i}{\partial x^{j_1}...\partial x^{j_s}}$$

we have from (2.4.2) by repeated differentiation

$$(2.4.7) \qquad \phi^\alpha \frac{\partial f^i_{,j_1...j_s}}{\partial a^\alpha} = 0$$

We denote

$$\mu_s = \text{rank}\left(\frac{\partial f^i}{\partial a^\alpha} \quad \frac{\partial f^i_{,j}}{\partial a^\alpha} \quad ... \quad \frac{\partial f^i_{,j_1...j_s}}{\partial a^\alpha}\right) \tag{2.4.8}$$

We get thus a sequence of ranks

$$\mu_0 \leq \mu_1 \leq ... \leq \mu_s \leq ... \leq r \tag{2.4.9}$$

THEOREM 2.4.4. *The number of essential parameters in terms of which the functions $f^i(x,a)$ are expressible is equal to the maximum number attained in the sequence (2.4.9).*

PROOF. If any $\mu_s$ is equal to $r$, $\phi^\alpha = 0$ and $r$ parameters are essential.

Suppose that $\mu_{s-1} < r$ and $\mu_s = \mu_{s-1}$. From this statement, it follows that

$$\frac{\partial f^i_{,j_1...j_s}}{\partial a^\alpha} = \sum_{t<s} \lambda^{k_1...k_t}_{j_1...j_s} \frac{\partial f^i_{,k_1...k_t}}{\partial a^\alpha} \tag{2.4.10}$$

Differentiating equation (2.4.10) with respect to $x^p$, we have

$$\frac{\partial f^i_{,j_1...j_s p}}{\partial a^\alpha} = \sum_{t<s} \lambda^{k_1...k_t}_{j_1...j_s} \frac{\partial f^i_{,k_1...k_t p}}{\partial a^\alpha}$$

$$= \sum_{t<s-1} \lambda^{k_1...k_t}_{j_1...j_s} \frac{\partial f^i_{,k_1...k_t p}}{\partial a^\alpha} + \sum_{t=s-1, u<s} \lambda^{k_1...k_t}_{j_1...j_s} \lambda^{l_1...l_u}_{k_1...k_t p} \frac{\partial f^i_{,l_1...l_u}}{\partial a^\alpha}$$

Therefore, $\mu_{s+1} = \mu_{s-1}$.

If $\mu_s < r$, the system of linear equations (2.4.2), (2.4.7) has rank $\mu_s$, then $r - \mu_s$ of the functions $\phi^\alpha$ may be choosen arbitrary and the others are then determined. Let $\phi^1, ..., \phi^{r-\mu_s}$ be functions of $a$ but not the $x$. Then

$$\phi^\sigma = \lambda^\sigma_\rho \phi^\rho \quad \rho = 1,...,r-\mu_s \quad \sigma = r-\mu_s+1,...,r$$

where $\lambda$ are functions of $x$ and $a$. There are accordingly $r - \mu_s$ independent equations (2.4.2). The commutator of any two equations equated to 0 admits $f^i$ as solution and hence it must be a linear combination of the given equations. Therefore, these equations form a complete system. Hence the functions are expressible in terms of $\mu_s$ essential parameters. □

CHAPTER 3

# Lie Group

## 3.1. Lie Group Composition

Let $G$ be an $r$-parameter Lie group with neutral element $e$. Let
$$a \to (a^1, ..., a^n)$$
be local coordinates on Lie group $G$. Let
$$a_3^L = \varphi^L(a_1, a_2) \quad a_3 = a_1 a_2 \tag{3.1.1}$$
be an operation on Lie group $G$.

All operators that we introduce below have relation to left and right shift representations and will have an additional index $r$ or $l$ telling us what kind of shift they describe.

First of all we introduce operators
$$A_r{}_L^K(a,b) = \frac{\partial \varphi^K(a,b)}{\partial b^L} \tag{3.1.2}$$
$$A_l{}_L^K(a,b) = \frac{\partial \varphi^K(a,b)}{\partial a^L} \tag{3.1.3}$$
as derivative of left and right shifts.

THEOREM 3.1.1.
$$A_r(a, bc) A_r(b, c) = A_r(ab, c) \tag{3.1.4}$$
$$A_l(ab, c) A_l(a, b) = A_l(a, bc) \tag{3.1.5}$$

PROOF. The operation is associative a(bc)=(ab)c Using chain rule we can find the derivatives of this equation
$$\frac{\partial a(bc)}{\partial bc} \frac{\partial bc}{\partial c} = \frac{\partial (ab)c}{\partial c} \tag{3.1.6}$$
$$\frac{\partial a(bc)}{\partial a} = \frac{\partial (ab)c}{\partial ab} \frac{\partial ab}{\partial a} \tag{3.1.7}$$
(3.1.4) folows from (3.1.6). (3.1.5) folows from (3.1.7). □

Because $ae = a$ and $eb = b$ we get equations
$$A_l{}_L^K(a, e) = \delta_L^K \tag{3.1.8}$$
$$A_r{}_L^K(e, b) = \delta_L^K \tag{3.1.9}$$

THEOREM 3.1.2. *Operator $A_l(a,b)$ has the inverse operator $A_l(ab, b^{-1})$*
$$A_l^{-1}(a, b) = A_l(ab, b^{-1}) \tag{3.1.10}$$

Operator $A_r(b,c)$ has the inverse operator $A_r(b^{-1},bc)$

(3.1.11) $$A_r^{-1}(b,c) = A_r(b^{-1},bc)$$

PROOF. According (3.1.4) we get $(a = b^{-1})$

(3.1.12) $$A_r(b^{-1},bc)A_r(b,c) = A_r(b^{-1}b,c) = A_r(e,c) = \delta$$

(3.1.11) follows from (3.1.12).

According (3.1.5) we get $(c = b^{-1})$

(3.1.13) $$A_l(ab,b^{-1})A_l(a,b) = A_l(a,bb^{-1}) = A_l(a,e) = \delta$$

(3.1.10) follows from (3.1.13). $\square$

THEOREM 3.1.3. *Operators $A_l(a,b)$ and $A_r(a,b)$ are invertible.*

PROOF. This is consequence of theorem 3.1.2. $\square$

We introduce **Lie group basic operators**

(3.1.14) $$\psi_{r\,N}^L(a) = A_{r\,N}^L(a,e) \qquad \psi_r(a) = A_r(a,e)$$

(3.1.15) $$\psi_{l\,N}^L(b) = A_{l\,N}^L(e,b) \qquad \psi_l(b) = A_l(e,b)$$

From (3.1.8) and (3.1.9) immediately follows that

(3.1.16) $$\psi_{l\,L}^K(e) = \delta_L^K$$

(3.1.17) $$\psi_{r\,L}^K(e) = \delta_L^K$$

By definition basic operators linearly map the tangent plane $T_eG$ into the tangent plane $T_aG$.

THEOREM 3.1.4. *Operators $\psi_l$ and $\psi_r$ are invertible.*

PROOF. This is consequence of theorem 3.1.3 and definitions (3.1.15) and (3.1.14). $\square$

Because operators $\psi_r$ and $\psi_l$ have inverse operators we introduce operators

(3.1.18) $$\lambda_r(a) = \psi_r^{-1}(a)$$

(3.1.19) $$\lambda_l(a) = \psi_l^{-1}(a)$$

which map $T_aG \to T_eG$

THEOREM 3.1.5.

(3.1.20) $$\lambda_r(a) = A_r(a^{-1},a)$$

(3.1.21) $$\lambda_l(a) = A_l(a,a^{-1})$$

PROOF. According (3.1.19), (3.1.15) and (3.1.10) we get

$$\lambda_l(a) = A_l^{-1}(e,a) = A_l(a,a^{-1})$$

This proves (3.1.21).

According (3.1.18), (3.1.14) and (3.1.11) we get

$$\lambda_r(a) = A_r^{-1}(a,e) = A_r(a^{-1},a)$$

This proves (3.1.20). $\square$

THEOREM 3.1.6.

(3.1.22) $$A_r(a,b) = \psi_r(ab)\lambda_r(b)$$
(3.1.23) $$A_l(a,b) = \psi_l(ab)\lambda_l(a)$$

PROOF. Assume $c = e$ in (3.1.4)

(3.1.24) $$A_r(a,b)A_r(b,e) = A_r(ab,e)$$

From equations (3.1.14), (3.1.24), it follows that

(3.1.25) $$A_r(a,b)\psi_r(b) = \psi_r(ab)$$

From (3.1.25) follows

(3.1.26) $$A_r(a,b) = \psi_r(ab)\psi_r^{-1}(b)$$

(3.1.22) follows from (3.1.26) using (3.1.18).

Assume $a = e$ into (3.1.5)

(3.1.27) $$A_l(b,c)A_l(e,b) = A_l(e,bc)$$

From equations (3.1.15), (3.1.27), it follows that

(3.1.28) $$A_l(b,c)\psi_l(b) = \psi_l(bc)$$

From (3.1.28) follows

(3.1.29) $$A_l(a,b) = \psi_l(ab)\psi_l^{-1}(a)$$

(3.1.23) follows from (3.1.29) using (3.1.19). □

THEOREM 3.1.7. *The Lie group operation satisfies to differential equations*

(3.1.30) $$\frac{\partial \varphi^K(a,b)}{\partial b^L} = \psi_{rT}^K(ab)\lambda_{rL}^T(b) \qquad \frac{\partial ab}{\partial b} = \psi_r(ab)\lambda_r(b)$$

(3.1.31) $$\frac{\partial \varphi^K(a,b)}{\partial a^L} = \psi_{lT}^K(ab)\lambda_{lL}^T(a) \qquad \frac{\partial ab}{\partial a} = \psi_l(ab)\lambda_l(a)$$

PROOF. The equation (3.1.30) follows from equations (3.1.22) and (3.1.2). The equation (3.1.31) follows from equations (3.1.23) and (3.1.3). □

THEOREM 3.1.8.

(3.1.32) $$\frac{\partial a^{-1}}{\partial a} = -\psi_l(a^{-1})\lambda_r(a)$$

(3.1.33) $$\frac{\partial a^{-1}}{\partial a} = -\psi_r(a^{-1})\lambda_l(a)$$

PROOF. Differentiating the equation $e = a^{-1}a$ with respect to $a$, we get the equation

$$0 = \frac{\partial a^{-1}a}{\partial a^{-1}}\frac{\partial a^{-1}}{\partial a} + \frac{\partial a^{-1}a}{\partial a} =$$
$$= A_l(a^{-1},a)\frac{\partial a^{-1}}{\partial a} + A_r(a^{-1},a)$$
$$\frac{\partial a^{-1}}{\partial a} = -A_l^{-1}(a^{-1},a)A_r(a^{-1},a)$$

Using (3.1.20) and (3.1.21) we get

(3.1.34) $$\frac{\partial a^{-1}}{\partial a} = -\lambda_l^{-1}(a^{-1})\lambda_r(a)$$

(3.1.32) follows from (3.1.34).

Differentiating the equation $e = aa^{-1}$ with respect to $a$, we get the equation

$$0 = \frac{\partial aa^{-1}}{\partial a} + \frac{\partial aa^{-1}}{\partial a^{-1}}\frac{\partial a^{-1}}{\partial a} =$$

$$= A_l(a, a^{-1}) + A_r(a, a^{-1})\frac{\partial a^{-1}}{\partial a}$$

$$\frac{\partial a^{-1}}{\partial a} = -A_r^{-1}(a, a^{-1})A_l(a, a^{-1})$$

Using (3.1.20) and (3.1.21) we get

(3.1.35) $$\frac{\partial a^{-1}}{\partial a} = -\lambda_r^{-1}(a^{-1})\lambda_l(a)$$

(3.1.33) follows from (3.1.35). □

THEOREM 3.1.9.

(3.1.36) $$\frac{\partial a^{-1}b}{\partial a} = -\psi_l(a^{-1}b)\psi_r^{-1}(a)$$

(3.1.37) $$\frac{\partial ba^{-1}}{\partial a} = -\psi_r(ba^{-1})\psi_l^{-1}(a)$$

PROOF. Using the chain rule and equations (3.1.23), (3.1.32) we get

(3.1.38) $$\begin{aligned}\frac{\partial a^{-1}b}{\partial a} &= \frac{\partial a^{-1}b}{\partial a^{-1}}\frac{\partial a^{-1}}{\partial a} \\ &= -A_l(a^{-1}, b)\psi_l(a^{-1})\lambda_r(a) \\ &= -\psi_l(a^{-1}b)\psi_l^{-1}(a^{-1})\psi_l(a^{-1})\lambda_r(a)\end{aligned}$$

(3.1.36) follows from (3.1.38).

Using the chain rule and equations (3.1.22), (3.1.33) we get

(3.1.39) $$\begin{aligned}\frac{\partial ba^{-1}}{\partial a} &= \frac{\partial ba^{-1}}{\partial a^{-1}}\frac{\partial a^{-1}}{\partial a} \\ &= -A_r(b, a^{-1})\psi_r(a^{-1})\lambda_l(a) \\ &= -\psi_r(ba^{-1})\psi_r^{-1}(a^{-1})\psi_r(a^{-1})\lambda_l(a)\end{aligned}$$

(3.1.37) follows from (3.1.39). □

THEOREM 3.1.10.

(3.1.40) $$\frac{\partial abc}{\partial b} = \psi_l(abc)\lambda_l(ab)\psi_r(ab)\lambda_r(b)$$

(3.1.41) $$\frac{\partial abc}{\partial b} = \psi_r(abc)\lambda_r(bc)\psi_l(bc)\lambda_l(b)$$

PROOF. Using the chain rule and equations (3.1.30), (3.1.31) we get

(3.1.42) $$\frac{\partial abc}{\partial b} = \frac{\partial abc}{\partial ab}\frac{\partial ab}{\partial b} = \psi_l(abc)\lambda_l(ab)\psi_r(ab)\lambda_r(b)$$

The equation (3.1.40) follows from the equation (3.1.42). Using the chain rule and equations (3.1.30), (3.1.31) we get

(3.1.43) $$\frac{\partial abc}{\partial b} = \frac{\partial abc}{\partial bc}\frac{\partial bc}{\partial b} = \psi_r(abc)\lambda_r(bc)\psi_l(bc)\lambda_l(b)$$

3.2. 1-Parameter Group        19The equation (3.1.41) follows from the equation (3.1.43). □

THEOREM 3.1.11.

$$\frac{\partial aba^{-1}}{\partial a} = (\psi_l(aba^{-1}) - \psi_r(aba^{-1}))\lambda_l(a) \tag{3.1.44}$$

$$\frac{\partial aba^{-1}}{\partial b} = \psi_l(aba^{-1})\lambda_l(ab)\psi_r(ab)\lambda_r(b) \tag{3.1.45}$$

$$\frac{\partial aba^{-1}}{\partial b} = \psi_r(aba^{-1})\lambda_r(ba^{-1})\psi_l(ba^{-1})\lambda_l(b) \tag{3.1.46}$$

PROOF. The equation

$$\begin{aligned}\frac{\partial aba^{-1}}{\partial a} &= \frac{\partial a(ba^{-1})}{\partial a} + \frac{\partial (ab)a^{-1}}{\partial a^{-1}}\frac{\partial a^{-1}}{\partial a} \\ &= \psi_l(aba^{-1})\lambda_l(a) - \psi_r(aba^{-1})\lambda_r(a^{-1})\psi_r(a^{-1})\lambda_l(a)\end{aligned} \tag{3.1.47}$$

follows from equations (3.1.30), (3.1.31), (3.1.33). The equation

$$\frac{\partial aba^{-1}}{\partial a} = \psi_l(aba^{-1})\lambda_l(a) - \psi_r(aba^{-1})\lambda_l(a) \tag{3.1.48}$$

follows from equations (3.1.47), (3.1.18). The equation (3.1.44) follows from the equation (3.1.48).

The equation (3.1.45) follows from the equation (3.1.40). The equation (3.1.46) follows from the equation (3.1.41). □

THEOREM 3.1.12.

$$\left.\frac{\partial aba^{-1}}{\partial b}\right|_{b=e} = \lambda_l(a)\psi_r(a) \tag{3.1.49}$$

$$\left.\frac{\partial aba^{-1}}{\partial b}\right|_{b=e} = \lambda_r(a^{-1})\psi_l(a^{-1}) \tag{3.1.50}$$

$$\lambda_l(a)\psi_r(a) = \lambda_r(a^{-1})\psi_l(a^{-1}) \tag{3.1.51}$$

PROOF. The equation (3.1.49) follows from the equation (3.1.45). The equation (3.1.50) follows from the equation (3.1.46). The equation (3.1.51) follows from equations (3.1.49), (3.1.50). □

## 3.2. 1-Parameter Group

If the manifold of the group has the dimension 1, we have the group that depends on 1 parameter. In this case the operation on the group is just

$$c = \varphi(a,b)$$

where $a, b, c$ are numbers. In this case our notation will be simplified. We have not operators but functions

$$A_l(a,b) = \frac{\partial \varphi(a,b)}{\partial a}$$
$$A_r(a,b) = \frac{\partial \varphi(a,b)}{\partial b}$$

which by the definition satisfy equations

$$A_l(a, e) = 1$$
$$A_r(e, b) = 1$$

As well we introduce basic functions

$$\psi_r(a) = A_r(a, e)$$
$$\psi_l(b) = A_l(e, b)$$

Functions $\psi_r(a)$ and $\psi_l(a)$ never turn to 0 and

$$\lambda_r(a) = \frac{1}{\psi_r(a)}$$
$$\lambda_l(a) = \frac{1}{\psi_l(a)}$$

THEOREM 3.2.1. *The operation on the 1-parameter Lie group satisfies to differential equations*

$$(3.2.1) \qquad \frac{\partial \varphi(a, b)}{\partial b} = \frac{\psi_r(ab)}{\psi_r(b)}$$

PROOF. This is consequence of (3.1.30). □

THEOREM 3.2.2. *On the 1-parameter group we can introduce such coordinate $A$ that the operation $\Phi$ on the group has form*

$$(3.2.2) \qquad \Phi(A, B) = A + B$$

*and the identity of the group is $E = 0$.*

PROOF. We introduce a new variable $A$ such that

$$(3.2.3) \qquad dA = \frac{da}{\psi_r(a)}$$

Because $\psi_r(a) \neq 0$ 1-1 map $A = F(a)$ and its inverse map $a = f(A)$ exist. This means that if $A = F(a)$, $B = F(b)$, $C = F(c)$ and $c = \varphi(a, b)$ then there exist a function $\Phi$ such that

$$C = \Phi(A, B) = F(\varphi(a, b))$$

Now we have derivative

$$\frac{\partial \Phi(A, B)}{\partial B} = \frac{dC}{dc} \frac{\partial \varphi(a, b)}{\partial b} \frac{db}{dB}$$

Using (3.2.3) and (3.2.1) we get

$$\frac{\partial \Phi(A, B)}{\partial B} = \frac{1}{\psi_r(c)} \frac{\psi_r(c)}{\psi_r(b)} \psi_r(b) = 1$$

Thus we have

$$(3.2.4) \qquad \Phi(A, B) = \xi(A) + B$$

If we get the solution of equation (3.2.3) in the form

$$A = \int_e^a \frac{da}{\psi_r(a)}$$

then we see that the identity of group is $A = 0$. If we assume $B = 0$ in (3.2.4) we get

$$(3.2.5) \qquad A = \xi(A)$$

(3.2.2) follows from (3.2.4) and (3.2.5) □

### 3.3. Right Shift

For the right shift
$$b' = ba$$
the system of differential equations (3.1.30) gets form

$$(3.3.1) \qquad \frac{\partial b'^K}{\partial a^L} = \psi_{rT}^K(b')\lambda_{rL}^T(a)$$

Maps $b'^K$ are solutions of the system (3.3.1) and according to [2]-(2.4.14) they depend on $b^1, ..., b^n$ which we can assume as constants. Thus solution of the system (3.3.1) depends on $n$ arbitrary constants and therefore the system (3.3.1) is completely integrable. Condition of its integrability has the form

$$\frac{\partial \psi_{rT}^L(b')}{\partial b'^S}\frac{\partial b'^S}{\partial a^R}\lambda_{rP}^T(a) + \psi_{rT}^L(b')\frac{\partial \lambda_{rP}^T(a)}{\partial a^R} =$$
$$= \frac{\partial \psi_{rT}^L(b')}{\partial b'^S}\frac{\partial b'^S}{\partial a^P}\lambda_{rR}^T(a) + \psi_{rT}^L(b')\frac{\partial \lambda_{rR}^T(a)}{\partial a^P}$$

Using (3.3.1) we can write this condition in more simple form

$$\frac{\partial \psi_{rT}^L(b')}{\partial b'^S}\psi_{rV}^S(b')\lambda_{rR}^V(a)\lambda_{rP}^T(a) + \psi_{rT}^L(b')\frac{\partial \lambda_{rP}^T(a)}{\partial a^R}$$
$$= \frac{\partial \psi_{rT}^L(b')}{\partial b'^S}\psi_{rV}^S(b')\lambda_{rP}^V(a)\lambda_{rR}^T(a) + \psi_{rT}^L(b')\frac{\partial \lambda_{rR}^T(a)}{\partial a^P}$$

$$(3.3.2) \qquad \begin{aligned} &\frac{\partial \psi_{rT}^L(b')}{\partial b'^S}\psi_{rV}^S(b') - \frac{\partial \psi_{rV}^L(b')}{\partial b'^S}\psi_{rT}^S(b') \\ &= \psi_{rU}^L(b')\psi_{rT}^P(a)\psi_{rV}^R(a)\left(\frac{\partial \lambda_{rR}^U(a)}{\partial a^P} - \frac{\partial \lambda_{rP}^U(a)}{\partial a^R}\right) \end{aligned}$$

Let

$$(3.3.3) \qquad C_{rVT}^U = \psi_{rV}^R(a)\psi_{rT}^P(a)\left(\frac{\partial \lambda_{rR}^U(a)}{\partial a^P} - \frac{\partial \lambda_{rP}^U(a)}{\partial a^R}\right)$$

The equation

$$(3.3.4) \qquad \frac{\partial \psi_{rT}^L(b')}{\partial b'^S}\psi_{rV}^S(b') - \frac{\partial \psi_{rV}^L(b')}{\partial b'^S}\psi_{rT}^S(b') = C_{rVT}^U\psi_{rU}^L(b')$$

follows from equations (3.3.2), (3.3.3). If we differentiate this equation with respect to $a^P$ we get

$$\frac{\partial C_{rVT}^U}{\partial a^P}\psi_{rU}^L(b') = 0$$

because $\psi_{rU}^L(b')$ does not depend on $a$. At the same time $\psi_{rU}^L(b')$ are line independent because $det\|\psi_{rU}^L\| \neq 0$. Therefore

$$\frac{\partial C_{rVT}^U}{\partial a^P} = 0$$

and $C_{rTV}^U$ are constants which we call **right structural constant of Lie algebra**.
From (3.3.3), it follows that

$$C_{rTV}^U \lambda_{rP}^T(a) \lambda_{rR}^V(a) = \frac{\partial \lambda_{rP}^U(a)}{\partial a^R} - \frac{\partial \lambda_{rR}^U(a)}{\partial a^P} \tag{3.3.5}$$

We call (3.3.5) Maurer equation.

THEOREM 3.3.1. *Vector fields defined by differential operator*

$$X_{rV} = \psi_{rV}^S(a) \frac{\partial}{\partial a^S} \tag{3.3.6}$$

*are line independent and their commutator is*

$$(X_{rT}, X_{rV}) = C_{rTV}^U X_{rU}$$

PROOF. Line independence of vector fields follows from theorem 3.1.4. Then we see that according to (3.3.4)

$$(X_{rT}, X_{rV}) = (X_{rT} \psi_{rV}^D - X_{rV} \psi_{rT}^D) \frac{\partial}{\partial a^D} =$$
$$= \left( \psi_{rT}^P(a) \frac{\partial \psi_{rV}^D(a)}{\partial a^P} - \psi_{rV}^R(a) \frac{\partial \psi_{rT}^D(a)}{\partial a^R} \right) \frac{\partial}{\partial a^D} =$$
$$= C_{rTV}^U \psi_{rU}^D(a) \frac{\partial}{\partial a^D} = C_{rTV}^U X_{rU}$$

□

Let we have the homomorphism $f : G_1 \to G$ of the 1-parameter Lie group $G_1$ into the group $G$. Image of this group is the 1-parameter subgroup. If $t$ is coordinate on the group $G_1$ we can write $a = f(t)$ and find out differential equation for this subgroup. We assume in case of right shift that $a = f(t_1)$, $b = f(t_2)$, $c = ab = f(t)$, $t = t_1 + t_2$. Then we have

$$\frac{dc^K}{dt} = \frac{\partial c^K}{\partial b^L} \frac{db^L}{dt} = \psi_{rT}^K(c) \lambda_{rL}^T(b) \frac{db^L}{dt_2} \frac{dt_2}{dt}$$
$$\frac{dc^K}{dt} = \psi_{rT}^K(c) \lambda_{rL}^T(b) \frac{db^L}{dt_2}$$

Left part does not depend on $t_2$, therefore right part does not depend on $t_2$. We assume that

$$\lambda_{rL}^T(b) \frac{db^L}{dt_2} = \alpha^T$$

Thus we get system of differential equations

$$\frac{dc^K}{dt} = \psi_{rT}^K(c) \alpha^T$$

Because $\psi_r$ is derivative of right shift at identity of group this equation means that 1-parameter group is determined by vector $\alpha^T \in T_e G$ and transfers this vector along 1-parameter group without change. We call this vector field **right invariant vector field**. We introduce vector product on $T_e$ as

$$[\alpha, \beta]^T = C_{rRS}^T \alpha^R \beta^S \tag{3.3.7}$$

Space $T_e G$ equipped by such operation becomes Lie algebra $\mathfrak{g}_r$. We call it **right defined Lie algebra of Lie group**

THEOREM 3.3.2. *Space of right invariant vector fields has finite dimension equal of dimension of Lie group. It is Lie algebra with product equal to commutator of vector fields and this algebra is isomorphic to Lie algebra $\mathfrak{g}_r$.*

PROOF. It follows from (3.3.6) and (3.3.7) because $\alpha^K$ and $\beta^K$ are constants
$$(X_{rT}\alpha^T, X_{rV}\beta^V) = (X_{rT}, X_{rV})\alpha^T\beta^V =$$
$$= C_{rTV}^U X_{rU}\alpha^T\beta^V = [\alpha, \beta]^U X_{rU}$$

□

## 3.4. Left Shift

For the left shift
$$b' = ab$$
the system of differential equations (3.1.31) gets form

(3.4.1) $$\frac{\partial b'^K}{\partial a^L} = \psi_{lT}^K(b')\lambda_{lL}^T(a)$$

Functions $b'^K$ are solutions of system (3.4.1) and according to [2]-(2.4.12) they depend on $b^1, ..., b^n$ which we can assume as constants. Thus solution of system (3.4.1) depends on $n$ arbitrary constants and therefore the system (3.4.1) is completely integrable. Condition of its integrability has the form

$$\frac{\partial \psi_{lT}^L(b')}{\partial b'^S}\frac{\partial b'^S}{\partial a^R}\lambda_{lP}^T(a) + \psi_{lT}^L(b')\frac{\partial \lambda_{lP}^T(a)}{\partial a^R} =$$
$$= \frac{\partial \psi_{lT}^L(b')}{\partial b'^S}\frac{\partial b'^S}{\partial a^P}\lambda_{lR}^T(a) + \psi_{lT}^L(b')\frac{\partial \lambda_{lR}^T(a)}{\partial a^P}$$

Using (3.4.1) we can write this condition in more simple form

$$\frac{\partial \psi_{lT}^L(b')}{\partial b'^S}\psi_{lV}^S(b')\lambda_{lR}^V(a)\lambda_{lP}^T(a) + \psi_{lT}^L(b')\frac{\partial \lambda_{lP}^T(a)}{\partial a^R}$$
$$= \frac{\partial \psi_{lT}^L(b')}{\partial b'^S}\psi_{lV}^S(b')\lambda_{lP}^V(a)\lambda_{lR}^T(a) + \psi_{lT}^L(b')\frac{\partial \lambda_{lR}^T(a)}{\partial a^P}$$

(3.4.2)
$$\frac{\partial \psi_{lT}^L(b')}{\partial b'^S}\psi_{lV}^S(b') - \frac{\partial \psi_{lV}^L(b')}{\partial b'^S}\psi_{lT}^S(b')$$
$$= \psi_{lU}^L(b')\psi_{lT}^P(a)\psi_{lV}^R(a)\left(\frac{\partial \lambda_{lR}^U(a)}{\partial a^P} - \frac{\partial \lambda_{lP}^U(a)}{\partial a^R}\right)$$

Let

(3.4.3) $$C_{lVT}^U = \psi_{lV}^R(a)\psi_{lT}^P(a)\left(\frac{\partial \lambda_{lR}^U(a)}{\partial a^P} - \frac{\partial \lambda_{lP}^U(a)}{\partial a^R}\right)$$

The equation

(3.4.4) $$\frac{\partial \psi_{lT}^L(b')}{\partial b'^S}\psi_{lV}^S(b') - \frac{\partial \psi_{lV}^L(b')}{\partial b'^S}\psi_{lT}^S(b') = C_{lVT}^U\psi_{lU}^L(b')$$

follows from equations (3.4.2), (3.4.3). If we differentiate this equation with respect to $a^P$ we get

$$\frac{\partial C_{lVT}^U}{\partial a^P}\psi_{lU}^L(b') = 0$$

because $\psi_{lU}^L(b')$ does not depend on $a$. At the same time $\psi_{lU}^L(b')$ are line independent because $det\|\psi_{lU}^L\| \neq 0$. Therefore

$$\frac{\partial C_{lVT}^U}{\partial a^P} = 0$$

and $C_{lTV}^U$ are constants which we call **left structural constant of Lie algebra**. From (3.4.3), it follows that

(3.4.5) $$C_{lTV}^U \lambda_{lP}^T(a) \lambda_{lR}^V(a) = \frac{\partial \lambda_{lP}^U(a)}{\partial a^R} - \frac{\partial \lambda_{lR}^U(a)}{\partial a^P}$$

We call (3.4.5) Maurer equation.

THEOREM 3.4.1. *Vector fields defined by differential operator*

(3.4.6) $$X_{lV} = \psi_{lV}^S(a) \frac{\partial}{\partial a^S}$$

*are line independent and their commutator is*

$$(X_{lT}, X_{lV}) = C_{lTV}^U X_{lU}$$

PROOF. Line independence of vector fields follows theorem 3.1.4. Then we see that according to (3.4.4)

$$(X_{lT}, X_{lV}) = (X_{lT}\psi_{lV}^D - X_{lV}\psi_{lT}^D)\frac{\partial}{\partial a^D} =$$
$$= \left(\psi_{lT}^P(a)\frac{\partial \psi_{lV}^D(a)}{\partial a^P} - \psi_{lV}^R(a)\frac{\partial \psi_{lT}^D(a)}{\partial a^R}\right)\frac{\partial}{\partial a^D} =$$
$$= C_{lTV}^U \psi_{lU}^D(a)\frac{\partial}{\partial a^D} = C_{lTV}^U X_U$$

□

Let we have the homomorphism $f : G_1 \to G$ of the 1-parameter Lie group $G_1$ into the group $G$. Image of this group is the 1-parameter subgroup. If $t$ is coordinate on the group $G_1$ we can write $a = f(t)$ and find out differential equation for this subgroup. We assume in case of left shift that $a = f(t_1)$, $b = f(t_2)$, $c = ab = f(t)$, $t = t_1 + t_2$. Then we have

$$\frac{dc^K}{dt} = \frac{\partial c^K}{\partial b^L}\frac{db^L}{dt} = \psi_{lT}^K(c)\lambda_{lL}^T(b)\frac{db^L}{dt_2}\frac{dt_2}{dt}$$

$$\frac{dc^K}{dt} = \psi_{lT}^K(c)\lambda_{lL}^T(b)\frac{db^L}{dt_2}$$

Left part does not depend on $t_2$, therefore right part does not depend on $t_2$. We assume that

$$\lambda_{lL}^T(b)\frac{db^L}{dt_2} = \alpha^T$$

Thus we get system of differential equations

$$\frac{dc^K}{dt} = \psi_{lT}^K(c)\alpha^T$$

Because $\psi_l$ is derivative of right shift at identity of group this equation means that 1-parameter group is determined by vector $\alpha^T \in T_e G$ and transfers this vector along 1-parameter group without change. We call this vector field **left invariant vector field**. We introduce vector product on $T_e$ as

(3.4.7) $$[\alpha,\beta]^T = C_{l\,RS}^{\ \ T}\alpha^R\beta^S$$

Space $T_eG$ equipped by such operation becomes Lie algebra $\mathfrak{g}_l$. We call it **left defined Lie algebra of Lie group**

THEOREM 3.4.2. *Space of right invariant vector fields has finite dimension equal of dimension of Lie group. It is Lie algebra with product equal to commutator of vector fields and this algebra is isomorphic to Lie algebra $\mathfrak{g}_l$.*

PROOF. It follows from (3.4.6) and (3.4.7) because $\alpha^K$ and $\beta^K$ are constants

$$(X_{lT}\alpha^T, X_{lV}\beta^V) = (X_{lT}, X_{lV})\alpha^T\beta^V =$$
$$= C_{l\,TV}^{\ \ U}X_{lU}\alpha^T\beta^V = [\alpha,\beta]^U X_{lU}$$

□

### 3.5. Relation between Algebras Lie $\mathfrak{g}_l$ and $\mathfrak{g}_r$

We defined two different algebras Lie on space $T_eG$. Our goal now is to define relation between these algebras. We start to solve this problem from analysis of Taylor expansion of group operation.

THEOREM 3.5.1. *Up to second order of infinitesimal, structure of operation on group Lie*

(3.5.1) $$\varphi^K(a,b) = a^K + b^K - e^K + I_{LM}^K(a^M - e^M)(b^L - e^L)$$

*where we introduce* **infinitesimal generators of group Lie**

$$I_{LM}^K = I_{l\,LM}^{\ \ K} = I_{r\,ML}^{\ \ K}$$

PROOF. We can find Tailor expansion of product by both arguments. However like in sections 3.3 and 3.4 we can assume one argument as parameter and find Tailor expansion by another argument. In this case its coefficients will depend on former argument.

Thus according to (3.1.2) and (3.1.14), product in vicinity of $a \in G$ has Tailor expansion relative $b$

(3.5.2) $$\varphi^K(a,b) = a^K + \psi_{r\,L}^{\ \ K}(a)(b^L - e^L) + o(b^L - e^L)$$
$$ab = a + \psi_r(a)(b - e) + o(b - e)$$

where coefficients of expansion depend on $a$. At the same time $\psi_r(a)$ also has Tailor expansion in vicinity of $e$

(3.5.3) $$\psi_{r\,L}^{\ \ K}(a) = \delta_L^K + I_{r\,LM}^{\ \ K}(a^M - e^M) + o(a^M - e^M)$$

If we substitute (3.5.3) into (3.5.2) we get expansion

$$\varphi^K(a,b) = a^K + (\delta_L^K + I_{r\,LM}^{\ \ K}(a^M - e^M) + o(a^M - e^M))(b^L - e^L) + o(b^L - e^L)$$

(3.5.4) $\varphi^K(a,b) = a^K + b^K - e^K + I_{r\,LM}^{\ \ K}(a^M - e^M))(b^L - e^L) + o(a^M - e^M, b^L - e^L)$

According (3.1.3) and (3.1.15), product in vicinity of $b \in G$ has Tailor expansion relative $a$

(3.5.5) $$\varphi^K(a,b) = b^K + \psi_{l\,L}^{\ \ K}(b)(a^L - e^L) + o(a^L - e^L)$$
$$ab = b + \psi_l(b)(a - e) + o(a - e)$$

where coefficients of expansion depend on $b$. At the same time $\psi_l(b)$ also has Tailor expansion in vicinity of $e$

(3.5.6) $$\psi_l{}_L^K(b) = \delta_L^K + I_l{}_{LM}^K(b^M - e^M) + o(b^M - e^M)$$

If we substitute (3.5.6) into (3.5.5) we get expansion

$$\varphi^K(a,b) = b^K + (\delta_L^K + I_l{}_{LM}^K(b^M - e^M) + o(b^M - e^M))(a^L - a^L) + o(a^L - e^L)$$

(3.5.7) $\varphi^K(a,b) = b^K + a^K - e^K + I_l{}_{LM}^K(b^M - e^M))(a^L - e^L) + o(a^M - e^M, b^L - e^L)$

(3.5.4) and (3.5.7) are Tailor expansions of the same function amd they should be coincide. Comparing them we see that $I_l{}_{LM}^K = I_r{}_{ML}^K$. Then group Lie operation has Tailor expansion (3.5.1). $\square$

THEOREM 3.5.2. *Up to second order of infinitesimal, structure of operation on group Lie*

(3.5.8) $$\psi_r{}_L^K(a) = \delta_L^K + I_{ML}^K(a^M - e^M) + o(a^M - e^M)$$

(3.5.9) $$\psi_l{}_L^K(a) = \delta_L^K + I_{LM}^K(a^M - e^M) + o(a^M - e^M)$$

(3.5.10) $$\lambda_r{}_L^K(a) = \delta_L^K - I_{ML}^K(a^M - e^M) + o(a^M - e^M)$$

(3.5.11) $$\lambda_l{}_L^K(a) = \delta_L^K - I_{LM}^K(a^M - e^M) + o(a^M - e^M)$$

PROOF. (3.5.8) and (3.5.9) follow from differentiating of (3.5.1) with respect to $a$ or $b$ and definitions (3.1.14) and (3.1.15).

Now we assume that $\lambda_r{}_L^K(a)$ has Tailor expansion

(3.5.12) $$\lambda_r{}_L^K(a) = \delta_L^K + J_{LM}^K(a^M - e^M) + o(a^M - e^M)$$

Then substituting (3.5.12) and (3.5.8) into equation

$$\lambda_r{}_L^K(a)\psi_r{}_M^L(a) = \delta_M^K$$

we get up to order 1

$$(\delta_L^K + I_{NL}^K(a^N - e^N) + o(a^N - e^N))(\delta_M^L + J_{MP}^L(a^P - e^P) + o(a^P - e^P)) = \delta_M^K$$

$$\delta_M^K + (J_{MN}^K + I_{NM}^K)(a^N - e^N) + o(a^P - e^P) = \delta_M^K$$

Therefore $J_{MN}^K + I_{NM}^K = 0$ and substituting $J_{MN}^K$ into (3.5.12) we get (3.5.10).

The same way we can prove (3.5.11). $\square$

THEOREM 3.5.3.

(3.5.13) $$\left.\frac{\partial \psi_r{}_L^K(a)}{\partial a^M}\right|_{a=e} = I_{ML}^K$$

(3.5.14) $$\left.\frac{\partial \psi_l{}_L^K(a)}{\partial a^M}\right|_{a=e} = I_{LM}^K$$

(3.5.15) $$\left.\frac{\partial \lambda_r{}_L^K(a)}{\partial a^M}\right|_{a=e} = -I_{ML}^K$$

(3.5.16) $$\left.\frac{\partial \lambda_l{}_L^K(a)}{\partial a^M}\right|_{a=e} = -I_{LM}^K$$

PROOF. The theorem follows immediately from theorem 3.5.2. $\square$

On the base of theorem 3.5.3 we can get more detailed information about left and right group algebras Lie. Because definition (3.3.3) does not depend on $a$ we can estimate it when $a = e$. Using (3.1.16) and (3.5.16) we get

$$C^U_{lTV} = \delta^P_T \delta^R_V(-I^U_{PR} + I^U_{RP}) = I^U_{VT} - I^U_{TV} \tag{3.5.17}$$

Using (3.1.17) and (3.5.15) we get

$$C^U_{rTV} = \delta^P_T \delta^R_V(-I^U_{RP} + I^U_{PR}) = I^U_{TV} - I^U_{VT} \tag{3.5.18}$$

The equation

$$C^U_{rTV} = -C^U_{lTV}$$

follows from equations (3.5.17) and (3.5.18). Therefore, algebras $\mathfrak{g}_l$ and $\mathfrak{g}_r$ are anti isomorphic. Usually we asume that $\mathfrak{g}_l$ is Lie algebra of Lie group and use notation

$$C^U_{TV} = C^U_{lTV}$$

Let $\mathfrak{g} = T_e G$ be **Lie algebra of Lie group** $G$ with an operation

$$\alpha^L_3 = C^L_{MN} \alpha^M_1 \alpha^N_2 \quad \alpha_3 = C(\alpha_1, \alpha_2)$$

## 3.6. Group $Gl_n(F)$

Since we consider vector space of dimension $n$ over field $F$, then the set of transformations generates group $GL_n(F)$. Element of group $GL_n(F)$ is matrix $n \times n$. Product of matrices $a = (a^i_j)$ and $b = (b^i_j)$ is defined by the equation

$$(ab)^i_j = a^i_k b^k_j \tag{3.6.1}$$

THEOREM 3.6.1.

$$A_r{}^{\cdot k\,i}_{l\cdot j}(a, b) = a^k_j \delta^i_l \tag{3.6.2}$$

$$A_l{}^{\cdot k\,i}_{l\cdot j}(a, b) = b^i_l \delta^k_j \tag{3.6.3}$$

PROOF. The equation

$$A_r{}^{\cdot k\,i}_{l\cdot j}(a, b) = \frac{\partial a^k_m b^m_l}{\partial b^j_i} = a^k_m \delta^m_j \delta^i_l \tag{3.6.4}$$

follows from equations (3.1.2), (3.6.1). The equation (3.6.2) follows from the equation (3.6.4). The equation

$$A_l{}^{\cdot k\,i}_{l\cdot j}(a, b) = \frac{\partial a^k_m b^m_l}{\partial a^j_i} = \delta^k_j \delta^i_m b^m_l \tag{3.6.5}$$

follows from equations (3.1.3), (3.6.1). The equation (3.6.3) follows from the equation (3.6.5). □

THEOREM 3.6.2.

$$\psi_r{}^{\cdot k\,i}_{l\cdot j}(a) = a^k_j \delta^i_l \tag{3.6.6}$$

$$\psi_l{}^{\cdot k\,i}_{l\cdot j}(b) = b^i_l \delta^k_j \tag{3.6.7}$$

PROOF. The equation (3.6.6) follows from equations (3.1.14), (3.6.2). The equation (3.6.7) follows from equations (3.1.15), (3.6.3). □

Theorem 3.6.3.

$$\lambda_r{}^{\cdot k\,i}_{\,l\cdot j}(a) = a^{-1\cdot k}{}_j \delta^i_l \tag{3.6.8}$$

$$\lambda_l{}^{\cdot k\,i}_{\,l\cdot j}(b) = b^{-1\cdot i}{}_l \delta^k_j \tag{3.6.9}$$

Proof. The equation

$$\lambda_r{}^{\cdot k\,i}_{\,l\cdot j}(a) = A_r{}^{\cdot k\,i}_{\,l\cdot j}(a^{-1},a) = a^{-1\cdot k}{}_j \delta^i_l \tag{3.6.10}$$

follows from equations (3.1.20), (3.6.2). The equation (3.6.8) follows from the equation (3.6.10). The equation

$$\lambda_l{}^{\cdot k\,i}_{\,l\cdot j}(b) = A_l{}^{\cdot k\,i}_{\,l\cdot j}(b,b^{-1}) = b^{-1\cdot i}{}_l \delta^k_j \tag{3.6.11}$$

follows from equations (3.1.21), (3.6.3). The equation (3.6.9) follows from the equation (3.6.11). $\square$

Theorem 3.6.4.

$$I^{\cdot k\,i\ m}_{\,l\cdot j\cdot n} = \delta^i_n \delta^m_l \delta^k_j \tag{3.6.12}$$

$$C^{\cdot k\,i\ m}_{ll\cdot j\cdot n} = I^{\cdot k\,m\,i}_{\,l\cdot n\cdot j} - I^{\cdot k\,i\ m}_{\,l\cdot j\cdot n} = \delta^m_j \delta^i_l \delta^k_n - \delta^i_n \delta^m_l \delta^k_j \tag{3.6.13}$$

Proof. The equation (3.6.12) follows from equations (3.5.14), (3.6.7) and the equation

$$I^{\cdot k\,im}_{\,l\cdot jn} = \left.\frac{\partial \psi_l{}^{\cdot k\,i}_{\,l\cdot j}(a)}{\partial a^n_m}\right|_{a=e} = \left.\frac{\partial a^i_l \delta^k_j}{\partial a^n_m}\right|_{a=e} = \delta^i_n \delta^m_l \delta^k_j$$

The equation (3.6.13) follows from equations (3.5.17), (3.6.12). $\square$

CHAPTER 4

# Representation of Lie Group

## 4.1. Linear Representation of Lie Group

Let $F$ be field of characteristic 0, $V$ be $F$-vector space. According to the theorem [1]-4.4.3, linear transformation of $F$-vector space has form

(4.1.1) $$v' = a\,v$$

where $a \in \mathcal{L}(F, V, V)$. According to the theorem [1]-11.1.6, the set $GL(V)$ of linear transformations (4.1.1) is the group Lie which is called **linear transformation group**.

A representation[4.1]
$$f : G \longrightarrow V$$
of a Lie group $G$ in a continuous $F$-vector space $V$ is called **linear representation of Lie group** $G$ and has form
$$v' = f(a) \circ v$$
Let $\overline{\overline{e}}$ be the basis of $F$-vector space $V$. Then expansion of vector $v' = v'^{\alpha} e_{\alpha}$ has form
$$v'^{\alpha} = f^{\alpha}_{\beta}(a) v^{\beta}$$

THEOREM 4.1.1. *Let*

(4.1.2) $$v' = f(a) \circ v$$

*be linear representation of Lie group* $G$. *Maps* $f^{\alpha}_{\beta}$ *have the following properties*

(4.1.3) $$f^{\alpha}_{\beta}(e) = \delta^{\alpha}_{\beta}$$

(4.1.4) $$f^{\alpha}_{\beta}(a^{-1}) = f^{-1\cdot\alpha}{}_{\beta}(a)$$

(4.1.5) $$f^{\alpha}_{\beta}(ba) = f^{\alpha}_{\nu}(b) f^{\nu}_{\beta}(a)$$

PROOF. The equation (4.1.3) follows from the statement, that $a = e$ correspond to identity map.

Since the map (4.1.2) is representation of Lie group $G$, then the equation

(4.1.6) $$v''^{\alpha} = f^{\alpha}_{\beta}(ba) v^{\beta}$$

follows from the equation (4.1.2) and from the equation

(4.1.7) $$v''^{\alpha} = f^{\alpha}_{\beta}(b) v'^{\beta}$$

The equation

(4.1.8) $$v''^{\alpha} = f^{\alpha}_{\nu}(b) f^{\nu}_{\beta}(a) v^{\beta}$$

---

[4.1]We reviewed the general concept of representation of group in the section [2]-2.4.

follows immediately from equations (4.1.2), (4.1.7). The equation (4.1.5) follows from equations (4.1.6), (4.1.8). From equations (4.1.3), (4.1.5), it follows that

$$\delta^\alpha_\beta = f^\alpha_\nu(a^{-1}) f^\nu_\beta(a) \tag{4.1.9}$$

The equation (4.1.4) follows from equation (4.1.9). $\square$

THEOREM 4.1.2. *The linear representation of Lie group $G$*

$$v' = f(a) \circ v \tag{4.1.10}$$

*satisfies to system of differential equations*

$$\frac{\partial v'^\alpha}{\partial a^L} = I^\alpha_{\beta K} v'^\beta \lambda^K_{l\,L}(a) \tag{4.1.11}$$

*where* **infinitesimal generators of representation** *are defined by equation*

$$I^\alpha_{\beta\,M} = \left.\frac{\partial f^\alpha_\beta(b)}{\partial b^M}\right|_{b=e} \qquad I = \left.\frac{\partial f(b)}{\partial b}\right|_{b=e} \tag{4.1.12}$$

PROOF. Let $a, b \in G$. Since the map (4.1.10) is representation of Lie group $G$, then the equation

$$v''^\alpha = f^\alpha_\beta(ba) v^\beta$$

follows from the equation (4.1.10) and from the equation

$$v''^\alpha = f^\alpha_\beta(b) v'^\beta \tag{4.1.13}$$

Let

$$c = ba \tag{4.1.14}$$

We express $b$ through $a$ and $c$ and differentiate the equation (4.1.13) with respect to $a$

$$\begin{aligned} 0 &= \frac{\partial f^\nu_\beta(b) v'^\beta}{\partial b^K} \frac{\partial b^K}{\partial a^L} + \frac{\partial f^\nu_\beta(b) v'^\beta}{\partial v'^\alpha} \frac{\partial v'^\alpha}{\partial a^L} \\ &= \frac{\partial f^\nu_\beta(b)}{\partial b^K} v'^\beta \frac{\partial b^K}{\partial a^L} + f^\nu_\alpha(b) \frac{\partial v'^\alpha}{\partial a^L} \end{aligned} \tag{4.1.15}$$

The equation

$$\frac{\partial v'^\alpha}{\partial a^L} = -f^{-1\cdot\alpha}_{\ \ \nu}(b) \frac{\partial f^\nu_\beta(b)}{\partial b^K} v'^\beta \frac{\partial b^K}{\partial a^L} \tag{4.1.16}$$

follows from the equation (4.1.15). Left part does not depend on $b$, so and right part does not depend on $b$. We define right part when $b = e$. From equations (4.1.14), (3.1.37), when $ca^{-1} = e$, it follows that

$$\left.\frac{\partial b}{\partial a}\right|_{ca^{-1}=e} = \left.\frac{\partial ca^{-1}}{\partial a}\right|_{ca^{-1}=e} = -\left.\psi_r(ca^{-1})\right|_{ca^{-1}=e} \psi^{-1}_l(a) = -\lambda_l(a) \tag{4.1.17}$$

Equation

$$\left. f^{-1\cdot\alpha}_{\ \ \nu}(b) \frac{\partial f^\nu_\beta(b)}{\partial b^K} v'^\beta \right|_{b=e} = I^\alpha_{\beta\,M} v'^\beta \tag{4.1.18}$$

follows from equations (4.1.12). Equation (4.1.11) follows from equations (4.1.18), (4.1.17), (4.1.16). $\square$

THEOREM 4.1.3. *The system of differential equations* (4.1.11) *is completely integrable system iff*

(4.1.19) $$I^\alpha_{\beta K} I^\beta_{\nu P} - I^\alpha_{\beta P} I^\beta_{\nu K} = I^\alpha_{\nu T} C^T_{l P K}$$

PROOF. Equation

(4.1.20) $$\begin{aligned}\frac{\partial^2 v'^\alpha}{\partial a^L \partial a^M} &= I^\alpha_{\beta K} \frac{\partial v'^\beta}{\partial a^M} \lambda_l{}^K_L(a) + I^\alpha_{\beta K} v'^\beta \frac{\partial \lambda_l{}^K_L(a)}{\partial a^M} \\ &= I^\alpha_{\beta K} I^\beta_{\nu P} v'^\nu \lambda_l{}^P_M(a) \lambda_l{}^K_L(a) + I^\alpha_{\beta K} v'^\beta \frac{\partial \lambda_l{}^K_L(a)}{\partial a^M}\end{aligned}$$

follows from the equation (4.1.11). The system of differential equations (4.1.11) is completely integrable system iff

(4.1.21) $$\frac{\partial^2 v'^\alpha}{\partial a^L \partial a^M} = \frac{\partial^2 v'^\alpha}{\partial a^M \partial a^L}$$

The equation

(4.1.22) $$\begin{aligned}& I^\alpha_{\beta K} I^\beta_{\nu P} v'^\nu \lambda_l{}^P_M(a) \lambda_l{}^K_L(a) + I^\alpha_{\beta K} v'^\beta \frac{\partial \lambda_l{}^K_L(a)}{\partial a^M} \\ &= I^\alpha_{\beta P} I^\beta_{\nu K} v'^\nu \lambda_l{}^P_M(a) \lambda_l{}^K_L(a) + I^\alpha_{\beta K} v'^\beta \frac{\partial \lambda_l{}^K_M(a)}{\partial a^L}\end{aligned}$$

follows from equations (4.1.20), (4.1.21). The equation

(4.1.23) $$\begin{aligned}(I^\alpha_{\beta K} I^\beta_{\nu P} - I^\alpha_{\beta P} I^\beta_{\nu K}) v'^\nu \lambda_l{}^P_M(a) \lambda_l{}^K_L(a) &= I^\alpha_{\beta K} v'^\beta \left( \frac{\partial \lambda_l{}^K_M(a)}{\partial a^L} - \frac{\partial \lambda_l{}^K_L(a)}{\partial a^M} \right) \\ &= I^\alpha_{\beta T} v'^\beta C^T_{l P K} \lambda_l{}^P_M(a) \lambda_l{}^K_L(a)\end{aligned}$$

follows from equations (3.4.5), (4.1.22). Since the matrix $\lambda_l$ is non singular matrix, then the equation

(4.1.24) $$(I^\alpha_{\beta K} I^\beta_{\nu P} - I^\alpha_{\beta P} I^\beta_{\nu K}) v'^\nu = I^\alpha_{\nu T} v'^\nu C^T_{l P K}$$

follows from the equation (4.1.23). Since the vector $v$ is arbitrary vector, then the equation (4.1.19) follows from the equation (4.1.24). □

REMARK 4.1.4. Let $F$ be a field and $V$ be $F$-vector space of dimension $n$. Linear transformation $a_1$ of $F$-vector space $V$

(4.1.25) $$v'^\alpha_1 = v^\beta_1 a_{1\beta}{}^\alpha$$

written relative to basis $\overline{\overline{e}}_{V1}$ maps a vector

$$v = v^\alpha_1 e_{V1 \cdot \alpha}$$

into a vector

$$v' = v'^\alpha_1 e_{V1 \cdot \alpha}$$

Passive transformation[4.2] $b$

(4.1.26) $$e_{V2 \cdot \alpha} = b^\beta_\alpha e_{V1 \cdot \beta}$$

maps the basis $\overline{\overline{e}}_{V1}$ into basis $\overline{\overline{e}}_{V2}$. Coordinates of vectors $v$, $v'$ relative to the basis $\overline{\overline{e}}_{V2}$

$$v = v^\alpha_2 e_{V2 \cdot \alpha}$$
$$v' = v'^\alpha_2 e_{V2 \cdot \alpha}$$

---

[4.2] See definition of passive transformation in sections [1]-5.2, [4]-5.2.

have form
$$v_2^\alpha = v_1^\beta b^{-1\cdot\alpha}{}_\beta$$
$$v_2'^\alpha = v_2'^\beta b^{-1\cdot\alpha}{}_\beta$$

Therefore, the linear transformation (4.1.25) relative to the basis $\overline{\overline{e}}_{V2}$ has form

(4.1.27) $$v_2'^\alpha = v_2^\beta a_{2\beta}^\alpha$$

(4.1.28) $$a_{2\beta}^\alpha = b^{-1\cdot\alpha}{}_\nu a_{1\lambda}^\nu b_\beta^\lambda$$

□

THEOREM 4.1.5. *Let*

(4.1.29) $$v_1'^\alpha = f_\beta^\alpha(a_1) v_1^\beta$$

*be linear geometric object.*[4.3] *Passive transformation of left shift*

(4.1.30) $$a_1 \to g a_1$$

*on Lie group $G$ generates transformation of infinitesimal generators*[4.4]

(4.1.31) $$I_{2\cdot\mu K}^{\alpha} = f_\mu^\beta(g) f^{-1\cdot\alpha}{}_\rho(g) I_{1\cdot\beta J}^{\rho} \lambda_{lM}^{J}(g) \psi_{rK}^{M}(g)$$

PROOF. According to equations (4.1.29), passive transformation of left shift (4.1.30) on Lie group $G$ generates coordinate transformation of geometric object

(4.1.32) $$v_2^\alpha = f_\beta^\alpha(g^{-1}) v_1^\beta$$

(4.1.33) $$v_2'^\alpha = f_\beta^\alpha(g^{-1}) v_1'^\beta$$

where we assume

(4.1.34) $$v_2'^\alpha = f_\beta^\alpha(a_2) v_2^\beta$$

Equations

(4.1.35) $$\frac{\partial v_1'^\alpha}{\partial a_1^L} = I_{1\cdot\beta K}^{\alpha} v_1'^\beta \lambda_{lL}^{K}(a_1)$$

(4.1.36) $$\frac{\partial v_2'^\alpha}{\partial a_2^L} = I_{2\cdot\beta K}^{\alpha} v_2'^\beta \lambda_{lL}^{K}(a_2)$$

follow from the equation (4.1.11). Using the chain rule we get

(4.1.37) $$\frac{\partial v_2'^\pi}{\partial a_2^L} = \frac{\partial v_2'^\pi}{\partial v_1'^\rho} \frac{\partial v_1'^\rho}{\partial a_1^M} \frac{\partial a_1^M}{\partial a_2^L}$$

The equation

(4.1.38) $$\frac{\partial v_2'^\pi}{\partial v_1'^\rho} = f_\rho^\pi(g^{-1})$$

follows from the equation (4.1.33). The equation

(4.1.39) $$v_1'^\alpha = f_\beta^\alpha(g a_2 g^{-1}) v_1^\beta$$

---

[4.3] See definition of geometric object in sections [1]-5.3, [4]-5.3.

[4.4] The theorem is based on the statement that basis vector of representation of Lie group is covariant vector in parameter space ([5], p. 27) and contravariant vector in space of representation ([5], p. 26).

follows from equations (4.1.32), (4.1.33), (4.1.34). From equations (4.1.29), (4.1.39), it follows that passive transformation of left shift (4.1.30) generates transformation in parameter space

$$a_1 = g a_2 g^{-1} \tag{4.1.40}$$

The equation

$$\begin{aligned}\frac{\partial a_1^M}{\partial a_2^L} &= \psi_{lN}^{M}(ga_2g^{-1})\lambda_{lK}^{N}(ga_2)\psi_{rI}^{K}(ga_2)\lambda_{rL}^{I}(a_2) \\ &= \psi_{lN}^{M}(a_1)\lambda_{lK}^{N}(ga_2)\psi_{rI}^{K}(ga_2)\lambda_{rL}^{I}(a_2)\end{aligned} \tag{4.1.41}$$

follows from equations (3.1.45), (4.1.40). The equation

$$\begin{aligned}&I_{2\cdot\beta K}^{\ \alpha}v_2'^{\beta}\lambda_{lL}^{K}(a_2) \\ &= f_\rho^{\alpha}(g^{-1})I_{1\cdot\beta J}^{\ \rho}v_1'^{\beta}\lambda_{lM}^{J}(a_1)\psi_{lN}^{M}(a_1)\lambda_{lK}^{N}(ga_2)\psi_{rI}^{K}(ga_2)\lambda_{rL}^{I}(a_2)\end{aligned} \tag{4.1.42}$$

follows from equations (4.1.35), (4.1.36), (4.1.37), (4.1.38), (4.1.41). The equation

$$\begin{aligned}&I_{2\cdot\beta K}^{\ \alpha}f_\mu^{\beta}(g^{-1})v_1'^{\mu} \\ &= f_\rho^{\alpha}(g^{-1})I_{1\cdot\beta J}^{\ rho}v_1'^{\beta}\lambda_{lM}^{J}(ga_2)\psi_{rI}^{M}(ga_2)\lambda_{rL}^{I}(a_2)\psi_{lK}^{L}(a_2)\end{aligned} \tag{4.1.43}$$

follows from equations (3.1.19), (4.1.33), (4.1.42). Left part of the equation (4.1.43) does not depend on $a_2$, so and right part does not depend on $a_2$. The equation

$$\begin{aligned}&I_{2\cdot\beta K}^{\ \alpha}f_\mu^{\beta}(g^{-1})v_1'^{\mu} \\ &= f_\rho^{\alpha}(g^{-1})I_{1\cdot\beta J}^{\ \rho}v_1'^{\beta}\lambda_{lM}^{J}(g)\psi_{rK}^{M}(g)\end{aligned} \tag{4.1.44}$$

follows from the equation (4.1.43) when $a_2 = e$. Since $v_1'$ is arbitrary vector, then the equation

$$I_{2\cdot\mu K}^{\ \alpha}f^{-1\cdot\mu}{}_{\beta}(g) = f^{-1\cdot\alpha}{}_{\rho}(g)I_{1\cdot\beta J}^{\ \rho}\lambda_{lM}^{J}(g)\psi_{rK}^{M}(g) \tag{4.1.45}$$

follows from equations (4.1.4), (4.1.44). The equation (4.1.31) follows from the equation (4.1.45). □

## 4.2. Conjugated Representation

DEFINITION 4.2.1. Let $V$ be vector space over field $F$. Vector space $V' = \mathcal{L}(F; V; F)$ is called **conjugated vector space**. □

According to the definition 4.2.1, elements of conjugated vector space are linear maps

$$u : V \to F \tag{4.2.1}$$

The linear map (4.2.1) is called **linear functional** on $F$-vector space $V$. We write the value of linear functional as

$$u \circ a = (u, a)$$

THEOREM 4.2.2. *Let $\overline{\overline{e}} = (e_\alpha \in V, \alpha \in I)$ be basis of $F$-vector space $V$. Let $V'$ be $F$-vector space, conjugated to $F$-vector space $V$. The set of vectors $\overline{\overline{e}} = (e^\alpha \in V, \alpha \in I)$ such that*

$$(e^\alpha, e_\beta) = e^\alpha \circ e_\beta = \delta_\beta^\alpha \tag{4.2.2}$$

*is basis of $F$-vector space $V'$.*

PROOF. From the equation (4.2.2) it follows that

$$(4.2.3) \qquad (e^\alpha, a) = (e^\alpha, e_\beta a^\beta) = (e^\alpha, e_\beta) a^\beta = \delta^\alpha_\beta a^\beta = a^\alpha$$

Let

$$u : V \to F$$

be linear functional. Then

$$(4.2.4) \qquad (u, a) = (u, e_\alpha a^\alpha) = (u, e_\alpha) a^\alpha = u_\alpha a^\alpha$$

where

$$(u, e_\alpha) = u_\alpha$$

From equations (4.2.3), (4.2.4) it follows that

$$(u, a) = u_\alpha (e^\alpha, a)$$

According to definitions [3]-(2.1.9), [3]-(2.1.10)

$$u = u_\alpha e^\alpha$$

Therefore, the set $\overline{\overline{e}} = (e^\alpha \in V, \alpha \in I)$ is basis of $F$-vector space $V'$. □

COROLLARY 4.2.3. *Let $\overline{\overline{e}} = (e_\alpha \in V, \alpha \in I)$ be basis of $F$-vector space $V$. Then*

$$a^\alpha = (e^\alpha, a)$$

□

Basis $\overline{\overline{e}} = (e^\alpha \in V, \alpha \in I)$ is called **basis dual to basis** $\overline{\overline{e}} = (e_\alpha \in V, \alpha \in I)$.

Since the product in field is commutative, then there is no difference between left and right vector spaces. If we consider a vector space over division ring, then we has difference between left and right vector spaces. In the future, I suppose to explore group of transformations in vector space over division ring. So in this book, I want to use notation which is compatible with non commutative product.

According to the theorem [1]-4.4.3, since we consider right vector space, then we write a matrix of map on the left of coordinates of vector. In particular, this is true for matrix of linear transformation and for coordinates of linear functional. According to the theorem [1]-6.2.2, the set of linear functionals on right vector space is left vector space. According to the theorem [1]-4.4.3, we write a matrix of map on the right of coordinates of linear functional. To preserve functional notation, we write a map of conjugated vector space

$$[A] \circ u = u_\alpha A^\alpha_\beta e^\beta$$

DEFINITION 4.2.4. *Let $\overline{\overline{e}} = (e_\alpha \in V, \alpha \in I)$ be basis of $F$-vector space $V$. Let $\overline{\overline{e}} = (e^\alpha \in V, \alpha \in I)$ be basis dual to the basis $\overline{\overline{e}} = (e_\alpha \in V, \alpha \in I)$. Let*

$$(4.2.5) \qquad v_1'^\alpha = f_{1 \cdot \beta}^{\phantom{1 \cdot}\alpha}(a) v_1^\beta$$

*be a representation of a Lie group $G$ in a continuous $F$-vector space $V$. A representation*

$$v_{2 \cdot \alpha}' = v_{2 \cdot \beta} f_{2 \cdot \alpha}^{\phantom{2 \cdot \alpha}\beta}(a)$$

*of a Lie group $G$ in a continuous $F$-vector space $V'$ is called **representation conjugated to representation** (4.2.5), if*

$$(4.2.6) \qquad ([f_2(a)] \circ v_2, f_1(a) \circ v_1) = (u, v)$$

for any $a \in G$. □

THEOREM 4.2.5. *Let $\bar{\bar{e}} = (e_\alpha \in V, \alpha \in I)$ be basis of $F$-vector space $V$. Let $\bar{\bar{e}} = (e^\alpha \in V, \alpha \in I)$ be basis dual to the basis $\bar{\bar{e}} = (e_\alpha \in V, \alpha \in I)$. If a representation*

$$v_1'^\alpha = f_{1 \cdot \beta}^{\ \ \alpha}(a) v_1^\beta \tag{4.2.7}$$

*of a Lie group $G$ in a continuous $F$-vector space $V$ has conjugated representation*

$$v_{2 \cdot \alpha}' = v_{2 \cdot \beta} f_{2 \cdot \alpha}^{\ \ \beta}(a) \tag{4.2.8}$$

*of a Lie group $G$ in a continuous $F$-vector space $V'$, then*

$$f_2(a) = f_1^{-1}(a) = f_1(a^{-1}) \tag{4.2.9}$$

PROOF. The equation

$$f_1(a) \circ v_1 = e_\nu f_{1 \cdot \beta}^{\ \ \nu}(a) v_1^\beta \tag{4.2.10}$$

follows from the equation (4.2.7). The equation

$$[f_2(a)] \circ v_2 = v_{2 \cdot \alpha} f_{2 \cdot \lambda}^{\ \ \alpha}(a) e^\lambda \tag{4.2.11}$$

follows from the equation (4.2.8). The equation

$$\begin{aligned}
v_2 \circ v_1 = v_{2 \cdot \alpha} v_1^\alpha &= ([f_2(a)] \circ v_2) \circ (f_1(a) \circ v_1) \\
&= (v_{2 \cdot \alpha} f_{2 \cdot \mu}^{\ \ \alpha}(a) e^\mu) \circ (e_\nu f_{1 \cdot \beta}^{\ \ \nu}(a) v_1^\beta) \\
&= v_{2 \cdot \alpha} f_{2 \cdot \mu}^{\ \ \alpha}(a)(e^\mu \circ e_\nu) f_{1 \cdot \beta}^{\ \ \nu}(a) v_1^\beta \\
&= v_{2 \cdot \alpha} f_{2 \cdot \mu}^{\ \ \alpha}(a) \delta_\nu^\mu f_{1 \cdot \beta}^{\ \ \nu}(a) v_1^\beta \\
&= v_{2 \cdot \alpha} f_{2 \cdot \nu}^{\ \ \alpha}(a) f_{1 \cdot \beta}^{\ \ \nu}(a) v_1^\beta
\end{aligned} \tag{4.2.12}$$

follows from equations (4.2.6), (4.2.10), (4.2.11). Since $v$ is arbitrary vector, then the equation

$$v_{2 \cdot \beta} = v_{2 \cdot \alpha} f_{2 \cdot \nu}^{\ \ \alpha}(a) f_{1 \cdot \beta}^{\ \ \nu}(a) \tag{4.2.13}$$

follows from the equation (4.2.12). Since $v_2$ is arbitrary functional, then the equation

$$\delta_\beta^\alpha = f_{2 \cdot \nu}^{\ \ \alpha}(a) f_{1 \cdot \beta}^{\ \ \nu}(a) \tag{4.2.14}$$

follows from the equation (4.2.13). The equation (4.2.9) follows from the equation (4.2.14).

Since the map

$$a \to f_1(a) \tag{4.2.15}$$

is representation of group $G$, then the equation

$$\begin{aligned}
v_{2 \cdot \beta} f_{2 \cdot \alpha}^{\ \ \beta}(a_1) f_{2 \cdot \nu}^{\ \ \alpha}(a_2) &= v_{2 \cdot \beta} f_{1 \cdot \alpha}^{\ \ \beta}(a_1^{-1}) f_{1 \cdot \nu}^{\ \ \alpha}(a_2^{-1}) \\
&= v_{2 \cdot \beta} f_{1 \cdot \nu}^{\ \ \alpha}(a_2^{-1}) f_{1 \cdot \alpha}^{\ \ \beta}(a_1^{-1}) \\
&= v_{2 \cdot \beta} f_{1 \cdot \nu}^{\ \ \beta}(a_2^{-1} a_1^{-1}) = v_{2 \cdot \beta} f_{1 \cdot \nu}^{\ \ \beta}((a_1 a_2)^{-1}) \\
&= v_{2 \cdot \beta} f_{2 \cdot \nu}^{\ \ \beta}(a_1 a_2)
\end{aligned} \tag{4.2.16}$$

follows from equations (4.1.5), (4.2.9). From the equation (4.2.16), it follows that the map

$$a \to f_2(a)$$

is right side representation of group $G$. □

THEOREM 4.2.6. *Let*

$$v_1'^\alpha = f_{1.\beta}^{\ \alpha}(a) v_1^\beta \tag{4.2.17}$$

*be linear representation of group $G$ in $F$-vector space $V$. The representation*

$$v_{2\cdot\alpha}' = v_{2\cdot\beta} f_{2.\alpha}^{\ \beta}(a) \tag{4.2.18}$$

*conjugated to representation* (4.2.17) *satisfies to system of differential equations*

$$\frac{\partial v_{2\cdot\beta}'}{\partial a^L} = v_{2\cdot\alpha}' I_{2.\beta K}^{\ \alpha} \lambda_l{}^K_L(a) \tag{4.2.19}$$

*where infinitesimal generators of representation* (4.2.18) *are defined by equation*

$$I_{2.\beta L}^{\ \alpha} = -I_{1.\beta L}^{\ \alpha} \tag{4.2.20}$$

PROOF. The equation

$$I_{2.\beta M}^{\ i} = \left.\frac{\partial f_{2.\beta}^{\ \alpha}(a)}{\partial a^M}\right|_{a=e} = \left.\frac{\partial f_{1.\beta}^{\ \alpha}(a^{-1})}{\partial a^M}\right|_{a=e} = \left.\frac{\partial f_{1.\beta}^{\ \alpha}(a^{-1})}{\partial a^{-1\cdot N}} \frac{\partial a^{-1\cdot N}}{\partial a^M}\right|_{a=e} \tag{4.2.21}$$

follows from equations (4.1.12), (4.2.9). It is evident that

$$\left.\frac{\partial f_{1.\beta}^{\ \alpha}(a^{-1})}{\partial a^{-1\cdot N}}\right|_{a=e} = \left.\frac{\partial f_{1.\beta}^{\ \alpha}(a)}{\partial a^N}\right|_{a=e} = I_{1.\beta N}^{\ \alpha} \tag{4.2.22}$$

From the equation (3.1.33), it follows that

$$\left.\frac{\partial a^{-1\cdot N}}{\partial a^M}\right|_{a=e} = -\psi_{l\cdot K}^{\ N}(a^{-1}) \lambda_{r\cdot M}^{\ K}(a)\big|_{a=e} = -\delta^N_M \tag{4.2.23}$$

The equation (4.2.20) follows from equations (4.2.21), (4.2.22), (4.2.23).

Let $v_1$ be geometric object of representation (4.2.17). Let $v_2$ be geometric object of representation (4.2.18). Consider expression

$$v = v_{2\cdot\alpha} v_1^\alpha \tag{4.2.24}$$

The equation

$$v' = v_{2\cdot\alpha}' v_1'^\alpha = v_{2\cdot\beta} f_{2.\alpha}^{\ \beta}(a) f_{1.\nu}^{\ \alpha}(a) v_1^\nu = v_{2\cdot\beta} v_1^\beta = v \tag{4.2.25}$$

follows from equations (4.2.14), (4.2.17), (4.2.18), (4.2.24). The equation

$$\frac{\partial v'}{\partial a^L} = 0 \tag{4.2.26}$$

follows from the equation (4.2.25). The equation

$$\frac{\partial v'}{\partial a^L} = \frac{\partial v_{2\cdot\alpha}' v_1'^\alpha}{\partial a^L} = \frac{\partial v_{2\cdot\alpha}'}{\partial a^L} v_1'^\alpha + v_{2\cdot\alpha}' \frac{\partial v_1'^\alpha}{\partial a^L} \tag{4.2.27}$$

follows from the equation (4.2.24). The equation

$$\frac{\partial v_{2\cdot\alpha}'}{\partial a^L} v_1'^\alpha + v_{2\cdot\alpha}' \frac{\partial v_1'^\alpha}{\partial a^L} = 0 \tag{4.2.28}$$

follows from equations (4.2.26), (4.2.27). The equation

$$\frac{\partial v_{2\cdot\alpha}'}{\partial a^L} v_1'^\alpha + v_{2\cdot\alpha}' I_{1.\beta K}^{\ \alpha} v_1'^\beta \lambda_l{}^K_L(a) = 0 \tag{4.2.29}$$

follows from equations (4.1.11), (4.2.28). Since vector $v_1$ is arbitrary vector, then the equation

$$\frac{\partial v'_{2 \cdot \beta}}{\partial a^L} = -v'_{2 \cdot \alpha} I_{1 \cdot \beta K}^{\ \ \alpha} \lambda_l{}^K_L(a) \tag{4.2.30}$$

follows from the equation (4.2.29). The equation (4.2.19) follows from equations (4.2.20), (4.2.30). $\square$

### 4.3. Algebraic Properties of Linear Representation

THEOREM 4.3.1. *Let*

$$\begin{aligned} v'^{\alpha}_1 &= f^{\alpha}_{\beta}(a) v^j_1 \\ v'^{\alpha}_2 &= f^{\alpha}_{\beta}(a) v^j_2 \end{aligned} \tag{4.3.1}$$

*be linear geometric objects. Then their sum is also linear geometric object*

$$v'^{\alpha} = f^{\alpha}_{\beta}(a) v^{\beta} \tag{4.3.2}$$

PROOF. From the equation (4.3.1), it follows that

$$v'\alpha = v'^{\alpha}_1 + v'^{\alpha}_2 = f^{\alpha}_{\beta}(a) v^{\beta}_1 + f^{\alpha}_{\beta}(a) v^{\beta}_2 \tag{4.3.3}$$

The equation (4.3.2) follows from the equation (4.3.3). $\square$

THEOREM 4.3.2. *Let*

$$v'^{\alpha}_1 = f_{1 \cdot \beta}{}^{\alpha}(a) v^{\beta}_1$$

*be the representation of group $G$ in $F$-vector space $V_1$. Let*

$$v'^{\gamma}_2 = f_{2 \cdot \epsilon}{}^{\gamma}(a) v^{\epsilon}_2$$

*be the representation of group $G$ in $F$-vector space $V_2$. Let geometric object*

$$v'^{\alpha \gamma} = f_{1 \cdot \beta}{}^{\alpha}(a) f_{2 \cdot \epsilon}{}^{\gamma}(a) v^{\beta \epsilon}$$

*be tensor product of geometric objects $v_1 \in V_1$ and $v_2 \in V_2$*

$$f(a)(v_1 \otimes v_2) = f_1(a)(v_1) \otimes f_2(a)(v_2) \tag{4.3.4}$$

*The representation $f$ is called* **tensor product of representations** *and*

$$I = I_1 \otimes E_2 + E_1 \otimes I_2 \tag{4.3.5}$$

PROOF. First of all let us show that we get representation of group $G$ in vector space $V = V_1 \otimes V_2$. If we apply transformations for $a \in G$ and $b \in G$ consequently then we have

$$\begin{aligned} f(a)(f(b)(v_1 \otimes v_2)) &= f(a)(f_1(b)(v_1) \otimes f_2(b)(v_2)) \\ &= f_1(a)(f_1(b)(v_1)) \otimes f_2(a)(f_2(b)(v_2)) \\ &= f_1(ab)(v_1) \otimes f_2(ab)(v_2) \\ &= f(ab)(v_1 \otimes v_2) \end{aligned} \tag{4.3.6}$$

Therefore $f$ is also representation.

Differentiating equation

$$f(a) = f_1(a) \otimes f_2(a)$$

with respect to $a$ we get

$$\frac{\partial f}{\partial a} = \frac{\partial f_1}{\partial a} \otimes f_2 + f_1 \otimes \frac{\partial f_2}{\partial a} \tag{4.3.7}$$

We get (4.3.5) from (4.3.7) by definition when $a = e$. □

THEOREM 4.3.3. *Let vector space $V$ be direct sum of vector spaces $V_1$ and $V_2$*

$$V = V_1 \oplus V_2$$

*Let*

$$v_1'^\alpha = f_{1.\beta}^{\ \alpha}(a) v_1^\beta$$

*be the representation of group $G$ in vector space $V_1$. Let*

$$v_2'^\gamma = f_{2.\epsilon}^{\ \gamma}(a) v_2^\epsilon$$

*be the representation of group $G$ in vector space $V_2$. Let geometric object $v$ is direct sum of geometric objects $v_1 \in V_1$ and $v_2 \in V_2$*

(4.3.8) $$f(a) \circ (v_1 \oplus v_2) = (f_1(a) \circ v_1) \oplus (f_2(a) \circ v_2)$$

*Then we can define representation of group $G$ in vector space $V = V_1 \oplus V_2$ and*

$$I_P = \begin{pmatrix} I_{1\,b_1 P}^{\ a_1} & 0 \\ 0 & I_{2\,b_2 P}^{\ a_2} \end{pmatrix}$$

*The representation $f$ is called* **direct sum of representations**.

PROOF. We can just get derivative of mapping (4.3.8). □

## 4.4. Differential Equations of Representation

THEOREM 4.4.1. *The linear representation of Lie group $G$*

(4.4.1) $$v'^\alpha = f_\beta^\alpha(a) v^\beta$$

*satisfies to system of differential equations*

(4.4.2) $$\frac{\partial f_\beta^\alpha}{\partial a^L} = I_{\nu K}^\alpha f_\beta^\nu \lambda_{l\,L}^{\ K}(a)$$

PROOF. According to the theorem 4.1.2, the linear representation (4.4.1) satisfies to system of differential equations

(4.4.3) $$\frac{\partial v'^\alpha}{\partial a^L} = I_{\beta K}^\alpha v'^\beta \lambda_{l\,L}^{\ K}(a)$$

The equation

(4.4.4) $$\frac{\partial f_\beta^\alpha(a) v^\beta}{\partial a^L} = I_{\beta K}^\alpha f_\nu^\beta(a) v^\nu \lambda_{l\,L}^{\ K}(a)$$

follows from the equation (4.4.3). Since $v$ is arbitrary vector, then the equation (4.4.2) follows from the equation (4.4.4). □

THEOREM 4.4.2. *Let*

(4.4.5) $$v_1'^\alpha = f_{1.\beta}^{\ \alpha}(a) v_1^\beta$$

*be linear representation of group $G$ in $F$-vector space $V$. Then*

(4.4.6) $$I_{1.\beta K}^{\ \alpha} f_{1.\gamma}^{\ \beta} \lambda_{l\,L}^{\ K}(a) = f_{1.\rho}^{\ \alpha}(a) I_{1.\gamma J}^{\ \rho} \lambda_{r\,L}^{\ J}(a)$$

PROOF. According to the theorem 4.2.6, there exists linear representation

$$v'_{2 \cdot \alpha} = v_{2 \cdot \beta} f_{2 \cdot \alpha}^{\cdot \beta}(a) \tag{4.4.7}$$

conjugated to representation (4.4.5) and conjugated representation satisfies to system of differential equations

$$\frac{\partial v'_{2 \cdot \beta}}{\partial a^L} = v'_{2 \cdot \alpha} I_{2 \cdot \beta K}^{\cdot \alpha} \lambda_{l\,L}^{K}(a) \tag{4.4.8}$$

The equation

$$\frac{\partial v_{2 \cdot \gamma} f_{2 \cdot \beta}^{\cdot \gamma}(a)}{\partial a^L} = v_{2 \cdot \gamma} f_{2 \cdot \alpha}^{\cdot \gamma}(a) I_{2 \cdot \beta K}^{\cdot \alpha} \lambda_{l\,L}^{K}(a) \tag{4.4.9}$$

follows from the equation (4.4.8). Since $v_2$ is arbitrary vector, then the equation

$$\frac{\partial f_{2 \cdot \beta}^{\cdot \gamma}(a)}{\partial a^L} = f_{2 \cdot \alpha}^{\cdot \gamma}(a) I_{2 \cdot \beta K}^{\cdot \alpha} \lambda_{l\,L}^{K}(a) \tag{4.4.10}$$

follows from the equation (4.4.9). The equation

$$\frac{\partial f_{1 \cdot \beta}^{\cdot \gamma}(a^{-1})}{\partial a^L} = -f_{1 \cdot \alpha}^{\cdot \gamma}(a^{-1}) I_{1 \cdot \beta K}^{\cdot \alpha} \lambda_{l\,L}^{K}(a) \tag{4.4.11}$$

follows from equations (4.2.9), (4.2.20), (4.4.10). The equation

$$\frac{\partial f_{1 \cdot \beta}^{\cdot \gamma}(a^{-1})}{\partial a^L} = \frac{\partial f_{1 \cdot \beta}^{\cdot \gamma}(a^{-1})}{\partial a^{-1 \cdot M}} \frac{\partial a^{-1 \cdot M}}{\partial a^L} = -\frac{\partial f_{1 \cdot \beta}^{\cdot \gamma}(a^{-1})}{\partial a^{-1 \cdot M}} \psi_{r\,N}^{M}(a^{-1}) \lambda_{l\,L}^{N}(a) \tag{4.4.12}$$

follows from the equation (3.1.33). The equation

$$\begin{aligned} -\frac{\partial f_{1 \cdot \beta}^{\cdot \gamma}(a^{-1})}{\partial a^{-1 \cdot M}} &= -f_{1 \cdot \alpha}^{\cdot \gamma}(a^{-1}) I_{1 \cdot \beta K}^{\cdot \alpha} \lambda_{l\,L}^{K}(a) \psi_{l\,N}^{L}(a) \lambda_{r\,M}^{N}(a^{-1}) \\ &= -f_{1 \cdot \alpha}^{\cdot \gamma}(a^{-1}) I_{1 \cdot \beta K}^{\cdot \alpha} \lambda_{r\,M}^{K}(a^{-1}) \end{aligned} \tag{4.4.13}$$

follows from equations (3.1.18), (3.1.19), (4.4.11), (4.4.12). The equation

$$\frac{\partial f_{1 \cdot \beta}^{\cdot \gamma}(a)}{\partial a^M} = f_{1 \cdot \alpha}^{\cdot \gamma}(a) I_{1 \cdot \beta K}^{\cdot \alpha} \lambda_{r\,M}^{K}(a) \tag{4.4.14}$$

follows from the equation (4.4.13). According to the theorem 4.4.1, linear representation (4.4.5) satisfies to system of differential equations

$$\frac{\partial f_{1 \cdot \beta}^{\cdot \alpha}}{\partial a^L} = I_{\nu K}^{\alpha} f_{1 \cdot \beta}^{\cdot \nu} \lambda_{l\,L}^{K}(a) \tag{4.4.15}$$

The equation (4.4.6) follows from equations (4.4.13), (4.4.15). □

THEOREM 4.4.3. *Infinitesimal generators of linear representation are constants.*

PROOF. According to the theorem 4.1.5, passive transformation of left shift

$$a_1 \to g a_1$$

on Lie group $G$ generates transformation of infinitesimal generators

$$I_{2 \cdot \mu K}^{\cdot \alpha} = f_{\mu}^{\beta}(g) f^{-1 \cdot \alpha}_{\rho}(g) I_{1 \cdot \beta J}^{\cdot \rho} \lambda_{l\,M}^{J}(g) \psi_{r\,K}^{M}(g) \tag{4.4.16}$$

The equation

$$I_{1 \cdot \gamma J}^{\cdot \rho} = f^{-1 \cdot \rho}_{\alpha}(a) I_{1 \cdot \beta K}^{\cdot \alpha} f_{\gamma}^{\beta} \lambda_{l\,L}^{K}(a) \psi_{r\,J}^{L}(a) \tag{4.4.17}$$

follows from equations (3.1.18), (3.1.19), (4.4.6). The theorem follows from equations (4.4.16), (4.4.17). □

THEOREM 4.4.4. *The system of differential equations* (4.4.2) *is completely integrable system iff*

(4.4.18) $$I^\alpha_{\beta K} I^\beta_{\nu P} - I^\alpha_{\beta P} I^\beta_{\nu K} = I^\alpha_{\nu T} C^T_{lPK}$$

PROOF. Equation

(4.4.19) $$\frac{\partial^2 f^\alpha_\beta}{\partial a^L \partial a^M} = I^\alpha_{\nu K} \frac{\partial f^\nu_\beta}{\partial a^M} \lambda^K_{lL}(a) + I^\alpha_{\nu K} f^\nu_\beta \frac{\partial \lambda^K_{lL}(a)}{\partial a^M}$$
$$= I^\alpha_{\nu K} I^\nu_{\mu P} f^\mu_\beta \lambda^P_{lM}(a) \lambda^K_{lL}(a) + I^\alpha_{\nu K} f^\nu_\beta \frac{\partial \lambda^K_{lL}(a)}{\partial a^M}$$

follows from the equation (4.4.2). The system of differential equations (4.4.2) is completely integrable system iff

(4.4.20) $$\frac{\partial^2 f^\alpha_\beta}{\partial a^L \partial a^M} = \frac{\partial^2 f^\alpha_\beta}{\partial a^M \partial a^L}$$

The equation

(4.4.21) $$I^\alpha_{\nu K} I^\nu_{\mu P} f^\mu_\beta \lambda^P_{lM}(a) \lambda^K_{lL}(a) + I^\alpha_{\nu K} f^\nu_\beta \frac{\partial \lambda^K_{lL}(a)}{\partial a^M}$$
$$= I^\alpha_{\nu P} I^\nu_{\mu K} f^\mu_\beta \lambda^P_{lM}(a) \lambda^K_{lL}(a) + I^\alpha_{\nu K} f^\nu_\beta \frac{\partial \lambda^K_{lM}(a)}{\partial a^L}$$

follows from equations (4.4.19), (4.4.20). The equation

(4.4.22) $$(I^\alpha_{\nu K} I^\nu_{\mu P} - I^\alpha_{\nu P} I^\nu_{\mu K}) f^\mu_\beta \lambda^P_{lM}(a) \lambda^K_{lL}(a)$$
$$= I^\alpha_{\nu K} f^\nu_\beta \left( \frac{\partial \lambda^K_{lM}(a)}{\partial a^L} - \frac{\partial \lambda^K_{lL}(a)}{\partial a^M} \right)$$
$$= I^\alpha_{\nu K} f^\nu_\beta C^K_{lPK} \lambda^P_{lM}(a) \lambda^K_{lL}(a)$$

follows from equations (3.4.3), (4.4.21). Since the matrix $\lambda_l$ is non singular matrix and the map $f$ is arbitrary, then the equation (4.4.18) follows from the equation (4.4.22). □

## 4.5. Group $Gl_n(F)$

Consider the representation of group $GL_n(F)$ in congugated vector space

(4.5.1) $$u'_k = u_l a^{-1 \cdot l}_{\phantom{-1 \cdot} k}$$

To write the system of differential equations for representation (4.5.1), we repeat development of the system of differential equations (4.1.11) considering the representation (4.5.1).

The equation (4.1.15) has form

(4.5.2) $$0 = \frac{\partial u'_b b^{-1 \cdot b}_{\phantom{-1 \cdot} n}}{\partial b^k_l} \frac{\partial b^k_l}{\partial a^m_q} + \frac{\partial u'_b b^{-1 \cdot b}_{\phantom{-1 \cdot} n}}{\partial u'_c} \frac{\partial u'_c}{\partial a^m_q}$$
$$= u'_b \frac{\partial b^{-1 \cdot b}_{\phantom{-1 \cdot} n}}{\partial b^k_l} \frac{\partial b^k_l}{\partial a^m_q} + \frac{\partial u'_b}{\partial u'_c} b^{-1 \cdot b}_{\phantom{-1 \cdot} n} \frac{\partial u'_c}{\partial a^m_q}$$

The equation

$$
\begin{aligned}
0 & = \frac{\partial b_l^k b^{-1\cdot l}{}_m}{\partial b_r^p} = \frac{\partial b_l^k}{\partial b_r^p} b^{-1\cdot l}{}_m + b_l^k \frac{\partial b^{-1\cdot l}{}_m}{\partial b_r^p} = \delta_p^k \delta_l^r b^{-1\cdot l}{}_m + b_l^k \frac{\partial b^{-1\cdot l}{}_m}{\partial b_r^p} \\
& = \delta_p^k b^{-1\cdot r}{}_m + b_l^k \frac{\partial b^{-1\cdot l}{}_m}{\partial b_r^p}
\end{aligned}
\tag{4.5.3}
$$

follows from the equation $b_l^k b^{-1\cdot l}{}_m = \delta_m^k$. The equation

$$
\frac{\partial b^{-1\cdot l}{}_m}{\partial b_r^p} = -b^{-1\cdot l}{}_k \delta_p^k b^{-1\cdot r}{}_m = -b^{-1\cdot l}{}_p b^{-1\cdot r}{}_m
\tag{4.5.4}
$$

follows from the equation (4.5.3). The equation

$$
\frac{\partial b_l^k}{\partial a_n^m} = \frac{\partial c_p^k a^{-1\cdot p}{}_l}{\partial a_n^m} = c_p^k \frac{\partial a^{-1\cdot p}{}_l}{\partial a_n^m} = -c_p^k a^{-1\cdot p}{}_m a^{-1\cdot n}{}_l
\tag{4.5.5}
$$

follows from the equation $b = ca^{-1}$ and the equation (4.5.3).

The equation

$$
\begin{aligned}
0 & = u'_b b^{-1\cdot b}{}_k b^{-1\cdot l}{}_n c_p^k a^{-1\cdot p}{}_m a^{-1\cdot q}{}_l + \delta_b^c b^{-1\cdot b}{}_n \frac{\partial u'_c}{\partial a_q^m} \\
& = u'_b b^{-1\cdot b}{}_k b^{-1\cdot l}{}_n b_m^k a^{-1\cdot q}{}_l + b^{-1\cdot c}{}_n \frac{\partial u'_c}{\partial a_q^m} \\
& = u'_m b^{-1\cdot l}{}_n a^{-1\cdot q}{}_l + b^{-1\cdot c}{}_n \frac{\partial u'_c}{\partial a_q^m}
\end{aligned}
\tag{4.5.6}
$$

follows from equations (4.5.2), (4.5.4), (4.5.5). The equation

$$
\frac{\partial u'_c}{\partial a_q^m} = -u'_m a^{-1\cdot q}{}_c
\tag{4.5.7}
$$

follows from the equation (4.5.6).

The equation

$$
I_{\cdot-m\cdot p}^{l\cdot-r} = -\delta_p^l \delta_m^r
\tag{4.5.8}
$$

follows from the equation (4.5.4).

CHAPTER 5

# CHAPTER 6

# Index

# CHAPTER 7

# Special Symbols and Notations

# Линейное представление группы Ли

Александр Клейн

Aleks_Kleyn@MailAPS.org
`http://AleksKleyn.dyndns-home.com:4080/`
`http://sites.google.com/site/AleksKleyn/`
`http://arxiv.org/a/kleyn_a_1`
`http://AleksKleyn.blogspot.com/`

Аннотация. В книге рассмотрены дифференциальные уравнения для операции в группе Ли, а также для представлений группы Ли в векторном пространстве.

# Оглавление

Глава 1

# Введение

## 1.1. Об этой книге

Эта книга основана на дипломном проекте, написанном мною в одесском университете. Мне очень повезло, что мой преподаватель Гаврильченко Михаил Леонидович прочёл курс лекций "Группы и алгебры Ли". Эти лекции были основаны на книге [5]. Эти лекции стали фундаментом моих исследований на долгие годы. Опираясь на эти лекции, я записал уравнения группы Ли для конкретного случая представления группы $GL$ в векторном пространстве геометрических объектов.

Впоследствии группа $GL$ была заменена на произвольную группу Ли. Именно в таком виде эта теория была включена в книгу. Соответственно этому главы книги имеют следующее содержание. Представление группы Ли тесно связано с теорией дифференциальных уравнений. В главе 2 дано краткое изложение теории дифференциальных уравнений, посвящённое вполне интегрируемым системам дифференциальных уравнений. В главе 3 я записываю дифференциальные уравнения для группы Ли. Дифференциальные уравнения представления группы Ли записаны в главе 4.

## 1.2. Соглашения

Соглашение 1.2.1. *Пусть $A$ - $\Omega_1$-алгебра. Пусть $B$ - $\Omega_2$-алгебра. Запись*

$$A \longrightarrow_* \rightarrow B$$

*означает, что определено представление $\Omega_1$-алгебры $A$ в $\Omega_2$-алгебре $B$.* □

Без сомнения, у читателя могут быть вопросы, замечания, возражения. Я буду признателен любому отзыву.

Глава 2

# Дифференциальные уравнения

## 2.1. Вполне интегрируемые системы

Соглашение 2.1.1. *Пусть $F$ - полное поле. В этой книге, мы рассматриваем $F$-векторные пространства. Соответственно, мы рассматриваем дифференциальные уравнения над полем $F$.* □

Пусть[2.1] нам дана система $nm$ дифференциальных уравнений в частных производных

$$(2.1.1) \qquad \frac{\partial \theta^\alpha}{\partial x^i} = \psi_i^\alpha(\theta^1, ..., \theta^m, x^1, ..., x^n) = \psi_i^\alpha(\theta, x)$$

$$1 \leq \alpha \leq m \quad 1 \leq i \leq n \quad \theta^\alpha = \theta^\alpha(x^1, ..., x^n)$$

Эта система разрешена относительно всех частных производных и эквивалентна системе уравнений в полных дифференциалах

$$d\theta^\alpha = \psi_i^\alpha(\theta, x) dx^i$$

Чтобы найти условия интегрируемости, мы дифференцируем $(2.1.1)$ по $x^j$.

$$\frac{\partial^2 \theta^\alpha}{\partial x^j \partial x^i} = \frac{\partial \psi_i^\alpha}{\partial x^j} + \frac{\partial \psi_i^\alpha}{\partial \theta^\sigma} \frac{\partial \theta^\sigma}{\partial x^j} =$$
$$= \frac{\partial \psi_i^\alpha}{\partial x^j} + \frac{\partial \psi_i^\alpha}{\partial \theta^\sigma} \psi_j^\sigma$$
$$\frac{\partial^2 \theta^\alpha}{\partial x^i \partial x^j} = \frac{\partial \psi_j^\alpha}{\partial x^i} + \frac{\partial \psi_j^\alpha}{\partial \theta^\sigma} \psi_i^\sigma$$

Так как вторая производная непрерывной функции не зависит от порядка дифференцирования, мы ожидаем, что

$$(2.1.2) \qquad \frac{\partial \psi_i^\alpha}{\partial x^j} + \frac{\partial \psi_i^\alpha}{\partial \theta^\sigma} \psi_j^\sigma = \frac{\partial \psi_j^\alpha}{\partial x^i} + \frac{\partial \psi_j^\alpha}{\partial \theta^\sigma} \psi_i^\sigma$$

Определение 2.1.2. Мы называем систему дифференциальных уравнений $(2.1.1)$ **вполне интегрируемой системой**, если условие $(2.1.2)$ удовлетворено тождественно. □

Решение вполне интегрируемой системы $(2.1.1)$ можно разложить в ряд Тейлора

$$(2.1.3) \quad \theta^\alpha = C^\alpha + \left(\frac{\partial \theta^\alpha}{\partial x^i}\right)_{x=x_0} (x^i - x_0^i) + \frac{1}{2} \left(\frac{\partial^2 \theta^\alpha}{\partial x^j \partial x^i}\right)_{x=x_0} (x^i - x_0^i)(x^j - x_0^j) + ...$$

---

[2.1]Эта глава написана под влиянием [5].

Здесь $C^\alpha$ ($1 \le \alpha \le m$) - константы,
$$\left(\frac{\partial \theta^\alpha}{\partial x^i}\right)_{x=x_0} = \psi_i^\alpha(C^1,...,C^m, x_0^1,...,x_0^n)$$

Мы получим остальные коэффициенты дифференцированием правой части (2.1.1) и подстановкой $C^\alpha$ и $x_0^i$. Мы говорим о решении, когда ряд (2.1.3) сходится. Согласно теореме Коши и Ковалевской это справедливо, если $\psi_i^\alpha$ аналитичны по всем аргументам.

В окрестности точки $(x_0^1,...,x_0^n)$, для которой ряд (2.1.3) сходится, мы имеем решение

(2.1.4) $$\theta^\alpha = \varphi^\alpha(x^1,...,x^n, C^1,...,C^m)$$

определённое $m$ константами.

Если условие (2.1.2) не удовлетворенно тождественно, то система дифференциальных уравнений (2.1.1) не является вполне интегрируемой. Следовательно, функции $\theta^\alpha$ удовлетворяют (2.1.1) и (2.1.2). Если ограничения (2.1.2) независимы от системы $S_1$, то уравнение (2.1.2) накладывает ограничения на функции $\theta^\alpha$ и порождает систему $S_1$ уравнений. Мы дифференцируем уравнения системы $S_1$ по $x^k$ и заменим производные $\dfrac{\partial \theta^\alpha}{\partial x^i}$ их значениями из уравнений (2.1.1). Если новые уравнения не удовлетворены тождественно, мы получим новую систему $S_2$ уравнений. Мы можем повторить этот процесс снова.

Так как мы не можем иметь более чем $m$ ограничений, эта цепь обрывается. Если система дифференциальных уравнений (2.1.1) имеет решение, это решение удовлетворяет всем ограничениям, в противном случае эти ограничения противоречивы.

Теорема 2.1.3. *Необходимое и достаточное условие существования решения системы дифференциальных уравнений* (2.1.1) *- это существование такого $N$, что $S_1,...,S_N$ алгебраически совместимы, и $S_{N+1}$ является их следствием. Если мы получили $p$ ограничений, то общее решение зависит от $m - p$ произвольных констант.*

Доказательство. Мы доказали, что это условие необходимо. Теперь мы покажем, что это условие достаточно. Предположим, что мы получили $p$ ограничений

(2.1.5) $$\Phi^\gamma(\theta, x) = 0 \quad \gamma = 1,...,p$$

Так как
$$Rg\left\|\frac{\partial \Phi^\gamma}{\partial \theta^\alpha}\right\| = p$$

мы можем выразить $p$ переменных $\theta^\alpha$ через остальные, пользуясь (2.1.5). Если мы правильно перенумеруем переменные, мы получим

(2.1.6) $$\theta^\nu = \phi^\nu(\theta^{p+1},...,\theta^m, x) \quad \nu = 1,...,p$$

Подставим (2.1.6) в (2.1.1) и прлучим

(2.1.7) $$\frac{\partial \theta^\alpha}{\partial x^i} = \bar\psi_i^\alpha(\theta^{p+1},...,\theta^m, x^1,...,x^n) \quad \alpha = p+1,...,m$$

Условие (2.1.2) удовлетворенно. Система дифференциальных уравнений (2.1.7) вполне интегрируема и имеет решение, зависящее от $m - p$ произвольных констант. Первые $p$ уравнений удовлетворены в силу (2.1.6), когда $\theta^{p+1}, ..., \theta^m, x$ являются решениями (2.1.7). □

## 2.2. Линейный дифференциальный оператор

Пусть даны $r$ дифференциальных операторов

$$(2.2.1) \qquad X_a f = \xi_a^i \frac{\partial f}{\partial x^i}$$

Оператор (2.2.1) линеен

$$X_a(f + g) = X_a f + X_a g$$
$$X_a(\lambda f) = \lambda X_a f \quad \lambda = const$$

Произведение этих операторов имеет вид

$$X_a X_b f = X_a(X_b f) = \xi_a^i \frac{\partial}{\partial x^i}\left(\xi_b^j \frac{\partial f}{\partial x^j}\right)$$
$$= \xi_a^i \frac{\partial \xi_b^j}{\partial x^i}\frac{\partial f}{\partial x^j} + \xi_a^i \xi_b^j \frac{\partial^2 f}{\partial x^i \partial x^j} =$$
$$= \frac{\partial f}{\partial x^j} X_a \xi_b^j + \xi_a^i \xi_b^j \frac{\partial^2 f}{\partial x^i \partial x^j}$$

Мы использовали (2.2.1) на последнем шаге. Это не оператор вида (2.2.1). Однако так как второе слагаемое симметрично относительно $a$ и $b$, мы можем определить коммутатор

$$(2.2.2) \qquad (X_a, X_b)f = X_a X_b f - X_b X_a f = (X_a \xi_b^j - X_b \xi_a^j)\frac{\partial f}{\partial x^j}$$

Это линейная операция и мы зовём её скобкой Пуассона. Косая симметрия этого коммутатора следует из (2.2.2).

Простое вычисление показывает, что

$$(2.2.3) \qquad ((X_a, X_b), X_c)f = ((X_a, X_b)\xi_c^j - X_c(X_a \xi_b^j - X_b \xi_a^j))\frac{\partial f}{\partial x^j} =$$

$$= \left(\frac{\partial \xi_c^j}{\partial x^k} X_a \xi_b^k - \frac{\partial \xi_c^j}{\partial x^k} X_b \xi_a^k - \frac{\partial \xi_b^j}{\partial x^k} X_c \xi_a^k + \xi_c^p \xi_a^q \frac{\partial^2 \xi_b^j}{\partial x^p \partial x^q} + \frac{\partial \xi_a^j}{\partial x^k} X_c \xi_b^k + \xi_c^p \xi_b^q \frac{\partial^2 \xi_a^j}{\partial x^p \partial x^q}\right)\frac{\partial f}{\partial x^j}$$

Если мы изменим порядок параметров, мы окончательно получим уравнение

$$(2.2.4) \qquad ((X_a, X_b), X_c)f + ((X_b, X_c), X_a)f + ((X_c, X_a), X_b)f = 0$$

(2.2.4) - это условие Якоби.

Теорема 2.2.1. *Коммутатор линейной комбинации операторов*

$$X'_a = \lambda_a^b X_b$$

*где $\lambda_a^b$ - функции $x$, выражается через эти операторы и их коммутаторы*

$$(2.2.5) \qquad (X'_a, X'_c)f = \mu_{ac}^d X_d f - \mu_{ca}^d X_d f + \lambda_a^b \lambda_c^d (X_b, X_d)f$$

Доказательство. Коммутатор имеет вид

$$(X'_a, X'_c)f = (\lambda_a^b X_b, \lambda_c^d X_d)f$$
$$(2.2.6) \quad = \lambda_a^b X_b(\lambda_c^d X_d f) - \lambda_c^d X_d(\lambda_a^b X_b f)$$
$$= \lambda_a^b X_b \lambda_c^d X_d f + \lambda_a^b \lambda_c^d X_b X_d f - \lambda_c^d X_d \lambda_a^b X_b f - \lambda_c^d \lambda_a^b X_d X_b f$$

(2.2.5) следует из (2.2.6), если положить

$$\mu_{ac}^d = \lambda_a^b X_b \lambda_c^d$$

□

## 2.3. Полная система линейных дифференциальных уравнений в частных производных

Теперь мы хотим изучить систему дифференциальных уравнений

$$(2.3.1) \quad X_a f = 0 \quad X_a = \xi_a^i \frac{\partial}{\partial x^i} \quad 1 \le i \le n \quad 1 \le a \le r$$

Мы положим, что

$$(2.3.2) \quad \operatorname{rank} \|\xi_a^i\| = r \le n$$

Это значит, что уравнения (2.3.1) линейно независимы.

Если $r = n$, то, очевидно, единственное решение будет $f = \text{const}$. Утверждение, что каждое решение уравнений (2.3.1) является также решением уравнения

$$(2.3.3) \quad (X_a, X_b)f = 0$$

следует из (2.2.2). Если

$$(2.3.4) \quad (X_a, X_b)f = \gamma_{ab}^c X_c f$$

где $\gamma$ - функции $x$, то система дифференциальных уравнений (2.3.1) и (2.3.3) эквивалентна системе дифференциальных уравнений (2.3.1). Мы не можем знать заранее, накладывает ли (2.3.3) новые условия на (2.3.1). Присоединим к системе дифференциальных уравнений (2.3.1) те уравнения системы (2.3.3), для которых не справедливо (2.3.4). Получим $s$ уравнений, $s \ge r$. Для новой системы повторим тот же процесс, получим $t$ уравнений, $t \ge s$, и т. д..

В результате последовательности указанных действий мы можем получить $n$ уравнений. В этом случае система дифференциальных уравнений (2.3.1) имеет единственное решение $f = \text{const}$, так как

$$\frac{\partial f}{\partial x^i} = 0$$

В противном случае мы получим $u$ уравнений, $u \le n$, для которых выполняется условие (2.3.4). В этом случае полученная система дифференциальных уравнений называется **полной** порядка $u$. Решение системы дифференциальных уравнений (2.3.1) является решением полной системы дифференциальных уравнений.

Теорема 2.3.1. *Если система дифференциальных уравнений* (2.3.1) *полна, то полна также система дифференциальных уравнений*

$$(2.3.5) \quad X'_a f = \lambda_a^b X_b f = 0$$

*где $\lambda$ функции $x$*

$$(2.3.6) \qquad \det\|\lambda_a^b\| \neq 0 \quad 1 \leq a, b \leq r$$

Доказательство. Из (2.3.4) следует, что правая сторона (2.2.5) линейна относительно $X_a f$, которые в силу (2.3.5), (2.3.6) линейно выражается через операторы $X'_a f$. Следовательно, из (2.2.5) получим выражения вида (2.3.4), что и доказывает утверждение. $\square$

В силу (2.3.2), не нарушая общности, положим

$$(2.3.7) \qquad \det\|\xi_a^i\| \neq 0 \quad 1 \leq i, a \leq r$$

Следовательно, мы можем разрешить систему дифференциальных уравнений (2.3.1) относительно $\frac{\partial f}{\partial x^1}$, ..., $\frac{\partial f}{\partial x^r}$ и записать результат в виде

$$(2.3.8) \qquad X'_a f = \frac{\partial f}{\partial x^a} + \psi_a^t \frac{\partial f}{\partial x^t} = 0 \quad a = 1, ..., r \quad t = r+1, ..., n$$

Эти уравнения имеют вид (2.3.5) и, следовательно, образуют систему дифференциальных уравнений, эквивалентную системе (2.3.1).

Полная система дифференциальных уравнений, представленная в виде (2.3.8), называется **якобиевой**.

Теорема 2.3.2. *Полная система дифференциальных уравнений (2.3.1) имеет точно $n - r$ независимых решений.*

Доказательство. Аналогично равенству (2.3.4) мы имеем

$$(2.3.9) \qquad (X'_a, X'_b)f = \gamma'^c_{ab} X'_c f$$

Из (2.2.2) и (2.3.8), следует, что

$$(2.3.10) \qquad (X'_a, X'_b)f = (X_a \psi_b^t - X_b \psi_a^t)\frac{\partial f}{\partial x^t} \quad t = r+1, ..., n$$

Так как $\frac{\partial f}{\partial x^1}$, ..., $\frac{\partial f}{\partial x^r}$ не входят в (2.3.10), $\gamma'^c_{ab} = 0$. Следовательно, для полной системы в якобиевой форме

$$(X'_a, X'_b)f = 0$$
$$X_a \psi_b^t - X_b \psi_a^t = 0$$

$$(2.3.11) \qquad \frac{\partial \psi_b^p}{\partial x^a} + \psi_a^q \frac{\partial \psi_b^p}{\partial x^q} = \frac{\partial \psi_a^p}{\partial x^b} + \psi_b^q \frac{\partial \psi_a^p}{\partial x^q} \quad a, b = 1, ..., r \quad p, q = r+1, ..., n$$

Сравнивая (2.3.11) и (2.1.2), мы видим, что система дифференциальных уравнений

$$(2.3.12) \qquad \frac{\partial x^p}{\partial x^a} = \psi_a^p \quad a = 1, ..., r \quad p = r+1, ..., n$$

вполне интегрируемы. Согласно теореме 2.1.3 система дифференциальных уравнений (2.3.8), а следовательно, система дифференциальных уравнений (2.3.1) допускают $n - r$ независимых решений, и более $n - r$ независимых решений система дифференциальных уравнений (2.3.1) иметь не может. $\square$

## 2.4. Существенные параметры семейства функций

Определение 2.4.1. Пусть даны $n$ функций
$$f^i(x^1, ..., x^n, a^1, ..., a^r), i = 1, ..., n$$
от переменных $x^1, ..., x^n$ и от параметров $a^1, ..., a^r$. Мы предполагаем, что функции $f^i$ непрерывны по $x^1, ..., x^n$ и $a^1, ..., a^r$. Предполагается также, что непрерывны их производные любого порядка. Мы называем параметры $a^1, ..., a^r$ **существенными параметрами**, если не существует функций $A^1, ..., A^{r-1}$, зависящих только от $a^1, ..., a^r$, таких, что имеют место тождества
$$(2.4.1) \qquad f^i(x, a) = F^i(x, A)$$

$\square$

Теорема 2.4.2. *Для того, чтобы $r$ параметров $a^\alpha$ были существенны, необходимо и достаточно, чтобы функции $f^i$ не удовлетворяли никакому уравнению вида*
$$(2.4.2) \qquad \phi^\alpha \frac{\partial f^i}{\partial a^\alpha} = 0$$
*где $\phi^\alpha \neq 0$.*

Доказательство. Предположим, что параметры $a$ не являются существенными. Тогда существуют $A$ такие, что (2.4.1) удовлетворено и
$$\operatorname{rank} \left\| \frac{\partial A^\sigma}{\partial a^\alpha} \right\| < r$$
Следовательно, существуют функции $\phi^\alpha(a) \neq 0$ такие, что
$$\phi^\alpha \frac{\partial A^\sigma}{\partial a^\alpha} = 0$$
Таким образом, функции $A^1, ..., A^{r-1}$ и любая их функция $\Phi(A)$ также удовлетворяют (2.4.2). Функции $F^i$ в (2.4.1) являются примерами функций $\Phi$. Следовательно, если $f^i$ не зависят существенно от $a^\alpha$, они удовлетворяют системе (2.4.2).

Обратно, положим, что функции $f^i$ удовлетворяют уравнению (2.4.1) для некоторых функций $\phi^\alpha$. Это уравнение имеет $r - 1$ независимых решения $A^1$, ..., $A^{r-1}$, которые являются функциями одних $a^\alpha$ и любое решение системы (2.4.2) является функцией $A^1, ..., A^{r-1}$. Следовательно, каждая из функций $f^i$ является функцией $x^1, ..., x^n$ и $A^1, ..., A^{r-1}$ и $a^1, ..., a^r$ не являются существенными. $\square$

Замечание 2.4.3. Из теоремы 2.4.2 следует, что если функции $f^i$ удовлетворяют полной системе $s$ дифференциальных уравнений
$$(2.4.3) \qquad \phi^\alpha_\sigma \frac{\partial f^i}{\partial a^\alpha} = 0 \quad \alpha = 1, ..., r \quad \sigma = 1, ..., s < r$$
то $f^i$ являются функциями $x$ и $r - s$ независимых решений системы (2.4.3), являющихся функциями одних $a$, и $f^i$ выражаются, следовательно, через $r - s$ существенных параметров. $\square$

Мы можем интерпретировать систему уравнений (2.4.2) как систему линейных уравнений относительно $\phi^\alpha$. Из линейной алгебры следует, что решения системы (2.4.3) порождают векторное пространство, и только при условии

$$\mu_0 = \operatorname{rank}\left(\frac{\partial f^i}{\partial a^\alpha}\right) = r \tag{2.4.4}$$

система линейных уравнений (2.4.3) имеет единственное решение $\phi^\alpha = 0$. Таким образом, из замечания 2.4.3 следует, что параметры существенны, если выполнено условие (2.4.4).

Допустим $\mu_0 < r$. Дифференцируя систему дифференциальных уравнений (2.4.2) по $x^j$, мы получим

$$\phi^\alpha \frac{\partial f^i}{\partial a^\alpha} = 0 \quad \phi^\alpha \frac{\partial f^i_{,j}}{\partial a^\alpha} = 0 \quad f^i_{,j} = \frac{\partial f^i}{\partial x^j} \tag{2.4.5}$$

Пусть

$$\mu_1 = \operatorname{rank}\left(\frac{\partial f^i}{\partial a^\alpha} \quad \frac{\partial f^i_{,j}}{\partial a^\alpha}\right) \tag{2.4.6}$$

Очевидно, что $\mu_1 \geq \mu_0$. Если $\mu_1 = r$, то система линейных уравнений (2.4.5) допускает единственное решение $\phi^\alpha = 0$. Следовательно, параметры существенны.

Полагая

$$f^i_{,j_1...j_s} = \frac{\partial^s f^i}{\partial x^{j_1}...\partial x^{j_s}}$$

последовательными дифференцированиями (2.4.2) получим

$$\phi^\alpha \frac{\partial f^i_{,j_1...j_s}}{\partial a^\alpha} = 0 \tag{2.4.7}$$

Обозначим

$$\mu_s = \operatorname{rank}\left(\frac{\partial f^i}{\partial a^\alpha} \quad \frac{\partial f^i_{,j}}{\partial a^\alpha} \quad ... \quad \frac{\partial f^i_{,j_1...j_s}}{\partial a^\alpha}\right) \tag{2.4.8}$$

Таким образом мы получим последовательность рангов

$$\mu_0 \leq \mu_1 \leq ... \leq \mu_s \leq ... \leq r \tag{2.4.9}$$

Теорема 2.4.4. *Число существенных параметров, через которые выражаются функции $f^i(x, a)$, равно максимальному числу, содержащемуся в последовательности* (2.4.9).

Доказательство. Если какое-нибудь $\mu_s$ равно $r$, то $\phi^\alpha = 0$ и $r$ параметров существенны.

Предположим, $\mu_{s-1} < r$ и $\mu_s = \mu_{s-1}$. Из этого утверждения следует

$$\frac{\partial f^i_{,j_1...j_s}}{\partial a^\alpha} = \sum_{t<s} \lambda^{k_1...k_t}_{j_1...j_s} \frac{\partial f^i_{,k_1...k_t}}{\partial a^\alpha} \tag{2.4.10}$$

Дифференцируя равенство (2.4.10) по $x^p$, мы получим

$$\frac{\partial f^i_{,j_1\ldots j_s p}}{\partial a^\alpha} = \sum_{t<s} \lambda^{k_1\ldots k_t}_{j_1\ldots j_s} \frac{\partial f^i_{,k_1\ldots k_t p}}{\partial a^\alpha}$$

$$= \sum_{t<s-1} \lambda^{k_1\ldots k_t}_{j_1\ldots j_s} \frac{\partial f^i_{,k_1\ldots k_t p}}{\partial a^\alpha} + \sum_{t=s-1, u<s} \lambda^{k_1\ldots k_t}_{j_1\ldots j_s} \lambda^{l_1\ldots l_u}_{k_1\ldots k_t p} \frac{\partial f^i_{,l_1\ldots l_u}}{\partial a^\alpha}$$

Следовательно, $\mu_{s+1} = \mu_{s-1}$.

Так как $\mu_s < r$, и система линейных уравнений (2.4.2), (2.4.7) имеет ранг $\mu_s$, то $r - \mu_s$ из функций $\phi^\alpha$ могут быть выбраны произвольно, причём остальные полностью определяются этим выбором. Возьмём $\phi^1, ..., \phi^{r-\mu_s}$ функциями одних только $a$. Тогда

$$\phi^\sigma = \lambda^\sigma_\rho \phi^\rho \quad \rho = 1, ..., r - \mu_s \quad \sigma = r - \mu_s + 1, ..., r$$

где $\lambda$ - функция от $x$ и $a$. Соответственно этому существует $r - \mu_s$ независимых уравнения (2.4.2). Коммутатор любых двух из этих уравнений, приравненный 0, имеет $f^i$ в качестве решений, поэтому он является линейной комбинацией данных уравнений. Следовательно, эти уравнения образуют полную систему. Таким образом, данные функции выражаются через $\mu_s$ существенных параметров. $\square$

Глава 3

# Группа Ли

### 3.1. Операция на группе Ли

Пусть $G$ является $r$-параметрической группой Ли с нейтральным элементом $e$. Пусть
$$a \to (a^1, ..., a^n)$$
локальные координаты в группе Ли $G$. Пусть

(3.1.1) $\qquad a_3^L = \varphi^L(a_1, a_2) \quad a_3 = a_1 a_2$

операция на группе Ли $G$.

Все операторы, которые мы определим ниже связаны с представлениями левого и правого сдвига и будут иметь дополнительный индекс $r$ или $l$, говорящий нам, какой тип сдвига они описывают.

Сперва мы введём операторы

(3.1.2) $\qquad A_r{}^K_L(a,b) = \dfrac{\partial \varphi^K(a,b)}{\partial b^L}$

(3.1.3) $\qquad A_l{}^K_L(a,b) = \dfrac{\partial \varphi^K(a,b)}{\partial a^L}$

как производная левого и правого сдвигов.

Теорема 3.1.1.

(3.1.4) $\qquad A_r(a,bc)A_r(b,c) = A_r(ab,c)$

(3.1.5) $\qquad A_l(ab,c)A_l(a,b) = A_l(a,bc)$

Доказательство. Операция ассоциативна a(bc)=(ab)c Пользуясь правилом дифференцирования сложной функции, мы можем найти производные этого равенства

(3.1.6) $\qquad \dfrac{\partial a(bc)}{\partial bc}\dfrac{\partial bc}{\partial c} = \dfrac{\partial (ab)c}{\partial c}$

(3.1.7) $\qquad \dfrac{\partial a(bc)}{\partial a} = \dfrac{\partial (ab)c}{\partial ab}\dfrac{\partial ab}{\partial a}$

(3.1.4) следует из (3.1.6). (3.1.5) следует из (3.1.7). $\qquad \square$

Так как $ae = a$ и $eb = b$, мы получаем равенства

(3.1.8) $\qquad A_l{}^K_L(a,e) = \delta^K_L$

(3.1.9) $\qquad A_r{}^K_L(e,b) = \delta^K_L$

Теорема 3.1.2. *Оператор $A_l(a,b)$ имеет обратный оператор $A_l(ab, b^{-1})$*

(3.1.10) $\qquad A_l^{-1}(a,b) = A_l(ab, b^{-1})$

*Оператор $A_r(b,c)$ имеет обратный оператор $A_r(b^{-1}, bc)$*

$$(3.1.11) \qquad A_r^{-1}(b,c) = A_r(b^{-1}, bc)$$

Доказательство. Согласно (3.1.4) мы имеем $(a = b^{-1})$

$$(3.1.12) \qquad A_r(b^{-1}, bc) A_r(b, c) = A_r(b^{-1}b, c) = A_r(e, c) = \delta$$

(3.1.11) следует из (3.1.12).

Согласно (3.1.5) мы имеем $(c = b^{-1})$

$$(3.1.13) \qquad A_l(ab, b^{-1}) A_l(a, b) = A_l(a, bb^{-1}) = A_l(a, e) = \delta$$

(3.1.10) следует из (3.1.13). $\square$

Теорема 3.1.3. *Операторы $A_l(a,b)$ и $A_r(a,b)$ обратимы.*

Доказательство. Это следствие теоремы 3.1.2. $\square$

Мы определим **базовые операторы группы Ли**

$$(3.1.14) \qquad \psi_r{}^L_N(a) = A_r{}^L_N(a, e) \qquad\qquad \psi_r(a) = A_r(a, e)$$
$$(3.1.15) \qquad \psi_l{}^L_N(b) = A_l{}^L_N(e, b) \qquad\qquad \psi_l(b) = A_l(e, b)$$

Из (3.1.8) и (3.1.9) немедленно следует, что

$$(3.1.16) \qquad \psi_l{}^K_L(e) = \delta^K_L$$
$$(3.1.17) \qquad \psi_r{}^K_L(e) = \delta^K_L$$

По определению базовые операторы линейно отображают касательную плоскость $T_eG$ в касательную плоскость $T_aG$.

Теорема 3.1.4. *Операторы $\psi_l$ и $\psi_r$ обратимы.*

Доказательство. Это следствие теоремы 3.1.3 и определений (3.1.15) и (3.1.14). $\square$

Так как операторы $\psi_r$ и $\psi_l$ имеют обратные операторы мы определим операторы

$$(3.1.18) \qquad \lambda_r(a) = \psi_r^{-1}(a)$$
$$(3.1.19) \qquad \lambda_l(a) = \psi_l^{-1}(a)$$

которые отображают $T_aG \to T_eG$

Теорема 3.1.5.

$$(3.1.20) \qquad \lambda_r(a) = A_r(a^{-1}, a)$$
$$(3.1.21) \qquad \lambda_l(a) = A_l(a, a^{-1})$$

Доказательство. Согласно (3.1.19), (3.1.15) и (3.1.10) мы получим

$$\lambda_l(a) = A_l^{-1}(e, a) = A_l(a, a^{-1})$$

Это доказывает (3.1.21).

Согласно (3.1.18), (3.1.14) и (3.1.11) мы получим

$$\lambda_r(a) = A_r^{-1}(a, e) = A_r(a^{-1}, a)$$

Это доказывает (3.1.20). $\square$

Теорема 3.1.6.

(3.1.22) $$A_r(a,b) = \psi_r(ab)\lambda_r(b)$$

(3.1.23) $$A_l(a,b) = \psi_l(ab)\lambda_l(a)$$

Доказательство. Положим $c = e$ в (3.1.4)

(3.1.24) $$A_r(a,b)A_r(b,e) = A_r(ab,e)$$

Из равенств (3.1.14), (3.1.24) следует

(3.1.25) $$A_r(a,b)\psi_r(b) = \psi_r(ab)$$

Из (3.1.25) следует

(3.1.26) $$A_r(a,b) = \psi_r(ab)\psi_r^{-1}(b)$$

(3.1.22) следует из (3.1.26) на основании (3.1.18).

Положим $a = e$ в (3.1.5)

(3.1.27) $$A_l(b,c)A_l(e,b) = A_l(e,bc)$$

Из равенств (3.1.15), (3.1.27) следует

(3.1.28) $$A_l(b,c)\psi_l(b) = \psi_l(bc)$$

Из (3.1.28) следует

(3.1.29) $$A_l(a,b) = \psi_l(ab)\psi_l^{-1}(a)$$

(3.1.23) следует из (3.1.29) на основании (3.1.19). $\square$

Теорема 3.1.7. *Операция группы Ли удовлетворяет дифференциальному уравнению*

(3.1.30) $$\frac{\partial \varphi^K(a,b)}{\partial b^L} = \psi_r{}^K_T(ab)\lambda_r{}^T_L(b) \qquad \frac{\partial ab}{\partial b} = \psi_r(ab)\lambda_r(b)$$

(3.1.31) $$\frac{\partial \varphi^K(a,b)}{\partial a^L} = \psi_l{}^K_T(ab)\lambda_l{}^T_L(a) \qquad \frac{\partial ab}{\partial a} = \psi_l(ab)\lambda_l(a)$$

Доказательство. Равенство (3.1.30) является следствием равенств (3.1.22) и (3.1.2). Равенство (3.1.31) является следствием равенств (3.1.23) и (3.1.3). $\square$

Теорема 3.1.8.

(3.1.32) $$\frac{\partial a^{-1}}{\partial a} = -\psi_l(a^{-1})\lambda_r(a)$$

(3.1.33) $$\frac{\partial a^{-1}}{\partial a} = -\psi_r(a^{-1})\lambda_l(a)$$

Доказательство. Дифференцируя равенство $e = a^{-1}a$ по $a$, мы получим равенство

$$0 = \frac{\partial a^{-1}a}{\partial a^{-1}}\frac{\partial a^{-1}}{\partial a} + \frac{\partial a^{-1}a}{\partial a} =$$
$$= A_l(a^{-1},a)\frac{\partial a^{-1}}{\partial a} + A_r(a^{-1},a)$$
$$\frac{\partial a^{-1}}{\partial a} = -A_l^{-1}(a^{-1},a)A_r(a^{-1},a)$$

Пользуясь (3.1.20) и (3.1.21), мы получим

$$\frac{\partial a^{-1}}{\partial a} = -\lambda_l^{-1}(a^{-1})\lambda_r(a) \tag{3.1.34}$$

(3.1.32) следует из (3.1.34).

Дифференцируя равенство $e = aa^{-1}$ по $a$, мы получим равенство

$$0 = \frac{\partial aa^{-1}}{\partial a} + \frac{\partial aa^{-1}}{\partial a^{-1}}\frac{\partial a^{-1}}{\partial a} =$$
$$= A_l(a, a^{-1}) + A_r(a, a^{-1})\frac{\partial a^{-1}}{\partial a}$$
$$\frac{\partial a^{-1}}{\partial a} = -A_r^{-1}(a, a^{-1})A_l(a, a^{-1})$$

Пользуясь (3.1.20) и (3.1.21) мы получим

$$\frac{\partial a^{-1}}{\partial a} = -\lambda_r^{-1}(a^{-1})\lambda_l(a) \tag{3.1.35}$$

(3.1.33) следует из (3.1.35). $\square$

Теорема 3.1.9.

$$\frac{\partial a^{-1}b}{\partial a} = -\psi_l(a^{-1}b)\psi_r^{-1}(a) \tag{3.1.36}$$

$$\frac{\partial ba^{-1}}{\partial a} = -\psi_r(ba^{-1})\psi_l^{-1}(a) \tag{3.1.37}$$

Доказательство. Пользуясь правилом дифференцирования сложной функции и равенствами (3.1.23), (3.1.32), мы получим

$$\begin{aligned}\frac{\partial a^{-1}b}{\partial a} &= \frac{\partial a^{-1}b}{\partial a^{-1}}\frac{\partial a^{-1}}{\partial a} \\ &= -A_l(a^{-1}, b)\psi_l(a^{-1})\lambda_r(a) \\ &= -\psi_l(a^{-1}b)\psi_l^{-1}(a^{-1})\psi_l(a^{-1})\lambda_r(a)\end{aligned} \tag{3.1.38}$$

(3.1.36) следует из (3.1.38).

Пользуясь правилом дифференцирования сложной функции и равенствами (3.1.22), (3.1.33), мы получим

$$\begin{aligned}\frac{\partial ba^{-1}}{\partial a} &= \frac{\partial ba^{-1}}{\partial a^{-1}}\frac{\partial a^{-1}}{\partial a} \\ &= -A_r(b, a^{-1})\psi_r(a^{-1})\lambda_l(a) \\ &= -\psi_r(ba^{-1})\psi_r^{-1}(a^{-1})\psi_r(a^{-1})\lambda_l(a)\end{aligned} \tag{3.1.39}$$

(3.1.37) следует из (3.1.39). $\square$

Теорема 3.1.10.

$$\frac{\partial abc}{\partial b} = \psi_l(abc)\lambda_l(ab)\psi_r(ab)\lambda_r(b) \tag{3.1.40}$$

$$\frac{\partial abc}{\partial b} = \psi_r(abc)\lambda_r(bc)\psi_l(bc)\lambda_l(b) \tag{3.1.41}$$

Доказательство. Пользуясь правилом дифференцирования сложной функции и равенствами (3.1.30), (3.1.31) мы получим

$$\frac{\partial abc}{\partial b} = \frac{\partial abc}{\partial ab}\frac{\partial ab}{\partial b} = \psi_l(abc)\lambda_l(ab)\psi_r(ab)\lambda_r(b) \tag{3.1.42}$$

Равенство (3.1.40) является следствием равенства (3.1.42). Пользуясь правилом дифференцирования сложной функции и равенствами (3.1.30), (3.1.31) мы получим

$$\frac{\partial abc}{\partial b} = \frac{\partial abc}{\partial bc}\frac{\partial bc}{\partial b} = \psi_r(abc)\lambda_r(bc)\psi_l(bc)\lambda_l(b) \tag{3.1.43}$$

Равенство (3.1.41) является следствием равенства (3.1.43). □

Теорема 3.1.11.

$$\frac{\partial aba^{-1}}{\partial a} = (\psi_l(aba^{-1}) - \psi_r(aba^{-1}))\lambda_l(a) \tag{3.1.44}$$

$$\frac{\partial aba^{-1}}{\partial b} = \psi_l(aba^{-1})\lambda_l(ab)\psi_r(ab)\lambda_r(b) \tag{3.1.45}$$

$$\frac{\partial aba^{-1}}{\partial b} = \psi_r(aba^{-1})\lambda_r(ba^{-1})\psi_l(ba^{-1})\lambda_l(b) \tag{3.1.46}$$

Доказательство. Равенство

$$\begin{aligned}\frac{\partial aba^{-1}}{\partial a} &= \frac{\partial a(ba^{-1})}{\partial a} + \frac{\partial (ab)a^{-1}}{\partial a^{-1}}\frac{\partial a^{-1}}{\partial a} \\ &= \psi_l(aba^{-1})\lambda_l(a) - \psi_r(aba^{-1})\lambda_r(a^{-1})\psi_r(a^{-1})\lambda_l(a)\end{aligned} \tag{3.1.47}$$

является следствием равенств (3.1.30), (3.1.31), (3.1.33). Равенство

$$\frac{\partial aba^{-1}}{\partial a} = \psi_l(aba^{-1})\lambda_l(a) - \psi_r(aba^{-1})\lambda_l(a) \tag{3.1.48}$$

является следствием равенств (3.1.47), (3.1.18). Равенство (3.1.44) является следствием равенства (3.1.48).

Равенство (3.1.45) является следствием равенства (3.1.40). Равенство (3.1.46) является следствием равенства (3.1.41). □

Теорема 3.1.12.

$$\left.\frac{\partial aba^{-1}}{\partial b}\right|_{b=e} = \lambda_l(a)\psi_r(a) \tag{3.1.49}$$

$$\left.\frac{\partial aba^{-1}}{\partial b}\right|_{b=e} = \lambda_r(a^{-1})\psi_l(a^{-1}) \tag{3.1.50}$$

$$\lambda_l(a)\psi_r(a) = \lambda_r(a^{-1})\psi_l(a^{-1}) \tag{3.1.51}$$

Доказательство. Равенство (3.1.49) является следствием равенства (3.1.45). Равенство (3.1.50) является следствием равенства (3.1.46). Равенство (3.1.51) является следствием равенств (3.1.49), (3.1.50). □

### 3.2. 1-параметрическая группа

Если многообразие группы имеет размерность 1, мы имеем группу, зависящую от 1 параметра. В этом случае операция на группе имеет вид

$$c = \varphi(a, b)$$

где $a, b, c$ - числа. В этом случае наша запись будет упрощена. Нам даны не операторы, а функции

$$A_l(a, b) = \frac{\partial \varphi(a, b)}{\partial a}$$
$$A_r(a, b) = \frac{\partial \varphi(a, b)}{\partial b}$$

которые по определению удовлетворяют уравнению

$$A_l(a, e) = 1$$
$$A_r(e, b) = 1$$

Мы также определим базовую функцию

$$\psi_r(a) = A_r(a, e)$$
$$\psi_l(b) = A_l(e, b)$$

Функции $\psi_r(a)$ и $\psi_l(a)$ никогда не обращаются в 0 и

$$\lambda_r(a) = \frac{1}{\psi_r(a)}$$
$$\lambda_l(a) = \frac{1}{\psi_l(a)}$$

Теорема 3.2.1. *Операция на 1-параметрической группе Ли удовлетворяет дифференциальным уравнениям*

$$(3.2.1) \qquad \frac{\partial \varphi(a, b)}{\partial b} = \frac{\psi_r(ab)}{\psi_r(b)}$$

Доказательство. Это следствие равества (3.1.30). □

Теорема 3.2.2. *Мы можем определить координату $A$ на 1-параметрической группе таким образом, что операция $\Phi$ на группе имеет вид*

$$(3.2.2) \qquad \Phi(A, B) = A + B$$

*и $E = 0$ - единица группы.*

Доказательство. Мы определим новую переменную $A$ таким образом, что

$$(3.2.3) \qquad dA = \frac{da}{\psi_r(a)}$$

Так как $\psi_r(a) \neq 0$, существуют взаимнооднозначное отображение $A = F(a)$ и его обратное отображение $a = f(A)$. Это значит, что если $A = F(a)$, $B = F(b)$, $C = F(c)$ и $c = \varphi(a, b)$, то существует функция $\Phi$ такая, что

$$C = \Phi(A, B) = F(\varphi(a, b))$$

Мы получили производную

$$\frac{\partial \Phi(A, B)}{\partial B} = \frac{dC}{dc} \frac{\partial \varphi(a, b)}{\partial b} \frac{db}{dB}$$

Используя (3.2.3) и (3.2.1), мы получим
$$\frac{\partial \Phi(A,B)}{\partial B} = \frac{1}{\psi_r(c)} \frac{\psi_r(c)}{\psi_r(b)} \psi_r(b) = 1$$

Таким образом, мы получили

(3.2.4) $$\Phi(A,B) = \xi(A) + B$$

Если мы возьмём решение уравнения (3.2.3) в виде
$$A = \int_e^a \frac{da}{\psi_r(a)}$$

то мы увидим, что $A = 0$ - единица группы. Если мы положим $B = 0$ в (3.2.4), мы получим

(3.2.5) $$A = \xi(A)$$

(3.2.2) следует из (3.2.4) и (3.2.5) $\square$

### 3.3. Правый сдвиг

Для правого сдвига
$$b' = ba$$
система дифференциальных уравнений (3.1.30) принимает вид

(3.3.1) $$\frac{\partial b'^K}{\partial a^L} = \psi_r{}^K_T(b') \lambda_r{}^T_L(a)$$

Отображения $b'^K$ являются решением системы (3.3.1) и согласно [2]-(2.4.14) они зависят от $b^1, ..., b^n$, которые мы можем предположить постоянными. Таким образом, решение системы (3.3.1) зависит от $n$ произвольных констант и, следовательно, система (3.3.1) вполне интегрируема. Условие её интегрируемости имеет вид

$$\frac{\partial \psi_r{}^L_T(b')}{\partial b'^S} \frac{\partial b'^S}{\partial a^R} \lambda_r{}^T_P(a) + \psi_r{}^L_T(b') \frac{\partial \lambda_r{}^T_P(a)}{\partial a^R} =$$
$$= \frac{\partial \psi_r{}^L_T(b')}{\partial b'^S} \frac{\partial b'^S}{\partial a^P} \lambda_r{}^T_R(a) + \psi_r{}^L_T(b') \frac{\partial \lambda_r{}^T_R(a)}{\partial a^P}$$

Согласно (3.3.1), мы можем записать это условие в более простой форме

$$\frac{\partial \psi_r{}^L_T(b')}{\partial b'^S} \psi_r{}^S_V(b') \lambda_r{}^V_R(a) \lambda_r{}^T_P(a) + \psi_r{}^L_T(b') \frac{\partial \lambda_r{}^T_P(a)}{\partial a^R}$$
$$= \frac{\partial \psi_r{}^L_T(b')}{\partial b'^S} \psi_r{}^S_V(b') \lambda_r{}^V_P(a) \lambda_r{}^T_R(a) + \psi_r{}^L_T(b') \frac{\partial \lambda_r{}^T_R(a)}{\partial a^P}$$

(3.3.2)
$$\frac{\partial \psi_r{}^L_T(b')}{\partial b'^S} \psi_r{}^S_V(b') - \frac{\partial \psi_r{}^L_V(b')}{\partial b'^S} \psi_r{}^S_T(b')$$
$$= \psi_r{}^L_U(b') \psi_r{}^P_T(a) \psi_r{}^R_V(a) \left( \frac{\partial \lambda_r{}^U_R(a)}{\partial a^P} - \frac{\partial \lambda_r{}^U_P(a)}{\partial a^R} \right)$$

Положим

(3.3.3) $$C_r{}^U_{VT} = \psi_r{}^R_V(a) \psi_r{}^P_T(a) \left( \frac{\partial \lambda_r{}^U_R(a)}{\partial a^P} - \frac{\partial \lambda_r{}^U_P(a)}{\partial a^R} \right)$$

Уравнение

(3.3.4) $$\frac{\partial \psi_r{}^L_T(b')}{\partial b'^S} \psi_r{}^S_V(b') - \frac{\partial \psi_r{}^L_V(b')}{\partial b'^S} \psi_r{}^S_T(b') = C_r{}^U_{VT} \psi_r{}^L_U(b')$$

следует из уравнений (3.3.2), (3.3.3). Если мы продифференцируем это равенство по $a^P$, мы получим

$$\frac{\partial C_{rVT}^{\phantom{rVT}U}}{\partial a^P}\psi_{rU}^{\phantom{rU}L}(b') = 0$$

так как $\psi_{rU}^{\phantom{rU}L}(b')$ не зависит от $a$. В тоже время $\psi_{rU}^{\phantom{rU}L}(b')$ линейно независимы, так как $det\|\psi_{rU}^{\phantom{rU}L}\| \neq 0$. Следовательно,

$$\frac{\partial C_{rVT}^{\phantom{rVT}U}}{\partial a^P} = 0$$

и $C_{rTV}^{\phantom{rTV}U}$ являются константами, которые называются **правыми структурными константами алгебры Ли**. Из (3.3.3) следует, что

$$(3.3.5) \qquad C_{rTV}^{\phantom{rTV}U}\lambda_{rP}^{\phantom{rP}T}(a)\lambda_{rR}^{\phantom{rR}V}(a) = \frac{\partial \lambda_{rP}^{\phantom{rP}U}(a)}{\partial a^R} - \frac{\partial \lambda_{rR}^{\phantom{rR}U}(a)}{\partial a^P}$$

Мы называем (3.3.5) уравнением Маурера.

Теорема 3.3.1. *Векторные поля, определённые дифференциальным оператором*

$$(3.3.6) \qquad X_{rV} = \psi_{rV}^{\phantom{rV}S}(a)\frac{\partial}{\partial a^S}$$

*линейно независимы и их коммутатор имеет вид*

$$(X_{rT}, X_{rV}) = C_{rTV}^{\phantom{rTV}U} X_{rU}$$

Доказательство. Линейная независимость векторных полей следует из теоремы 3.1.4. Тогда мы видим, что согласно (3.3.4)

$$(X_{rT}, X_{rV}) = (X_{rT}\psi_{rV}^{\phantom{rV}D} - X_{rV}\psi_{rT}^{\phantom{rT}D})\frac{\partial}{\partial a^D} =$$
$$= \left(\psi_{rT}^{\phantom{rT}P}(a)\frac{\partial \psi_{rV}^{\phantom{rV}D}(a)}{\partial a^P} - \psi_{rV}^{\phantom{rV}R}(a)\frac{\partial \psi_{rT}^{\phantom{rT}D}(a)}{\partial a^R}\right)\frac{\partial}{\partial a^D} =$$
$$= C_{rTV}^{\phantom{rTV}U}\psi_{rU}^{\phantom{rU}D}(a)\frac{\partial}{\partial a^D} = C_{rTV}^{\phantom{rTV}U} X_{rU}$$

$\square$

Пусть задан гомоморфизм $f: G_1 \to G$ 1-параметрической группы Ли $G_1$ в группу $G$. Образ этой группы является 1-параметрической подгруппой. Если $t$ - координата на группе $G_1$, мы можем записать $a = f(t)$ и найти дифференциальное уравнение для этой подгруппы. Мы положим в случае правого сдвига, что $a = f(t_1)$, $b = f(t_2)$, $c = ab = f(t)$, $t = t_1 + t_2$. Тогда мы получим

$$\frac{dc^K}{dt} = \frac{\partial c^K}{\partial b^L}\frac{db^L}{dt} = \psi_{rT}^{\phantom{rT}K}(c)\lambda_{rL}^{\phantom{rL}T}(b)\frac{db^L}{dt_2}\frac{dt_2}{dt}$$
$$\frac{dc^K}{dt} = \psi_{rT}^{\phantom{rT}K}(c)\lambda_{rL}^{\phantom{rL}T}(b)\frac{db^L}{dt_2}$$

Левая часть не зависит от $t_2$, следовательно, правая часть не зависит от $t_2$. Мы положим, что

$$\lambda_{rL}^{\phantom{rL}T}(b)\frac{db^L}{dt_2} = \alpha^T$$

Таким образом, мы получим систему дифференциальных уравнений
$$\frac{dc^K}{dt} = \psi_{r\,T}^{\,K}(c)\alpha^T$$

Так как $\psi_r$ - производная правого сдвига в единице группы, это равенство означает, что 1-параметрическая группа определена вектором $\alpha^T \in T_eG$ и переносит этот вектор вдоль 1-параметрической группы без изменений. Мы называем это векторное поле **правоинвариантным векторным полем**. Мы определим векторное произведение $T_e$ как

(3.3.7) $$[\alpha,\beta]^T = C_r{}_{RS}^{\,T}\alpha^R\beta^S$$

Пространство $T_eG$, снабжённое такой операцией, становится алгеброй Ли $\mathfrak{g}_r$. Мы называем её **определённой справа алгеброй Ли группы Ли**

Теорема 3.3.2. *Пространство правоинвариантных векторных полей имеет конечную размерность, равную размерности группы Ли. Это алгебра Ли с произведением равным коммутатору векторных полей и эта алгебра изоморфна алгебре Ли $\mathfrak{g}_r$.*

Доказательство. Это следует из (3.3.6) и (3.3.7), так как $\alpha^K$ и $\beta^K$ являются константами
$$(X_{rT}\alpha^T, X_{rV}\beta^V) = (X_{rT}, X_{rV})\alpha^T\beta^V =$$
$$= C_r{}_{TV}^{\,U}X_{rU}\alpha^T\beta^V = [\alpha,\beta]^U X_{rU}$$
□

## 3.4. Левый сдвиг

Для левого сдвига
$$b' = ab$$
система дифференциальных уравнений (3.1.31) принимает вид

(3.4.1) $$\frac{\partial b'^K}{\partial a^L} = \psi_{l\,T}^{\,K}(b')\lambda_{l\,L}^{\,T}(a)$$

Функции $b'^K$ являются решением системы (3.4.1) и согласно [2]-(2.4.12) они зависят от $b^1,...,b^n$, которые мы можем предположить постоянными. Таким образом, решение системы (3.4.1) зависит от $n$ произвольных констант и, следовательно, система (3.4.1) вполне интегрируема. Условие её интегрируемости имеет вид
$$\frac{\partial \psi_{l\,T}^{\,L}(b')}{\partial b'^S}\frac{\partial b'^S}{\partial a^R}\lambda_{l\,P}^{\,T}(a) + \psi_{l\,T}^{\,L}(b')\frac{\partial \lambda_{l\,P}^{\,T}(a)}{\partial a^R} =$$
$$= \frac{\partial \psi_{l\,T}^{\,L}(b')}{\partial b'^S}\frac{\partial b'^S}{\partial a^P}\lambda_{l\,R}^{\,T}(a) + \psi_{l\,T}^{\,L}(b')\frac{\partial \lambda_{l\,R}^{\,T}(a)}{\partial a^P}$$

Согласно (3.4.1), мы можем записать это условие в более простой форме
$$\frac{\partial \psi_{l\,T}^{\,L}(b')}{\partial b'^S}\psi_{l\,V}^{\,S}(b')\lambda_{l\,R}^{\,V}(a)\lambda_{l\,P}^{\,T}(a) + \psi_{l\,T}^{\,L}(b')\frac{\partial \lambda_{l\,P}^{\,T}(a)}{\partial a^R}$$
$$= \frac{\partial \psi_{l\,T}^{\,L}(b')}{\partial b'^S}\psi_{l\,V}^{\,S}(b')\lambda_{l\,P}^{\,V}(a)\lambda_{l\,R}^{\,T}(a) + \psi_{l\,T}^{\,L}(b')\frac{\partial \lambda_{l\,R}^{\,T}(a)}{\partial a^P}$$

$$\frac{\partial \psi_{lT}^L(b')}{\partial b'^S}\psi_{lV}^S(b') - \frac{\partial \psi_{lV}^L(b')}{\partial b'^S}\psi_{lT}^S(b')$$

(3.4.2)
$$= \psi_{lU}^L(b')\psi_{lT}^P(a)\psi_{lV}^R(a)\left(\frac{\partial \lambda_{lR}^U(a)}{\partial a^P} - \frac{\partial \lambda_{lP}^U(a)}{\partial a^R}\right)$$

Положим

(3.4.3)
$$C_{lVT}^U = \psi_{lV}^R(a)\psi_{lT}^P(a)\left(\frac{\partial \lambda_{lR}^U(a)}{\partial a^P} - \frac{\partial \lambda_{lP}^U(a)}{\partial a^R}\right)$$

Уравнение

(3.4.4)
$$\frac{\partial \psi_{lT}^L(b')}{\partial b'^S}\psi_{lV}^S(b') - \frac{\partial \psi_{lV}^L(b')}{\partial b'^S}\psi_{lT}^S(b') = C_{lVT}^U\psi_{lU}^L(b')$$

следует из уравнений (3.4.2), (3.4.3). Если мы продифференцируем это равенство по $a^P$, мы получим

$$\frac{\partial C_{lVT}^U}{\partial a^P}\psi_{lU}^L(b') = 0$$

так как $\psi_{lU}^L(b')$ не зависит от $a$. В тоже время $\psi_{lU}^L(b')$ линейно независимы, так как $det\|\psi_{lU}^L\| \neq 0$. Следовательно,

$$\frac{\partial C_{lVT}^U}{\partial a^P} = 0$$

и $C_{lTV}^U$ являются константами, которые называются **левыми структурными константами алгебры Ли**. Из (3.4.3) следует, что

(3.4.5)
$$C_{lTV}^U \lambda_{lP}^T(a)\lambda_{lR}^V(a) = \frac{\partial \lambda_{lP}^U(a)}{\partial a^R} - \frac{\partial \lambda_{lR}^U(a)}{\partial a^P}$$

Мы называем (3.4.5) уравнением Маурера.

Теорема 3.4.1. *Векторные поля, определённые дифференциальным оператором*

(3.4.6)
$$X_{lV} = \psi_{lV}^S(a)\frac{\partial}{\partial a^S}$$

*линейно независимы и коммутатор имеет вид*

$$(X_{lT}, X_{lV}) = C_{lTV}^U X_{lU}$$

Доказательство. Линейная независимость векторных полей следует из теоремы 3.1.4. Тогда мы видим, что согласно (3.4.4)

$$(X_{lT}, X_{lV}) = (X_{lT}\psi_{lV}^D - X_{lV}\psi_{lT}^D)\frac{\partial}{\partial a^D} =$$
$$= \left(\psi_{lT}^P(a)\frac{\partial \psi_{lV}^D(a)}{\partial a^P} - \psi_{lV}^R(a)\frac{\partial \psi_{lT}^D(a)}{\partial a^R}\right)\frac{\partial}{\partial a^D} =$$
$$= C_{lTV}^U\psi_{lU}^D(a)\frac{\partial}{\partial a^D} = C_{lTV}^U X_U$$

□

Пусть задан гомоморфизм $f : G_1 \to G$ 1-параметрической группы Ли $G_1$ в группу $G$. Образ этой группы является 1-параметрической подгруппой. Если $t$

- координата на группе $G_1$, мы можем записать $a = f(t)$ и найти дифференциальное уравнение для этой подгруппы. Мы положим в случае левого сдвига, что $a = f(t_1)$, $b = f(t_2)$, $c = ab = f(t)$, $t = t_1 + t_2$. Тогда мы получим

$$\frac{dc^K}{dt} = \frac{\partial c^K}{\partial b^L}\frac{db^L}{dt} = \psi_l{}_T^K(c)\lambda_l{}_L^T(b)\frac{db^L}{dt_2}\frac{dt_2}{dt}$$

$$\frac{dc^K}{dt} = \psi_l{}_T^K(c)\lambda_l{}_L^T(b)\frac{db^L}{dt_2}$$

Левая часть не зависит от $t_2$, следовательно, правая часть не зависит от $t_2$. Мы положим, что

$$\lambda_l{}_L^T(b)\frac{db^L}{dt_2} = \alpha^T$$

Таким образом, мы получим систему дифференциальных уравнений

$$\frac{dc^K}{dt} = \psi_l{}_T^K(c)\alpha^T$$

Так как $\psi_l$ - производная правого сдвига в единице группы, это равенство означает, что 1-параметрическая группа определена вектором $\alpha^T \in T_eG$ и переносит этот вектор вдоль 1-параметрической группы без изменений. Мы называем это векторное поле **левоинвариантным векторным полем**. Мы определим векторное произведение $T_e$ как

(3.4.7) $$[\alpha, \beta]^T = C_l{}_{RS}^T \alpha^R \beta^S$$

Пространство $T_eG$, снабжённое такой операцией, становится алгеброй Ли $\mathfrak{g}_l$. Мы называем её **определённой слева алгеброй Ли группы Ли**

Теорема 3.4.2. *Пространство правоинвариантных векторных полей имеет конечную размерность, равную размерности группы Ли. Это алгебра Ли с произведением равным коммутатору векторных полей и эта алгебра изоморфна алгебре Ли $\mathfrak{g}_l$.*

Доказательство. Это следует из (3.4.6) и (3.4.7), так как $\alpha^K$ и $\beta^K$ являются константами

$$(X_{lT}\alpha^T, X_{lV}\beta^V) = (X_{lT}, X_{lV})\alpha^T\beta^V =$$
$$= C_l{}_{TV}^U X_{lU}\alpha^T\beta^V = [\alpha, \beta]^U X_{lU}$$

$\square$

## 3.5. Отношение между алгебрами Ли $\mathfrak{g}_l$ и $\mathfrak{g}_r$

Мы определили две различные алгебры Ли на пространстве $T_eG$. Теперь наша задача - это определить отношение между этими алгебрами. Мы начнём решать эту задачу с анализа разложения Тейлора групповой операции.

Теорема 3.5.1. *С точностью до бесконечно малых второго порядка, структура операции на группе Ли*

(3.5.1) $$\varphi^K(a, b) = a^K + b^K - e^K + I_{LM}^K(a^M - e^M)(b^L - e^L)$$

*где мы определили* **инфинитезимальные образующие группы Ли**

$$I_{LM}^K = I_l{}_{LM}^K = I_r{}_{ML}^K$$

Доказательство. Мы можем найти ряд Тейлора произведения по обоим аргументам. Однако подобно разделам 3.3 и 3.4, мы можем рассматривать один аргумент как параметр и найти ряд Тейлора по другому аргументу. В этом случае его коэффициенты будут зависеть от первого аргумента.

Таким образом, согласно (3.1.2) и (3.1.14), произведение в окрестности $a \in G$ имеет ряд Тейлора относительно $b$

$$(3.5.2) \qquad \varphi^K(a,b) = a^K + \psi_{r\,L}^{\,K}(a)(b^L - e^L) + o(b^L - e^L)$$
$$ab = a + \psi_r(a)(b - e) + o(b - e)$$

где коэффициенты разложения зависят от $a$. В тоже время $\psi_r(a)$ также имеет ряд Тейлора в окрестности $e$

$$(3.5.3) \qquad \psi_{r\,L}^{\,K}(a) = \delta_L^K + I_{r\,LM}^{\,K}(a^M - e^M) + o(a^M - e^M)$$

Если мы подставим (3.5.3) в (3.5.2), мы получим разложение

$$\varphi^K(a,b) = a^K + (\delta_L^K + I_{r\,LM}^{\,K}(a^M - e^M) + o(a^M - e^M))(b^L - e^L) + o(b^L - e^L)$$

$$(3.5.4)\ \varphi^K(a,b) = a^K + b^K - e^K + I_{r\,LM}^{\,K}(a^M - e^M))(b^L - e^L) + o(a^M - e^M, b^L - e^L)$$

Согласно (3.1.3) и (3.1.15), произведение в окрестности $b \in G$ имеет ряд Тейлора относительно $a$

$$(3.5.5) \qquad \varphi^K(a,b) = b^K + \psi_{l\,L}^{\,K}(b)(a^L - e^L) + o(a^L - e^L)$$
$$ab = b + \psi_l(b)(a - e) + o(a - e)$$

где коэффициенты разложения зависят от $b$. В тоже время $\psi_l(b)$ также имеет ряд Тейлора в окрестности $e$

$$(3.5.6) \qquad \psi_{l\,L}^{\,K}(b) = \delta_L^K + I_{l\,LM}^{\,K}(b^M - e^M) + o(b^M - e^M)$$

Если мы подставим (3.5.6) в (3.5.5) мы получим разложение

$$\varphi^K(a,b) = b^K + (\delta_L^K + I_{l\,LM}^{\,K}(b^M - e^M) + o(b^M - e^M))(a^L - a^L) + o(a^L - e^L)$$

$$(3.5.7)\ \varphi^K(a,b) = b^K + a^K - e^K + I_{l\,LM}^{\,K}(b^M - e^M))(a^L - e^L) + o(a^M - e^M, b^L - e^L)$$

(3.5.4) и (3.5.7) являются разложением в ряд Тейлора одной и той же функции и они должны совпадать. Сравнивая их, мы видим, что $I_{l\,LM}^{\,K} = I_{r\,ML}^{\,K}$. Тогда операция группы Ли имеет ряд Тейлора (3.5.1). $\square$

Теорема 3.5.2. *С точностью до бесконечно малых второго порядка, структура операции на группе Ли*

$$(3.5.8) \qquad \psi_{r\,L}^{\,K}(a) = \delta_L^K + I_{ML}^{K}(a^M - e^M) + o(a^M - e^M)$$

$$(3.5.9) \qquad \psi_{l\,L}^{\,K}(a) = \delta_L^K + I_{LM}^{K}(a^M - e^M) + o(a^M - e^M)$$

$$(3.5.10) \qquad \lambda_{r\,L}^{\,K}(a) = \delta_L^K - I_{ML}^{K}(a^M - e^M) + o(a^M - e^M)$$

$$(3.5.11) \qquad \lambda_{l\,L}^{\,K}(a) = \delta_L^K - I_{LM}^{K}(a^M - e^M) + o(a^M - e^M)$$

Доказательство. (3.5.8) и (3.5.9) следуют из дифференцирования равенства (3.5.1) по $a$ или $b$ и определений (3.1.14) и (3.1.15).

Теперь мы положим, что $\lambda_r{}_L^K(a)$ имеет ряд Тейлора

(3.5.12) $$\lambda_r{}_L^K(a) = \delta_L^K + J_{LM}^K(a^M - e^M) + o(a^M - e^M)$$

Тогда подставляя (3.5.12) и (3.5.8) в уравнение

$$\lambda_r{}_L^K(a)\psi_r{}_M^L(a) = \delta_M^K$$

мы получим с точностью до порядка 1

$$(\delta_L^K + I_{NL}^K(a^N - e^N) + o(a^N - e^N))(\delta_M^L + J_{MP}^L(a^P - e^P) + o(a^P - e^P)) = \delta_M^K$$
$$\delta_M^K + (J_{MN}^K + I_{NM}^K)(a^N - e^N) + o(a^P - e^P) = \delta_M^K$$

Следовательно, $J_{MN}^K + I_{NM}^K = 0$ и, подставляя $J_{MN}^K$ в (3.5.12), мы получим (3.5.10).

Таким же образом, мы можем доказать (3.5.11). □

Теорема 3.5.3.

(3.5.13) $$\left.\frac{\partial \psi_r{}_L^K(a)}{\partial a^M}\right|_{a=e} = I_{ML}^K$$

(3.5.14) $$\left.\frac{\partial \psi_l{}_L^K(a)}{\partial a^M}\right|_{a=e} = I_{LM}^K$$

(3.5.15) $$\left.\frac{\partial \lambda_r{}_L^K(a)}{\partial a^M}\right|_{a=e} = -I_{ML}^K$$

(3.5.16) $$\left.\frac{\partial \lambda_l{}_L^K(a)}{\partial a^M}\right|_{a=e} = -I_{LM}^K$$

Доказательство. Теорема следует непосредственно из теоремы 3.5.2. □

На основе теоремы 3.5.3 мы можем получить более детальную информацию о левой и правой группах алгебры Ли. Так как определение (3.3.3) не зависит от $a$, мы можем оценить его, когда $a = e$. Согласно (3.1.16) и (3.5.16), мы получим

(3.5.17) $$C_{l\,TV}^U = \delta_T^P \delta_V^R(-I_{PR}^U + I_{RP}^U) = I_{VT}^U - I_{TV}^U$$

Согласно (3.1.17) и (3.5.15), мы получим

(3.5.18) $$C_{r\,TV}^U = \delta_T^P \delta_V^R(-I_{RP}^U + I_{PR}^U) = I_{TV}^U - I_{VT}^U$$

Равенство

$$C_{r\,TV}^U = -C_{l\,TV}^U$$

следует из равенств (3.5.17) и (3.5.18). Следовательно, алгебры $\mathfrak{g}_l$ и $\mathfrak{g}_r$ являются антиизоморфными. Обычно мы подразумеваем, что $\mathfrak{g}_l$ является алгеброй Ли группы Ли и use notation

$$C_{TV}^U = C_{l\,TV}^U$$

Пусть $\mathfrak{g} = T_eG$ является **алгеброй Ли группы Ли** $G$ с операцией

$$\alpha_3^L = C_{MN}^L \alpha_1^M \alpha_2^N \quad \alpha_3 = C(\alpha_1, \alpha_2)$$

### 3.6. Группа $Gl_n(F)$

Если мы рассматриваем векторное пространство размерности $n$ над полем $F$, то множество преобразований порождает группу $GL_n(F)$. Элементы группы $GL_n(F)$ - это матрицы $n \times n$. Произведение матриц $a = (a^i_j)$ и $b = (b^i_j)$ определено равенством

$$(ab)^i_j = a^i_k b^k_j \tag{3.6.1}$$

Теорема 3.6.1.

$$A_r {}^{\cdot k\,i}_{l\cdot j}(a,b) = a^k_j \delta^i_l \tag{3.6.2}$$

$$A_l {}^{\cdot k\,i}_{l\cdot j}(a,b) = b^i_l \delta^k_j \tag{3.6.3}$$

Доказательство. Равенство

$$A_r {}^{\cdot k\,i}_{l\cdot j}(a,b) = \frac{\partial a^k_m b^m_l}{\partial b^j_i} = a^k_m \delta^m_j \delta^i_l \tag{3.6.4}$$

следует из равенств (3.1.2), (3.6.1). Равенство (3.6.2) следует из равенства (3.6.4). Равенство

$$A_l {}^{\cdot k\,i}_{l\cdot j}(a,b) = \frac{\partial a^k_m b^m_l}{\partial a^j_i} = \delta^k_j \delta^i_m b^m_l \tag{3.6.5}$$

следует из равенств (3.1.3), (3.6.1). Равенство (3.6.3) следует из равенства (3.6.5). □

Теорема 3.6.2.

$$\psi_r {}^{\cdot k\,i}_{l\cdot j}(a) = a^k_j \delta^i_l \tag{3.6.6}$$

$$\psi_l {}^{\cdot k\,i}_{l\cdot j}(b) = b^i_l \delta^k_j \tag{3.6.7}$$

Доказательство. Равенство (3.6.6) следует из равенств (3.1.14), (3.6.2). Равенство (3.6.7) следует из равенств (3.1.15), (3.6.3). □

Теорема 3.6.3.

$$\lambda_r {}^{\cdot k\,i}_{l\cdot j}(a) = a^{-1\cdot k}_{\quad\; j} \delta^i_l \tag{3.6.8}$$

$$\lambda_l {}^{\cdot k\,i}_{l\cdot j}(b) = b^{-1\cdot i}_{\quad\; l} \delta^k_j \tag{3.6.9}$$

Доказательство. Равенство

$$\lambda_r {}^{\cdot k\,i}_{l\cdot j}(a) = A_r {}^{\cdot k\,i}_{l\cdot j}(a^{-1}, a) = a^{-1\cdot k}_{\quad\; j} \delta^i_l \tag{3.6.10}$$

следует из равенств (3.1.20), (3.6.2). Равенство (3.6.8) следует из равенства (3.6.10). Равенство

$$\lambda_l {}^{\cdot k\,i}_{l\cdot j}(b) = A_l {}^{\cdot k\,i}_{l\cdot j}(b, b^{-1}) = b^{-1\cdot i}_{\quad\; l} \delta^k_j \tag{3.6.11}$$

следует из равенств (3.1.21), (3.6.3). Равенство (3.6.9) следует из равенства (3.6.11). □

Теорема 3.6.4.

$$I^{\cdot k\,i\;\;m}_{l\cdot j\cdot n} = \delta^i_n \delta^m_l \delta^k_j \tag{3.6.12}$$

$$C^{\cdot k\,i\;\;m}_{l l\cdot j\cdot n} = I^{\cdot k\,m\,i}_{l\cdot n\cdot j} - I^{\cdot k\,i\;\;m}_{l\cdot j\cdot n} = \delta^m_j \delta^i_l \delta^k_n - \delta^i_n \delta^m_l \delta^k_j \tag{3.6.13}$$

Доказательство. Равенство (3.6.12) следует из равенств (3.5.14), (3.6.7) и равенства
$$I^{\cdot\,k\,i\,m}_{l\cdot j\,n} = \left.\frac{\partial \psi_l^{\cdot\,k\,i}_{\cdot l\cdot j}(a)}{\partial a^n_m}\right|_{a=e} = \left.\frac{\partial a^i_l \delta^k_j}{\partial a^n_m}\right|_{a=e} = \delta^i_n \delta^m_l \delta^k_j$$
Равенство (3.6.13) следует из равенств (3.5.17), (3.6.12). $\square$

Глава 4

# Представление группы Ли

### 4.1. Линейное представление группы Ли

Пусть $F$ - поле характеристики 0, $V$ - $F$-векторное пространство. Согласно теореме [1]-4.4.3, линейное преобразование $F$-векторного пространства имеет вид

$$v' = a\,v \tag{4.1.1}$$

где $a \in \mathcal{L}(F, V, V)$. Согласно теореме [1]-11.1.6, множество $GL(V)$ линейных преобразований (4.1.1) - это группа Ли, которую мы будем называть **группой линейных преобразований**.

Представление[4.1]

$$f : G \longrightarrow\!\!\!* V$$

группы Ли $G$ в непрерывном $F$-векторном пространстве $V$ называется **линейным представлением группы Ли** $G$ и имеет вид

$$v' = f(a) \circ v$$

Пусть $\overline{\overline{e}}$ - базис $F$-векторного пространства $V$. Тогда разложение вектора $v' = v'^{\alpha} e_{\alpha}$ имеет вид

$$v'^{\alpha} = f^{\alpha}_{\beta}(a) v^{\beta}$$

ТЕОРЕМА 4.1.1. *Пусть*

$$v' = f(a) \circ v \tag{4.1.2}$$

*линейное представление группы Ли $G$. Отображения $f^{\alpha}_{\beta}$ обладают следующими свойствами*

$$f^{\alpha}_{\beta}(e) = \delta^{\alpha}_{\beta} \tag{4.1.3}$$

$$f^{\alpha}_{\beta}(a^{-1}) = f^{-1\,\cdot\,\alpha}_{\phantom{-1}\beta}(a) \tag{4.1.4}$$

$$f^{\alpha}_{\beta}(ba) = f^{\alpha}_{\nu}(b) f^{\nu}_{\beta}(a) \tag{4.1.5}$$

ДОКАЗАТЕЛЬСТВО. Равенство (4.1.3) следует из утверждения, что $a = e$ соответствует тождественному отображению.

Так как отображение (4.1.2) является представлением группы Ли $G$, то равенство

$$v''^{\alpha} = f^{\alpha}_{\beta}(ba) v^{\beta} \tag{4.1.6}$$

следует из равенства (4.1.2) и равенства

$$v''^{\alpha} = f^{\alpha}_{\beta}(b) v'^{\beta} \tag{4.1.7}$$

---

[4.1]Мы рассмотрели общую концепцию представления группы в разделе [2]-2.4.

Равенство

$$v''^\alpha = f^\alpha_\nu(b) f^\nu_\beta(a) v^\beta \tag{4.1.8}$$

следует непосредственно из равенств (4.1.2), (4.1.7). Равенство (4.1.5) следует из равенств (4.1.6), (4.1.8). Из равенств (4.1.3), (4.1.5) следует

$$\delta^\alpha_\beta = f^\alpha_\nu(a^{-1}) f^\nu_\beta(a) \tag{4.1.9}$$

Равенство (4.1.4) следует из равенства (4.1.9). $\square$

Теорема 4.1.2. *Линейное представление группы Ли $G$*

$$v' = f(a) \circ v \tag{4.1.10}$$

*удовлетворяет системе дифференциальных уравнений*

$$\frac{\partial v'^\alpha}{\partial a^L} = I^\alpha_{\beta K} v'^\beta \lambda_l{}^K_L(a) \tag{4.1.11}$$

*где* **бесконечно малые образующие представления** *определены равенством*

$$I^\alpha_{\beta\,M} = \left.\frac{\partial f^\alpha_\beta(b)}{\partial b^M}\right|_{b=e} \qquad I = \left.\frac{\partial f(b)}{\partial b}\right|_{b=e} \tag{4.1.12}$$

Доказательство. Пусть $a, b \in G$. Так как отображение (4.1.10) является представлением группы Ли $G$, то равенство

$$v''^\alpha = f^\alpha_\beta(ba) v^\beta$$

следует из равенства (4.1.10) и равенства

$$v''^\alpha = f^\alpha_\beta(b) v'^\beta \tag{4.1.13}$$

Пусть

$$c = ba \tag{4.1.14}$$

Мы выразим $b$ через $a$ и $c$ и продифференцируем равенство (4.1.13) по $a$

$$\begin{aligned}
0 &= \frac{\partial f^\nu_\beta(b) v'^\beta}{\partial b^K} \frac{\partial b^K}{\partial a^L} + \frac{\partial f^\nu_\beta(b) v'^\beta}{\partial v'^\alpha} \frac{\partial v'^\alpha}{\partial a^L} \\
&= \frac{\partial f^\nu_\beta(b)}{\partial b^K} v'^\beta \frac{\partial b^K}{\partial a^L} + f^\nu_\alpha(b) \frac{\partial v'^\alpha}{\partial a^L}
\end{aligned} \tag{4.1.15}$$

Равенство

$$\frac{\partial v'^\alpha}{\partial a^L} = -f^{-1\cdot\alpha}_{\ \ \nu}(b) \frac{\partial f^\nu_\beta(b)}{\partial b^K} v'^\beta \frac{\partial b^K}{\partial a^L} \tag{4.1.16}$$

следует из равенства (4.1.15). Левая часть не зависит от $b$, поэтому и правая часть не зависит от $b$. Определим правую часть, когда $b = e$. Из равенств (4.1.14), (3.1.37), когда $ca^{-1} = e$, следует, что

$$\left.\frac{\partial b}{\partial a}\right|_{ca^{-1}=e} = \left.\frac{\partial ca^{-1}}{\partial a}\right|_{ca^{-1}=e} = -\left.\psi_r(ca^{-1})\right|_{ca^{-1}=e} \ \psi_l^{-1}(a) = -\lambda_l(a) \tag{4.1.17}$$

Равенство

$$\left. f^{-1\cdot\alpha}_{\ \ \nu}(b) \frac{\partial f^\nu_\beta(b)}{\partial b^K} v'^\beta \right|_{b=e} = I^\alpha_{\beta\,M} v'^\beta \tag{4.1.18}$$

следует из равенств (4.1.12). Равенство (4.1.11) следует из равенств (4.1.18), (4.1.17), (4.1.16). $\square$

Теорема 4.1.3. *Система дифференциальных уравнений* (4.1.11) *вполне интегрируема тогда и только тогда, когда*

$$(4.1.19) \qquad I^\alpha_{\beta K} I^\beta_{\nu P} - I^\alpha_{\beta P} I^\beta_{\nu K} = I^\alpha_{\nu T} C^T_{l\, PK}$$

Доказательство. Равенство

$$(4.1.20) \quad \begin{aligned} \frac{\partial^2 v'^\alpha}{\partial a^L \partial a^M} &= I^\alpha_{\beta K} \frac{\partial v'^\beta}{\partial a^M} \lambda_l{}^K_L(a) + I^\alpha_{\beta K} v'^\beta \frac{\partial \lambda_l{}^K_L(a)}{\partial a^M} \\ &= I^\alpha_{\beta K} I^\beta_{\nu P} v'^\nu \lambda_l{}^P_M(a) \lambda_l{}^K_L(a) + I^\alpha_{\beta K} v'^\beta \frac{\partial \lambda_l{}^K_L(a)}{\partial a^M} \end{aligned}$$

следует из равенства (4.1.11). Система дифференциальных уравнений (4.1.11) вполне интегрируема тогда и только тогда, когда

$$(4.1.21) \qquad \frac{\partial^2 v'^\alpha}{\partial a^L \partial a^M} = \frac{\partial^2 v'^\alpha}{\partial a^M \partial a^L}$$

Равенство

$$(4.1.22) \quad \begin{aligned} & I^\alpha_{\beta K} I^\beta_{\nu P} v'^\nu \lambda_l{}^P_M(a) \lambda_l{}^K_L(a) + I^\alpha_{\beta K} v'^\beta \frac{\partial \lambda_l{}^K_L(a)}{\partial a^M} \\ &= I^\alpha_{\beta P} I^\beta_{\nu K} v'^\nu \lambda_l{}^P_M(a) \lambda_l{}^K_L(a) + I^\alpha_{\beta K} v'^\beta \frac{\partial \lambda_l{}^K_M(a)}{\partial a^L} \end{aligned}$$

следует из равенств (4.1.20), (4.1.21). Равенство

$$(4.1.23) \quad \begin{aligned} (I^\alpha_{\beta K} I^\beta_{\nu P} - I^\alpha_{\beta P} I^\beta_{\nu K}) v'^\nu \lambda_l{}^P_M(a) \lambda_l{}^K_L(a) &= I^\alpha_{\beta K} v'^\beta \left( \frac{\partial \lambda_l{}^K_M(a)}{\partial a^L} - \frac{\partial \lambda_l{}^K_L(a)}{\partial a^M} \right) \\ &= I^\alpha_{\beta T} v'^\beta C^T_{l\, PK} \lambda_l{}^P_M(a) \lambda_l{}^K_L(a) \end{aligned}$$

следует из равенств (3.4.5), (4.1.22). Так как матрица $\lambda_l$ невырождена, то равенство

$$(4.1.24) \qquad (I^\alpha_{\beta K} I^\beta_{\nu P} - I^\alpha_{\beta P} I^\beta_{\nu K}) v'^\nu = I^\alpha_{\nu T} v'^\nu C^T_{l\, PK}$$

следует из равенства (4.1.23). Так как вектор $v$ произволен, то равенство (4.1.19) следует из равенства (4.1.24). □

Замечание 4.1.4. Пусть $F$ - поле и $V$ - $F$-векторное пространство размерности $n$. Линейное преобразование $a_1$ $F$-векторного пространства $V$

$$(4.1.25) \qquad v'^\alpha_1 = v^\beta_1 a_1{}^\alpha_\beta$$

записанное относительно базиса $\overline{\overline{e}}_{V1}$, отображает вектор

$$v = v^\alpha_1 e_{V1 \cdot \alpha}$$

в вектор

$$v' = v'^\alpha_1 e_{V1 \cdot \alpha}$$

Пассивное преобразование[4.2] $b$

$$(4.1.26) \qquad e_{V2 \cdot \alpha} = b^\beta_\alpha e_{V1 \cdot \beta}$$

отображает базис $\overline{\overline{e}}_{V1}$ в базис $\overline{\overline{e}}_{V2}$. Координаты векторов $v$, $v'$ относительно базиса $\overline{\overline{e}}_{V2}$

$$v = v^\alpha_2 e_{V2 \cdot \alpha}$$
$$v' = v'^\alpha_2 e_{V2 \cdot \alpha}$$

---

[4.2]Смотри определение пассивного преобразования в разделах [1]-5.2, [4]-5.2.

имеют вид

$$v_2^\alpha = v_1^\beta b^{-1\cdot\alpha}{}_\beta$$
$$v_2'^\alpha = v_2'^\beta b^{-1\cdot\alpha}{}_\beta$$

Следовательно, линейное преобразование (4.1.25) относительно базиса $\overline{\overline{e}}_{V2}$ имеет вид

(4.1.27) $$v_2'^\alpha = v_2^\beta a_{2\beta}{}^\alpha$$

(4.1.28) $$a_{2\beta}{}^\alpha = b^{-1\cdot\alpha}{}_\nu a_{1\lambda}{}^\nu b_\beta^\lambda$$

$\square$

Теорема 4.1.5. *Пусть*

(4.1.29) $$v_1'^\alpha = f_\beta^\alpha(a_1) v_1^\beta$$

*линейный геометрический объект.*[4.3] *Пассивное преобразование левого сдвига*

(4.1.30) $$a_1 \to g a_1$$

*на группе Ли $G$ порождает преобразование бесконечно малых образующих*[4.4]

(4.1.31) $$I_{2\cdot\mu K}{}^\alpha = f_\mu^\beta(g) f^{-1\cdot\alpha}{}_\rho(g) I_{1\cdot\beta J}{}^\rho \lambda_{lM}^J(g) \psi_{rK}^M(g)$$

Доказательство. Согласно равенствам (4.1.29), пассивное преобразование левого сдвига (4.1.30) на группе Ли $G$ порождает преобразование координат геометрического объекта

(4.1.32) $$v_2^\alpha = f_\beta^\alpha(g^{-1}) v_1^\beta$$

(4.1.33) $$v_2'^\alpha = f_\beta^\alpha(g^{-1}) v_1'^\beta$$

где мы полагаем

(4.1.34) $$v_2'^\alpha = f_\beta^\alpha(a_2) v_2^\beta$$

Равенства

(4.1.35) $$\frac{\partial v_1'^\alpha}{\partial a_1^L} = I_{1\cdot\beta K}{}^\alpha v_1'^\beta \lambda_{lL}^K(a_1)$$

(4.1.36) $$\frac{\partial v_2'^\alpha}{\partial a_2^L} = I_{2\cdot\beta K}{}^\alpha v_2'^\beta \lambda_{lL}^K(a_2)$$

следуют из равенства (4.1.11). Пользуясь правилом дифференцирования сложной функции, мы получим

(4.1.37) $$\frac{\partial v_2'^\pi}{\partial a_2^L} = \frac{\partial v_2'^\pi}{\partial v_1'^\rho} \frac{\partial v_1'^\rho}{\partial a_1^M} \frac{\partial a_1^M}{\partial a_2^L}$$

Равенство

(4.1.38) $$\frac{\partial v_2'^\pi}{\partial v_1'^\rho} = f_\rho^\pi(g^{-1})$$

---

[4.3]Смотри определение геометрического объекта в разделах [1]-5.3, [4]-5.3.

[4.4]Утверждение теоремы основано на утверждении, что базовый вектор представления группы Ли является ковариантным вектором в пространстве параметров ([5], с. 36) и контравариантным вектором в пространстве представления ([5], с. 35).

следует из равенства (4.1.33). Равенство

$$v_1'^\alpha = f^\alpha_\beta(ga_2g^{-1})v_1^\beta \tag{4.1.39}$$

следует из равенств (4.1.32), (4.1.33), (4.1.34). Из равенств (4.1.29), (4.1.39) следует, что пассивное преобразование левого сдвига (4.1.30) порождает преобразование в пространстве параметров

$$a_1 = ga_2g^{-1} \tag{4.1.40}$$

Равенство

$$\begin{aligned}\frac{\partial a_1^M}{\partial a_2^L} &= \psi_l{}^M_N(ga_2g^{-1})\lambda_l{}^N_K(ga_2)\psi_r{}^K_I(ga_2)\lambda_r{}^I_L(a_2)\\ &= \psi_l{}^M_N(a_1)\lambda_l{}^N_K(ga_2)\psi_r{}^K_I(ga_2)\lambda_r{}^I_L(a_2)\end{aligned} \tag{4.1.41}$$

следует из равенств (3.1.45), (4.1.40). Равенство

$$\begin{aligned}&I_{2\cdot}{}^\alpha_{\beta K}v_2'^\beta\lambda_l{}^K_L(a_2)\\ &= f^\alpha_\rho(g^{-1})I_{1\cdot}{}^\rho_{\beta J}v_1'^\beta\lambda_l{}^J_M(a_1)\psi_l{}^M_N(a_1)\lambda_l{}^N_K(ga_2)\psi_r{}^K_I(ga_2)\lambda_r{}^I_L(a_2)\end{aligned} \tag{4.1.42}$$

следует из равенств (4.1.35), (4.1.36), (4.1.37), (4.1.38), (4.1.41). Равенство

$$\begin{aligned}&I_{2\cdot}{}^\alpha_{\beta K}f^\beta_\mu(g^{-1})v_1'^\mu\\ &= f^\alpha_\rho(g^{-1})I_{1\cdot}{}^{rho}_{\beta J}v_1'^\beta\lambda_l{}^J_M(ga_2)\psi_r{}^M_I(ga_2)\lambda_r{}^I_L(a_2)\psi_l{}^L_K(a_2)\end{aligned} \tag{4.1.43}$$

следует из равенств (3.1.19), (4.1.33), (4.1.42). Левая часть равенства (4.1.43) не зависит от $a_2$, поэтому и правая часть не зависит от $a_2$. Равенство

$$\begin{aligned}&I_{2\cdot}{}^\alpha_{\beta K}f^\beta_\mu(g^{-1})v_1'^\mu\\ &= f^\alpha_\rho(g^{-1})I_{1\cdot}{}^\rho_{\beta J}v_1'^\beta\lambda_l{}^J_M(g)\psi_r{}^M_K(g)\end{aligned} \tag{4.1.44}$$

следует из равенства (4.1.43), когда $a_2 = e$. Так как $v_1'$ - произвольный вектор, то равенства

$$I_{2\cdot}{}^\alpha_{\mu K}f^{-1\cdot\mu}_\beta(g) = f^{-1\cdot\alpha}_\rho(g)I_{1\cdot}{}^\rho_{\beta J}\lambda_l{}^J_M(g)\psi_r{}^M_K(g) \tag{4.1.45}$$

следуют из равенств (4.1.4), (4.1.44). Равенство (4.1.31) следует из равенства (4.1.45). □

## 4.2. Сопряжённое представление

ОПРЕДЕЛЕНИЕ 4.2.1. Пусть $V$ - векторное пространство над полем $F$. Векторное пространство $V' = \mathcal{L}(F; V; F)$ называется **сопряжённым векторным пространством**. □

Согласно определению 4.2.1, элементы сопряжённого векторного пространства являются линейными отображениями

$$u : V \to F \tag{4.2.1}$$

Линейное отображение (4.2.1) называется **линейным функционалом** на векторном пространстве $V$. Мы будем записывать значение линейного функционала в виде

$$u \circ a = (u, a)$$

Теорема 4.2.2. *Пусть $\overline{\overline{e}} = (e_\alpha \in V, \alpha \in I)$ - базис $F$-векторного пространства $V$. Пусть $V'$ - $F$-векторное пространство, сопряжённое $F$-векторному пространству $V$. Множество векторов $\overline{\overline{e}} = (e^\alpha \in V, \alpha \in I)$ такое, что*

$$(4.2.2) \qquad (e^\alpha, e_\beta) = e^\alpha \circ e_\beta = \delta^\alpha_\beta$$

*является базисом $F$-векторного пространства $V'$.*

Доказательство. Из равенства (4.2.2) следует, что

$$(4.2.3) \qquad (e^\alpha, a) = (e^\alpha, e_\beta a^\beta) = (e^\alpha, e_\beta) a^\beta = \delta^\alpha_\beta a^\beta = a^\alpha$$

Пусть

$$u : V \to F$$

линейный функционал. Тогда

$$(4.2.4) \qquad (u, a) = (u, e_\alpha a^\alpha) = (u, e_\alpha) a^\alpha = u_\alpha a^\alpha$$

где

$$(u, e_\alpha) = u_\alpha$$

Из равенств (4.2.3), (4.2.4) следует, что

$$(u, a) = u_\alpha (e^\alpha, a)$$

Согласно определениям [3]-(2.1.9), [3]-(2.1.10)

$$u = u_\alpha e^\alpha$$

Следовательно, множество $\overline{\overline{e}} = (e^\alpha \in V, \alpha \in I)$ является базисом $F$-векторного пространства $V'$. □

Следствие 4.2.3. *Пусть $\overline{\overline{e}} = (e_\alpha \in V, \alpha \in I)$ - базис $F$-векторного пространства $V$. Тогда*

$$a^\alpha = (e^\alpha, a)$$

□

Базис $\overline{\overline{e}} = (e^\alpha \in V, \alpha \in I)$ называется **базисом, двойственным базису** $\overline{\overline{e}} = (e_\alpha \in V, \alpha \in I)$.

Поскольку произведение в поле коммутативно, то не существует различия между левым и правым векторным пространствами. Если мы рассматриваем векторное пространство над телом, то мы различаем левое и правое векторное пространство. Я планирую в будущем исследовать группы преобразований в векторном пространстве над телом. Поэтому в этой книге я хочу пользоваться записью, которая будет совместима с некоммутативным произведением.

Согласно теореме [1]-4.4.3, если мы рассматриваем правое векторное пространство, то мы пишем матрицу отображения слева от координат вектора. В частности, это верно для матрицы линейного преобразования и для координат линейного функционала. Согласно теореме [1]-6.2.2, множество линейных функционалов на правом векторном пространстве порождает левое векторное пространство. Согласно теореме [1]-4.4.3, мы пишем матрицу отображения справа от координат линейного функционала. Чтобы сохранить функциональную запись, мы будем записывать линейное отображение сопряжённого векторного пространства в виде

$$[A] \circ u = u_\alpha A^\alpha_\beta e^\beta$$

Определение 4.2.4. Пусть $\overline{\overline{e}} = (e_\alpha \in V, \alpha \in I)$ - базис $F$-векторного пространства $V$. Пусть $\overline{\overline{e}} = (e^\alpha \in V, \alpha \in I)$ - базис, двойственный базису $\overline{\overline{e}} = (e_\alpha \in V, \alpha \in I)$. Пусть

$$v_1'^\alpha = f_{1\cdot\beta}^{\ \alpha}(a) v_1^\beta \tag{4.2.5}$$

представление группы Ли $G$ в непрерывном $F$-векторном пространстве $V$. Представление

$$v'_{2\cdot\alpha} = v_{2\cdot\beta} f_{2\cdot\alpha}^{\ \beta}(a)$$

группы Ли $G$ в непрерывном $F$-векторном пространстве $V'$ называется **представлением, сопряжённым представлению** (4.2.5), если

$$([f_2(a)] \circ v_2, f_1(a) \circ v_1) = (u, v) \tag{4.2.6}$$

для любого $a \in G$. $\square$

Теорема 4.2.5. *Пусть $\overline{\overline{e}} = (e_\alpha \in V, \alpha \in I)$ - базис $F$-векторного пространства $V$. Пусть $\overline{\overline{e}} = (e^\alpha \in V, \alpha \in I)$ - базис, двойственный базису $\overline{\overline{e}} = (e_\alpha \in V, \alpha \in I)$. Если представление*

$$v_1'^\alpha = f_{1\cdot\beta}^{\ \alpha}(a) v_1^\beta \tag{4.2.7}$$

*группы Ли $G$ в непрерывном $F$-векторном пространстве $V$ имеет сопряжённое представление*

$$v'_{2\cdot\alpha} = v_{2\cdot\beta} f_{2\cdot\alpha}^{\ \beta}(a) \tag{4.2.8}$$

*группы Ли $G$ в непрерывном $F$-векторном пространстве $V'$, то*

$$f_2(a) = f_1^{-1}(a) = f_1(a^{-1}) \tag{4.2.9}$$

Доказательство. Равенство

$$f_1(a) \circ v_1 = e_\nu f_{1\cdot\beta}^{\ \nu}(a) v_1^\beta \tag{4.2.10}$$

следует из равенства (4.2.7). Равенство

$$[f_2(a)] \circ v_2 = v_{2\cdot\alpha} f_{2\cdot\lambda}^{\ \alpha}(a) e^\lambda \tag{4.2.11}$$

следует из равенства (4.2.8). Равенство

$$\begin{aligned}
v_2 \circ v_1 &= v_{2\cdot\alpha} v_1^\alpha = ([f_2(a)] \circ v_2) \circ (f_1(a) \circ v_1) \\
&= (v_{2\cdot\alpha} f_{2\cdot\mu}^{\ \alpha}(a) e^\mu) \circ (e_\nu f_{1\cdot\beta}^{\ \nu}(a) v_1^\beta) \\
&= v_{2\cdot\alpha} f_{2\cdot\mu}^{\ \alpha}(a)(e^\mu \circ e_\nu) f_{1\cdot\beta}^{\ \nu}(a) v_1^\beta \\
&= v_{2\cdot\alpha} f_{2\cdot\mu}^{\ \alpha}(a) \delta^\mu_\nu f_{1\cdot\beta}^{\ \nu}(a) v_1^\beta \\
&= v_{2\cdot\alpha} f_{2\cdot\nu}^{\ \alpha}(a) f_{1\cdot\beta}^{\ \nu}(a) v_1^\beta
\end{aligned} \tag{4.2.12}$$

следует из равенств (4.2.6), (4.2.10), (4.2.11). Так как $v$ - произвольный вектор, то уравнение

$$v_{2\cdot\beta} = v_{2\cdot\alpha} f_{2\cdot\nu}^{\ \alpha}(a) f_{1\cdot\beta}^{\ \nu}(a) \tag{4.2.13}$$

следует из равенства (4.2.12). Так как $v_2$ - произвольный функционал, то уравнение

$$\delta^\alpha_\beta = f_{2\cdot\nu}^{\ \alpha}(a) f_{1\cdot\beta}^{\ \nu}(a) \tag{4.2.14}$$

следует из равенства (4.2.13). Равенство (4.2.9) следует из равенства (4.2.14).

Так как отображение

$$a \to f_1(a) \tag{4.2.15}$$

является представлением группы $G$, то равенство

$$\begin{aligned}(4.2.16)\quad v_{2\cdot\beta}f_{2\cdot\alpha}{}^{\beta}(a_1)f_{2\cdot\nu}{}^{\alpha}(a_2) &= v_{2\cdot\beta}f_{1\cdot\alpha}{}^{\beta}(a_1^{-1})f_{1\cdot\nu}{}^{\alpha}(a_2^{-1}) \\ &= v_{2\cdot\beta}f_{1\cdot\nu}{}^{\alpha}(a_2^{-1})f_{1\cdot\alpha}{}^{\beta}(a_1^{-1}) \\ &= v_{2\cdot\beta}f_{1\cdot\nu}{}^{\beta}(a_2^{-1}a_1^{-1}) = v_{2\cdot\beta}f_{1\cdot\nu}{}^{\beta}((a_1a_2)^{-1}) \\ &= v_{2\cdot\beta}f_{2\cdot\nu}{}^{\beta}(a_1a_2)\end{aligned}$$

следует из равенств (4.1.5), (4.2.9). Из равенства (4.2.16) следует, что отображение

$$a \to f_2(a)$$

является правосторонним представлением группы $G$. $\square$

Теорема 4.2.6. *Пусть*

$$v_1'^{\alpha} = f_{1\cdot\beta}{}^{\alpha}(a)v_1^{\beta} \tag{4.2.17}$$

*линейное представление группы $G$ в $F$-векторном пространстве $V$. Представление*

$$v_{2\cdot\alpha}' = v_{2\cdot\beta}f_{2\cdot\alpha}{}^{\beta}(a) \tag{4.2.18}$$

*сопряжённое представлению (4.2.17), удовлетворяет системе дифференциальных уравнений*

$$\frac{\partial v_{2\cdot\beta}'}{\partial a^L} = v_{2\cdot\alpha}'I_{2\cdot\beta K}{}^{\alpha}\lambda_l{}_L^K(a) \tag{4.2.19}$$

*где бесконечно малые образующие представления (4.2.18) определены равенством*

$$I_{2\cdot\beta L}{}^{\alpha} = -I_{1\cdot\beta L}{}^{\alpha} \tag{4.2.20}$$

Доказательство. Равенство

$$I_{2\cdot\beta M}{}^{i} = \left.\frac{\partial f_{2\cdot\beta}{}^{\alpha}(a)}{\partial a^M}\right|_{a=e} = \left.\frac{\partial f_{1\cdot\beta}{}^{\alpha}(a^{-1})}{\partial a^M}\right|_{a=e} = \left.\frac{\partial f_{1\cdot\beta}{}^{\alpha}(a^{-1})}{\partial a^{-1\cdot N}}\frac{\partial a^{-1\cdot N}}{\partial a^M}\right|_{a=e} \tag{4.2.21}$$

следует из равенств (4.1.12), (4.2.9). Очевидно, что

$$\left.\frac{\partial f_{1\cdot\beta}{}^{\alpha}(a^{-1})}{\partial a^{-1\cdot N}}\right|_{a=e} = \left.\frac{\partial f_{1\cdot\beta}{}^{\alpha}(a)}{\partial a^N}\right|_{a=e} = I_{1\cdot\beta N}{}^{\alpha} \tag{4.2.22}$$

Из равенства (3.1.33) следует

$$\left.\frac{\partial a^{-1\cdot N}}{\partial a^M}\right|_{a=e} = -\psi_{l\cdot K}{}^{N}(a^{-1})\lambda_{r\cdot M}{}^{K}(a)\big|_{a=e} = -\delta_M^N \tag{4.2.23}$$

Равенство (4.2.20) следует из равенств (4.2.21), (4.2.22), (4.2.23).

Пусть $v_1$ - геометрический объект представления (4.2.17). Пусть $v_2$ - геометрический объект представления (4.2.18). Рассмотрим выражение

$$v = v_{2\cdot\alpha}v_1^{\alpha} \tag{4.2.24}$$

Равенство

$$v' = v_{2\cdot\alpha}'v_1'^{\alpha} = v_{2\cdot\beta}f_{2\cdot\alpha}{}^{\beta}(a)f_{1\cdot\nu}{}^{\alpha}(a)v_1^{\nu} = v_{2\cdot\beta}v_1^{\beta} = v \tag{4.2.25}$$

является следствием равенств (4.2.14), (4.2.17), (4.2.18), (4.2.24). Равенство

$$(4.2.26) \qquad \frac{\partial v'}{\partial a^L} = 0$$

следует из равенства (4.2.25). Равенство

$$(4.2.27) \qquad \frac{\partial v'}{\partial a^L} = \frac{\partial v'_{2 \cdot \alpha} v'^\alpha_1}{\partial a^L} = \frac{\partial v'_{2 \cdot \alpha}}{\partial a^L} v'^\alpha_1 + v'_{2 \cdot \alpha} \frac{\partial v'^\alpha_1}{\partial a^L}$$

следует из равенства (4.2.24). Равенство

$$(4.2.28) \qquad \frac{\partial v'_{2 \cdot \alpha}}{\partial a^L} v'^\alpha_1 + v'_{2 \cdot \alpha} \frac{\partial v'^\alpha_1}{\partial a^L} = 0$$

следует из равенств (4.2.26), (4.2.27). Равенство

$$(4.2.29) \qquad \frac{\partial v'_{2 \cdot \alpha}}{\partial a^L} v'^\alpha_1 + v'_{2 \cdot \alpha} I_{1 \cdot \beta K}^{\;\;\alpha} v'^\beta_1 \lambda_l{}^K_L(a) = 0$$

следует из равенств (4.1.11), (4.2.28). Так как вектор $v_1$ - произволен, то равенство

$$(4.2.30) \qquad \frac{\partial v'_{2 \cdot \beta}}{\partial a^L} = -v'_{2 \cdot \alpha} I_{1 \cdot \beta K}^{\;\;\alpha} \lambda_l{}^K_L(a)$$

следует из равенства (4.2.29). Равенство (4.2.19) следует из равенств (4.2.20), (4.2.30). $\square$

## 4.3. Алгебраические свойства линейного представления

Теорема 4.3.1. *Пусть*

$$(4.3.1) \qquad \begin{aligned} v'^\alpha_1 &= f^\alpha_\beta(a) v^j_1 \\ v'^\alpha_2 &= f^\alpha_\beta(a) v^j_2 \end{aligned}$$

*линейные геометрические объекты. Тогда их сумма также является линейным геометрическим объектом*

$$(4.3.2) \qquad v'^\alpha = f^\alpha_\beta(a) v^\beta$$

Доказательство. Из равенства (4.3.1) следует, что

$$(4.3.3) \qquad v'\alpha = v'^\alpha_1 + v'^\alpha_2 = f^\alpha_\beta(a) v^\beta_1 + f^\alpha_\beta(a) v^\beta_2$$

Равенство (4.3.2) следует из равенства (4.3.3). $\square$

Теорема 4.3.2. *Пусть*

$$v'^\alpha_1 = f_{1 \cdot \beta}^{\;\;\alpha}(a) v^\beta_1$$

*представление группы $G$ в $F$-векторном пространстве $V_1$. Пусть*

$$v'^\gamma_2 = f_{2 \cdot \epsilon}^{\;\;\gamma}(a) v^\epsilon_2$$

*представление группы $G$ в $F$-векторном пространстве $V_2$. Пусть геометрический объект*

$$v'^{\alpha\gamma} = f_{1 \cdot \beta}^{\;\;\alpha}(a) f_{2 \cdot \epsilon}^{\;\;\gamma}(a) v^{\beta\epsilon}$$

*является тензорным произведением геометрических объектов $v_1 \in V_1$ и $v_2 \in V_2$*

$$(4.3.4) \qquad f(a)(v_1 \otimes v_2) = f_1(a)(v_1) \otimes f_2(a)(v_2)$$

Мы будем называть представление $f$ **тензорным произведением представлений** и

(4.3.5) $$I = I_1 \otimes E_2 + E_1 \otimes I_2$$

Доказательство. Прежде всего покажем, что мы получили представление группы $G$ в векторном пространстве $V = V_1 \otimes V_2$. Если мы применим преобразования для $a \in G$ и $b \in G$ последовательно, то мы имеем

(4.3.6)
$$\begin{aligned}
f(a)(f(b)(v_1 \otimes v_2)) &= f(a)(f_1(b)(v_1) \otimes f_2(b)(v_2)) \\
&= f_1(a)(f_1(b)(v_1)) \otimes f_2(a)(f_2(b)(v_2)) \\
&= f_1(ab)(v_1) \otimes f_2(ab)(v_2) \\
&= f(ab)(v_1 \otimes v_2)
\end{aligned}$$

Следовательно, $f$ тоже представление.

Дифференцируя равенство

$$f(a) = f_1(a) \otimes f_2(a)$$

по $a$ мы получим

(4.3.7) $$\frac{\partial f}{\partial a} = \frac{\partial f_1}{\partial a} \otimes f_2 + f_1 \otimes \frac{\partial f_2}{\partial a}$$

Мы получим (4.3.5) из (4.3.7) по определению, когда $a = e$. □

Теорема 4.3.3. *Пусть векторное пространство $V$ является прямой суммой векторном пространств $V_1$ и $V_2$*

$$V = V_1 \oplus V_2$$

*Пусть*

$$v_1'^\alpha = f_{1 \cdot \beta}{}^\alpha(a) v_1^\beta$$

*представление группы $G$ в векторном пространстве $V_1$. Пусть*

$$v_2'^\gamma = f_{2 \cdot \epsilon}{}^\gamma(a) v_2^\epsilon$$

*представление группы $G$ в векторном пространстве $V_2$. Пусть геометрический объект $v$ является прямой суммой геометрических объектов $v_1 \in V_1$ и $v_2 \in V_2$*

(4.3.8) $$f(a) \circ (v_1 \oplus v_2) = (f_1(a) \circ v_1) \oplus (f_2(a) \circ v_2)$$

*Тогда определено представление группы $G$ в векторном пространстве $V = V_1 \oplus V_2$ и*

$$I_P = \begin{pmatrix} I_1{}^{a_1}_{b_1 P} & 0 \\ 0 & I_2{}^{a_2}_{b_2 P} \end{pmatrix}$$

Мы будем называть представление $f$ **прямой суммой представлений**.

Доказательство. Мы просто дифференцируем отображение (4.3.8). □

### 4.4. Дифференциальные уравнения представления

Теорема 4.4.1. *Линейное представление группы Ли $G$*

$$v'^\alpha = f^\alpha_\beta(a) v^\beta \tag{4.4.1}$$

*удовлетворяет системе дифференциальных уравнений*

$$\frac{\partial f^\alpha_\beta}{\partial a^L} = I^\alpha_{\nu K} f^\nu_\beta \lambda_l{}^K_L(a) \tag{4.4.2}$$

Доказательство. Согласно теореме 4.1.2, линейное представление (4.4.1) удовлетворяет системе дифференциальных уравнений

$$\frac{\partial v'^\alpha}{\partial a^L} = I^\alpha_{\beta K} v'^\beta \lambda_l{}^K_L(a) \tag{4.4.3}$$

Уравнение

$$\frac{\partial f^\alpha_\beta(a) v^\beta}{\partial a^L} = I^\alpha_{\beta K} f^\beta_\nu(a) v^\nu \lambda_l{}^K_L(a) \tag{4.4.4}$$

является следствием равенства (4.4.3). Так как $v$ - произвольный вектор, то равенства (4.4.2), являются следствием равенства (4.4.4). $\square$

Теорема 4.4.2. *Пусть*

$$v'^\alpha_1 = f_{1\cdot}{}^\alpha_\beta(a) v^\beta_1 \tag{4.4.5}$$

*линейное представление группы $G$ в $F$-векторном пространстве $V$. Тогда*

$$I_{1\cdot}{}^\alpha_{\beta K} f_{1\cdot}{}^\beta_\gamma \lambda_l{}^K_L(a) = f_{1\cdot}{}^\alpha_\rho(a) I_{1\cdot}{}^\rho_{\gamma J} \lambda_r{}^J_L(a) \tag{4.4.6}$$

Доказательство. Согласно теореме 4.2.6, существует линейное представление

$$v'_{2\cdot\alpha} = v_{2\cdot\beta} f_{2\cdot}{}^\beta_\alpha(a) \tag{4.4.7}$$

сопряжённое представлению (4.4.5) и сопряжённое представление удовлетворяет системе дифференциальных уравнений

$$\frac{\partial v'_{2\cdot\beta}}{\partial a^L} = v'_{2\cdot\alpha} I_{2\cdot}{}^\alpha_{\beta K} \lambda_l{}^K_L(a) \tag{4.4.8}$$

Уравнение

$$\frac{\partial v_{2\cdot\gamma} f_{2\cdot}{}^\gamma_\beta(a)}{\partial a^L} = v_{2\cdot\gamma} f_{2\cdot}{}^\gamma_\alpha(a) I_{2\cdot}{}^\alpha_{\beta K} \lambda_l{}^K_L(a) \tag{4.4.9}$$

является следствием равенства (4.4.8). Так как $v_2$ - произвольный вектор, то равенство

$$\frac{\partial f_{2\cdot}{}^\gamma_\beta(a)}{\partial a^L} = f_{2\cdot}{}^\gamma_\alpha(a) I_{2\cdot}{}^\alpha_{\beta K} \lambda_l{}^K_L(a) \tag{4.4.10}$$

является следствием равенства (4.4.9). Равенство

$$\frac{\partial f_{1\cdot}{}^\gamma_\beta(a^{-1})}{\partial a^L} = -f_{1\cdot}{}^\gamma_\alpha(a^{-1}) I_{1\cdot}{}^\alpha_{\beta K} \lambda_l{}^K_L(a) \tag{4.4.11}$$

является следствием равенств (4.2.9), (4.2.20), (4.4.10). Равенство

$$\frac{\partial f_{1\cdot}{}^\gamma_\beta(a^{-1})}{\partial a^L} = \frac{\partial f_{1\cdot}{}^\gamma_\beta(a^{-1})}{\partial a^{-1\cdot M}} \frac{\partial a^{-1\cdot M}}{\partial a^L} = -\frac{\partial f_{1\cdot}{}^\gamma_\beta(a^{-1})}{\partial a^{-1\cdot M}} \psi_r{}^M_N(a^{-1}) \lambda_l{}^N_L(a) \tag{4.4.12}$$

является следствием равенства (3.1.33). Равенство

$$\tag{4.4.13} \begin{aligned} -\frac{\partial f_{1\cdot\beta}^{\ \gamma}(a^{-1})}{\partial a^{-1\cdot M}} &= -f_{1\cdot\alpha}^{\ \gamma}(a^{-1})I_{1\cdot\beta K}^{\alpha}\lambda_{l\,L}^{\ K}(a)\psi_{l\,N}^{\ L}(a)\lambda_{r\,M}^{\ N}(a^{-1}) \\ &= -f_{1\cdot\alpha}^{\ \gamma}(a^{-1})I_{1\cdot\beta K}^{\alpha}\lambda_{r\,M}^{\ K}(a^{-1}) \end{aligned}$$

следует из равенств (3.1.18), (3.1.19), (4.4.11), (4.4.12). Равенство

$$\tag{4.4.14} \frac{\partial f_{1\cdot\beta}^{\ \gamma}(a)}{\partial a^M} = f_{1\cdot\alpha}^{\ \gamma}(a)I_{1\cdot\beta K}^{\alpha}\lambda_{r\,M}^{\ K}(a)$$

следует из равенства (4.4.13). Согласно теореме 4.4.1, линейное представление (4.4.5) удовлетворяет системе дифференциальных уравнений

$$\tag{4.4.15} \frac{\partial f_{1\cdot\beta}^{\ \alpha}}{\partial a^L} = I_{\nu K}^{\alpha} f_{1\cdot\beta}^{\ \nu}\lambda_{l\,L}^{\ K}(a)$$

Равенство (4.4.6) является следствием равенств (4.4.13), (4.4.15). □

Теорема 4.4.3. *Бесконечно малые образующие линейного представления являются константами.*

Доказательство. Согласно теореме 4.1.5, пассивное преобразование левого сдвига

$$a_1 \to ga_1$$

на группе Ли $G$ порождает преобразование бесконечно малых образующих

$$\tag{4.4.16} I_{2\cdot\mu K}^{\alpha} = f_{\mu}^{\beta}(g)f^{-1\cdot\alpha}_{\ \rho}(g)I_{1\cdot\beta J}^{\rho}\lambda_{l\,M}^{\ J}(g)\psi_{r\,K}^{\ M}(g)$$

Равенство

$$\tag{4.4.17} I_{1\cdot\gamma J}^{\rho} = f^{-1\cdot\rho}_{\ \alpha}(a)I_{1\cdot\beta K}^{\alpha}f_{\gamma}^{\beta}\lambda_{l\,L}^{\ K}(a)\psi_{r\,J}^{\ L}(a)$$

следует из равенств (3.1.18), (3.1.19), (4.4.6). Теорема следует из равенств (4.4.16), (4.4.17). □

Теорема 4.4.4. *Система дифференциальных уравнений (4.4.2) вполне интегрируема тогда и только тогда, когда*

$$\tag{4.4.18} I_{\beta K}^{\alpha}I_{\nu P}^{\beta} - I_{\beta P}^{\alpha}I_{\nu K}^{\beta} = I_{\nu T}^{\alpha}C_{l\,PK}^{\ T}$$

Доказательство. Равенство

$$\tag{4.4.19} \begin{aligned} \frac{\partial^2 f_{\beta}^{\alpha}}{\partial a^L \partial a^M} &= I_{\nu K}^{\alpha}\frac{\partial f_{\beta}^{\nu}}{\partial a^M}\lambda_{l\,L}^{\ K}(a) + I_{\nu K}^{\alpha}f_{\beta}^{\nu}\frac{\partial \lambda_{l\,L}^{\ K}(a)}{\partial a^M} \\ &= I_{\nu K}^{\alpha}I_{\mu P}^{\nu}f_{\beta}^{\mu}\lambda_{l\,M}^{\ P}(a)\lambda_{l\,L}^{\ K}(a) + I_{\nu K}^{\alpha}f_{\beta}^{\nu}\frac{\partial \lambda_{l\,L}^{\ K}(a)}{\partial a^M} \end{aligned}$$

следует из равенства (4.4.2). Система дифференциальных уравнений (4.4.2) вполне интегрируема тогда и только тогда, когда

$$\tag{4.4.20} \frac{\partial^2 f_{\beta}^{\alpha}}{\partial a^L \partial a^M} = \frac{\partial^2 f_{\beta}^{\alpha}}{\partial a^M \partial a^L}$$

Равенство

$$\tag{4.4.21} \begin{aligned} &I_{\nu K}^{\alpha}I_{\mu P}^{\nu}f_{\beta}^{\mu}\lambda_{l\,M}^{\ P}(a)\lambda_{l\,L}^{\ K}(a) + I_{\nu K}^{\alpha}f_{\beta}^{\nu}\frac{\partial \lambda_{l\,L}^{\ K}(a)}{\partial a^M} \\ &= I_{\nu P}^{\alpha}I_{\mu K}^{\nu}f_{\beta}^{\mu}\lambda_{l\,M}^{\ P}(a)\lambda_{l\,L}^{\ K}(a) + I_{\nu K}^{\alpha}f_{\beta}^{\nu}\frac{\partial \lambda_{l\,M}^{\ K}(a)}{\partial a^L} \end{aligned}$$

следует из равенств (4.4.19), (4.4.20). Равенство

$$
\begin{aligned}
(I^\alpha_{\nu K} I^\nu_{\mu P} - I^\alpha_{\nu P} I^\nu_{\mu K}) f^\mu_\beta \lambda^P_{l\,M}(a) \lambda^K_{l\,L}(a) \\
(4.4.22) \quad = I^\alpha_{\nu K} f^\nu_\beta \left( \frac{\partial \lambda^K_{l\,M}(a)}{\partial a^L} - \frac{\partial \lambda^K_{l\,L}(a)}{\partial a^M} \right) \\
= I^\alpha_{\nu K} f^\nu_\beta C^K_{l\,PK} \lambda^P_{l\,M}(a) \lambda^K_{l\,L}(a)
\end{aligned}
$$

следует из равенств (3.4.3), (4.4.21). Так как матрица $\lambda_l$ невырожденна и отображение $f$ произвольно, то равенство (4.4.18) следует из равенства (4.1.23). $\square$

### 4.5. Группа $Gl_n(F)$

Рассмотрим представление группы $GL_n(F)$ в сопряжённом векторном пространстве

$$(4.5.1) \qquad u'_k = u_l a^{-1\cdot l}{}_k$$

Чтобы записать систему дифференциальных уравнений для представления (4.5.1), мы повторим вывод системы дифференциальных уравнений (4.1.11) учитывая представление (4.5.1).

Уравнение (4.1.15) имеет вид

$$
\begin{aligned}
(4.5.2) \quad 0 &= \frac{\partial u'_b b^{-1\cdot b}{}_n}{\partial b^k_l} \frac{\partial b^k_l}{\partial a^m_q} + \frac{\partial u'_b b^{-1\cdot b}{}_n}{\partial u'_c} \frac{\partial u'_c}{\partial a^m_q} \\
&= u'_b \frac{\partial b^{-1\cdot b}{}_n}{\partial b^k_l} \frac{\partial b^k_l}{\partial a^m_q} + \frac{\partial u'_b}{\partial u'_c} b^{-1\cdot b}{}_n \frac{\partial u'_c}{\partial a^m_q}
\end{aligned}
$$

Равенство

$$
\begin{aligned}
(4.5.3) \quad 0 &= \frac{\partial b^k_l b^{-1\cdot l}{}_m}{\partial b^p_r} = \frac{\partial b^k_l}{\partial b^p_r} b^{-1\cdot l}{}_m + b^k_l \frac{\partial b^{-1\cdot l}{}_m}{\partial b^p_r} = \delta^k_p \delta^r_l b^{-1\cdot l}{}_m + b^k_l \frac{\partial b^{-1\cdot l}{}_m}{\partial b^p_r} \\
&= \delta^k_p b^{-1\cdot r}{}_m + b^k_l \frac{\partial b^{-1\cdot l}{}_m}{\partial b^p_r}
\end{aligned}
$$

является следствием равенства $b^k_l b^{-1\cdot l}{}_m = \delta^k_m$. Равенство

$$(4.5.4) \qquad \frac{\partial b^{-1\cdot l}{}_m}{\partial b^p_r} = -b^{-1\cdot l}{}_k \delta^k_p b^{-1\cdot r}{}_m = -b^{-1\cdot l}{}_p b^{-1\cdot r}{}_m$$

является следствием равенства (4.5.3). Равенство

$$(4.5.5) \qquad \frac{\partial b^k_l}{\partial a^m_n} = \frac{\partial c^k_p a^{-1\cdot p}{}_l}{\partial a^m_n} = c^k_p \frac{\partial a^{-1\cdot p}{}_l}{\partial a^m_n} = -c^k_p a^{-1\cdot p}{}_m a^{-1\cdot n}{}_l$$

является следствием равенства $b = ca^{-1}$ и равенства (4.5.3).

Равенство

$$
\begin{aligned}
(4.5.6) \quad 0 &= u'_b b^{-1\cdot b}{}_k b^{-1\cdot l}{}_n c^k_p a^{-1\cdot p}{}_m a^{-1\cdot q}{}_l + \delta^c_b b^{-1\cdot b}{}_n \frac{\partial u'_c}{\partial a^m_q} \\
&= u'_b b^{-1\cdot b}{}_k b^{-1\cdot l}{}_n b^k_m a^{-1\cdot q}{}_l + b^{-1\cdot c}{}_n \frac{\partial u'_c}{\partial a^m_q} \\
&= u'_m b^{-1\cdot l}{}_n a^{-1\cdot q}{}_l + b^{-1\cdot c}{}_n \frac{\partial u'_c}{\partial a^m_q}
\end{aligned}
$$

является следствием равенств (4.5.2), (4.5.4), (4.5.5). Равенство

$$(4.5.7) \qquad \frac{\partial u'_c}{\partial a^m_q} = -u'_m a^{-1\cdot q}{}_c$$

является следствием равенства (4.5.6).

Равенство

$$(4.5.8) \qquad I_{\cdot-m\cdot p}^{l\cdot-r} = -\delta^l_p \delta^r_m$$

является следствием равенства (4.5.4).

Глава 5

# Список литературы

# Глава 6

# Предметный указатель

# Глава 7

# Специальные символы и обозначения